# Dissipative Control of General Linear Time-Delay Systems: Applications of the Kronecker-Seuret Decomposition*

Qian Feng[a,*], Wei Xing Zheng[b], Xiaoyu Wang[c], Feng Xiao[a,d,*]

[a] School of Control and Computer Engineering, North China Electric Power University, Beijing, China.
[b] School of Computer, Data and Mathematical Sciences, Western Sydney University, Sydney, Australia.
[c] Institute of Systems Science, Beijing Wuzi University, Beijing, China.
[d] State Key Laboratory of Alternate Electrical Power System with Renewable Energy Sources, North China Electric Power University, Beijing, China.

**Abstract**

Stabilizing linear time-invariant delay systems, particularly when addressing an arbitrary finite number of pointwise and distributed delays (DDs) under dissipative constraints, poses a significant challenge. Existing solutions are often hindered by theoretical limitations, numerical obstacles, or an inability to address the complexities of the delay integral kernels. In this paper, we propose a unified framework to tackle the above problem by employing the concept of the Kronecker-Seuret decomposition (KSD) for matrix-valued functions, which we have recently developed for the analysis of complex delay structures in conjunction with the Krasovskiĭ functional approach. Our strategy can simultaneously address two distinct control problems, where the matrix kernels of DDs can contain an arbitrary finite number of square-integrable functions. We show in detail how the KSD can factorize and approximate different kernel functions simultaneously within an exact lifted formulation of the original problem. Furthermore, the use of KSD also enables us to construct complete-type functionals, whose integral kernels can include an arbitrary finite number of weakly differentiable functions that are linearly independent, underpinned by the utilization of general integral inequalities derived from the least-squares principle. The solution to each synthesis problem comprises two theorems accompanied by an iterative algorithm, which can be utilized as a single package to compute controller gains, thus eliminating the need for nonlinear solvers. We present the testing results of two challenging examples, which could not be addressed by existing methods, to demonstrate the effectiveness of our methodology. Additionally, the paper reviews recent advancements in the research of time-delay systems, which may provide useful information for both emerging and established researchers.

## 1. Introduction and Review

The study of time-delay systems (TDSs) boasts a rich history, with foundational work dating back to the early 20th century [1, 2]. Research in this field spans a variety of disciplines, such as functional differential equations (FDEs) [3–6], infinite-dimensional systems (IDSs) [7–10], coupled differential-difference equations [11–13], and numerous emerging topics in the field of cybernetics [14–20]. TDSs are used to model transport, propagation, and aftereffects in the state, input, and output of dynamical systems, often using pointwise and distributed delays (DDs) as tools. The nature of pointwise delays is well described in [21] as a transport equation coupled with boundary conditions. Additionally, delays can arise from transport media with intricate structures. A distributed delay [22–24] over delay interval $[-r, 0]$ is denoted by integral

$$\int_{-r}^{0} D(\tau)\boldsymbol{x}(t+\tau)\mathrm{d}\tau, \quad r > 0 \tag{1}$$

with a matrix-valued function $D(\cdot)$, where the integral aggregates a segment of the past dynamics [25–27]. Systems with both pointwise and distributed delays find extensive applications in diverse fields such as networked control [28–30], thermal dynamics [31, 32], neural networks [33], systems with predictor controllers [34], and the modeling of control saturation [35] and event-triggered mechanisms [36]. However, analyzing systems with general delay structures presents significant challenges and remains an area of active research [37–39].

Methods for the control and observation of linear time-delay systems (LTDSs) have been developed in the time [40–44] or frequency domain [45–48], making use of techniques from real and complex analysis [49]. A noteworthy observation pertinent to this topic is that the primary impediments to research advancement are more computational than theoretical in nature.

*This work was partially supported by the National Natural Science Foundation of China under Grant Nos. 62303180, 62373150 and 62273145, and Fundamental Research Funds for Central Universities under Grant 2025JC009, China.

*Corresponding author
*Email addresses:* qianfeng@ncepu.edu.cn, qfen204@aucklanduni.ac.nz (Qian Feng), W.Zheng@westernsydney.edu.au (Wei Xing Zheng), wangxiaoyu@bwu.edu.cn (Xiaoyu Wang), fengxiao@ncepu.edu.cn (Feng Xiao)

This characteristic will be illustrated by reviewing several categories of approaches that follow.

**Infinite-dimensional Systems/Time Domain:** Leveraging the theories of one-parameter semigroups [50–52], an LTDS can be equivalently described by an abstract differential equation [53] on an infinite-dimensional Banach space [54, 55], resulting in a complete and elegant characterization. By extending the classical optimal control theory [56], infinite-dimensional controllers/observers can be obtained by solving the linear-quadratic (LQ) problem over an infinite time horizon. Yet, addressing the LQ problem [57] in infinite dimension involves solving operator Riccati equations [58–60] that are difficult to compute numerically [61–63]. In fact, even the unknown parameters of the controllers/observers in this framework are given as abstract operators rather than in explicit concrete forms.

**Operator Theory/Frequency Domain:** The infinite-dimensional framework for LTDSs can also be exploited in the frequency domain in various ways as shown in [8, 47, 64]. LTDSs can be addressed using techniques in complex analysis to analyze the spectrum [48, 64] or transfer matrix of the underlying systems [47, 65]. Several strategies fall into this category, such as the algebraic/factorization approach [66–68], rational approximations [69–71], averaging approximations [72, 73] and spectral decomposition [74–76]. A notable insight from an application of spectral decomposition reveals that for a stabilizable/detectable LTDS [77–79], controllers/observers can always be constructed by stabilizing their unstable components using methods for simple LTI systems. This relies on the property that the spectrum of a retarded LTDS can only contain a finite number of unstable characteristic roots [80]. Nonetheless, this approach hinges on the feasibility of numerically computing all the unstable characteristic roots, which in itself represents a challenging aspect in the research of LTDSs [81–83].

**Eigenvalue Optimization/Frequency Domain:** A prominent recent line in frequency-domain-based methods for the stabilization of LTDSs involves computing the spectral abscissa (SA), making use of the techniques for eigenvalue optimization [80, 84, 85]. The rapid progression of these works is primarily driven by the recent breakthroughs in non-smooth optimization algorithms [86–89]. However, DDs with general matrix kernels have rarely been investigated with this framework primarily due to the complexities in handling the Laplace transform of the integral in (1), as $\int_{-r}^{0} D(\tau)\mathrm{e}^{\tau s}\mathrm{d}\tau$ may not have a closed-form expression [90–93] in $s \in \mathbb{C}$.

**FDE/Time Domain:** As demonstrated by [3], a general LTDS can always be expressed as an FDE using a Lebesgue-Stieltjes integral [49] by the Riesz representation theorem, where the matrix-valued integral kernels are of bounded variation that are right-continuous. A stabilizing controller/observer can then be constructed by using the state-transition matrix of the underlying system [77–79]. However, the fundamental solution matrix of an LTDS is not always computable via the method of steps [3], especially when the system includes general DDs as in (1). Similar to the state-space infinite-dimensional framework, controllers/observers can be constructed by solving the LQ problem over an infinite time horizon. The solution is formulated in terms of Riccati partial differential/integral equations [94–97], which are still difficult to solve numerically. Another representative method within this category involves the construction of predictor controllers [34], also known as the reduction approach [98], designed to equivalently eliminate the effect of input/output delays [99–103] and accommodate linear systems with both state and input delays [104–107]. A similar idea, called (Modified) Smith predictor [108–110], can be carried out in the frequency domain as well, using various techniques in abstract algebra [46, 111]. Notably, new constructive controller synthesis methods were proposed, which are inspired by the principle of predictor controllers, for LTDSs with particular structures such as backstepping and forwarding [112–114].

**Krasovskiĭ Functional/Time Domain:** This category of methods is an extension of the Lyapunov second method [115, 116] for FDE [41, 117, 118]. Unlike a delay-free LTI system whose stability can be analyzed by a quadratic function, functionals with integrals [119] are required for the stability analysis of general LTDSs, due to the inherent infinite dimension of the system's states. In general, two types of methodologies have been developed based on the construction of Krasovskiĭ functionals (KFs). The first one relies on the analysis of the complete Krasovskiĭ functionals (CKFs) [118], the existence of which furnishes a necessary and sufficient condition for the stability of LTDSs. The condition for CKF existence can be formulated using the so-called delay Lyapunov matrix (DLM) [120–126], which is the solution to functional differential matrix equations with boundary conditions. However, the computation of DLM [127–130] is a nontrivial task for LTDSs with general DDs, which is comparable to the computation of operator Riccati equations as we mentioned earlier. Furthermore, necessary stability conditions [127, 128, 131, 132] characterized by DLM are contingent on the assumption that the underlying system is already stable. Therefore, they cannot be utilized to determine the stability of an LTDS, nor to design a stabilizing controller for an unstable LTDS [133]. An interesting approach combining both time and frequency domain techniques was proposed in [134] for the stabilization of LTDSs with DDs, based on the concept of smoothed SA [135, 136] and DLM. This method, which requires computing DLM and its derivatives, is innovative but lacks detailed demonstration of how it can be carried out for an LTDS with general DDs or non-commensurate delays [137].

A more direct method of using KFs is to construct the functionals with fixed structures [43, 44, 138] parameterized by unknown matrices. This forms the basis for developing sufficient conditions for the existence of controllers and observers, expressed as convex optimization problems [139]. This methodology has been shown to be effective for the stability analysis and stabilization of LTDSs [31, 32, 140–150], supported by the development of efficient numerical algorithms for semidefinite programming (SDP) [151]. For a comprehensive review of this topic, the monographs [43, 44] are recommended. The main drawback of this approach, however, is that only sufficient conditions could be established, where the induced conservatism will primarily depend on the generality of the structures of KFs and the integral inequalities [152] utilized to construct them. Additionally, the use of congruent transformations alone may be insufficient for translating the original stability analysis condition into convex controller/observer synthesis conditions, as KFs in such contexts are not merely quadratic forms, adding to the complexity of the process.

Based on our literature review, it is clear that there are no effective solutions to the stabilization of LTDSs with an arbitrary finite number of pointwise and general distributed delays. Most existing approaches only consider pointwise delays [153–159], or impose conservative constraints on the structure of the state space matrices [160, 161] or limit the number of delays [144, 162] or the structure [152] of the DD kernels in $D(\cdot)$. While some methods theoretically could address this issue as mentioned earlier, they often face numerical challenges, such as those involving solving operator [61–63] or partial/integral Riccati equations [94–97], or computing DLMs [122, 125, 126] and the characteristic roots [48, 163] of the underlying systems. This observation underscores the need for solutions that rely on tractable numerical computations. Recently, several new approaches [144, 152, 162] have been developed specifically to address the presence of general DDs in an LTDS in conjunction with the KF and SDP approach, under the assumption that the delay values are known. Notably, we have recently proposed the concept of Kronecker-Seuret decomposition (KSD) [164], which can handle any finite number of $\mathscr{L}^2$ functions in the $D(\tau)$ of integral (1), without limiting the number of delays. This methodology has been successfully utilized in solving the dissipative robust estimation problem of a general LTDS in [164], where the solution is formulated using convex SDP constraints that can be efficiently computed by standard SDP solvers [165–167].

In this paper, we introduce a comprehensive framework, based on the strategy we recently proposed in [164], to design dissipative controllers for LTDSs with general delay structures by constructing a complete-type KF. As emphasized above, we highlight that the control problem we address has not been adequately resolved in existing literature. The generality of the system model is ensured by incorporating an arbitrary finite number of pointwise and general distributed delays at the states, inputs, and outputs, where the kernel matrices of all DDs can contain any number of $\mathscr{L}^2$ functions over bounded intervals. Thanks to the incorporation of these general delays, our system can be regarded as a realization of the general LTDS denoted by a Lebesgue-Stieltjes integral. Moreover, the system's differential equations are defined in accordance with the Carathéodory conditions [3, section 2.6] with respect to (w.r.t) the Lebesgue measure [49], which is well-suited for accommodating certain discontinuities and glitches in the system parameters.

To address the general delays in our system in conjunction with the KF approach, we adopted the KSD concept [164] that significantly generalizes our previous research on DDs with constant delay values [152, 162]. KSD allows for the decomposition of integral kernel matrices of any DD containing $\mathscr{L}^2$ functions as products between constant matrices and a list of basis functions with specific properties. The basis functions in our KSD are organized into a vector with three distinct components, where the first component is approximated by the rest of the components using the least-squares principle [55, page 182]. This setup gives users the liberty to factorize any portion of the $\mathscr{L}^2$ integral kernels or approximate them with any number of weakly differentiable and linearly independent functions (WDLIFs) along with their weak derivatives [168]. Benefiting from the structures of KSD, we utilize the third component of the basis function vector as the integral kernels for our general KF, which can contain any number of WDLIFs. This feature ensures the generality of our KF compared to existing results, as the generality of KFs is primarily determined by the complexity of their integral kernels. Finally, the construction of our KF is carried out using integral inequalities derived from the least-squares principle [55, page 182]. The use of all the preceding mathematical tools is intended to minimize the conservatism in the proposed framework.

Once our KF is successfully constructed, we formulate the controller synthesis conditions in a theorem as a solution to the dissipative control problem. Importantly, this theorem furnishes a framework that can address two distinct synthesis problems: those involving delays in either the input or the controller. Next, a second theorem is proposed to convexify the bilinear matrix inequality (BMI) in the synthesis conditions using [169–171, Projection Lemma], without compromising the matrix parameters in our KF. To further minimize conservatism, we set forth algorithms that can solve the BMI iteratively, which can be initialized by a feasible solution to the second theorem. As a result, our approach avoids using nonlinear SDP solvers [172–174] or other difficult numerical schemes [61, 86, 97, 130], while being able to address the presence of the general delay structures without resorting to simplifications.

The study of various cybernetic systems, such as networked control systems [28] and neural networks [33], has demonstrated that advanced methods for general LTDSs can be adapted to efficiently solve practical control problems when the underlying systems are associated with delay structures [35, 175–177]. Given the generality of our LTDS with the incorporation of dissipative constraints, the investigated problem could cover a wide range of LTDS-related problems in engineering contexts. Consequently, the proposed framework may serve as a blueprint and provide forward guidance for developing new solutions to open problems that arise from LTDS-related systems, such as coupled PDE-ODE systems [178–182], and sampled-data systems [183–186].

The remainder of the paper is organized into five sections. Section 2 opens with a motivating example to elucidate the core concept of the KSD, which is specifically developed for addressing the DDs in an LTDS in conjunction with the KF approach. Next, the formulation of our dissipative stabilization problem with the resulting CLS is presented in Section 3 by leveraging the KSD concept, and the solution to the synthesis problem is set out in Section 4. In Appendix F, we present a tutorial to demonstrate how KSD works with the proposed framework using specific examples. Moreover, Section 5 elaborates on how this framework is extended to design controllers with delays when input delays are absent. Finally, the testing and simulation results of two numerical examples are provided in Section 6 prior to the final conclusion. Some essential lemmas and proofs are included in the appendix to ensure minimal disruptions to the flow of the presentation.

**Notation.** Let $\mathbb{N} := \{1, 2, 3, \cdots\}$ denote the natural numbers and $\mathbb{N}_\nu := \{1, \ldots, \nu\} \subset \mathbb{N}$, and $\mathbb{S}^n := \{X \in \mathbb{R}^{n\times n} : X = X^\top\}$ the set of symmetric matrices, and $\mathcal{C}(\mathcal{X}; \mathbb{R}^n)$ the space of continuous functions with $\mathcal{X} \subseteq \mathbb{R}$. Set $\mathcal{M}\left(\mathcal{X}; \mathbb{R}^d\right)$ includes all measurable functions from $\mathcal{X}$ to $\mathbb{R}^d$ endowed with the Borel algebra. We use the standard Euclidean norm $\|\mathbf{x}\|_2 = \|\mathbf{x}\| := \sqrt{\mathbf{x}^\top \mathbf{x}}$, the $p$-seminorm $\|\boldsymbol{f}(\cdot)\|_p := \left(\int_{\mathcal{X}} \|\boldsymbol{f}(x)\|_2^p \,\mathsf{d}x\right)^{1/p}$, and define $\mathcal{L}^p(\mathcal{X}; \mathbb{R}^n) := \{\boldsymbol{f}(\cdot) \in \mathcal{M}(\mathcal{X}; \mathbb{R}^n) : \|\boldsymbol{f}(\cdot)\|_p < +\infty\}$ and the Sobolev space $\mathcal{W}^{1,2}(\mathcal{X}; \mathbb{R}^n) = \{\boldsymbol{f}(\cdot) \in \mathcal{L}^2(\mathcal{X}; \mathbb{R}^n) : \boldsymbol{f}'(\cdot) \in \mathcal{L}^2(\mathcal{X}; \mathbb{R}^n)\}$, where $\boldsymbol{f}'(\cdot)$ is the weak derivative. We write $\mathsf{Sy}(X) := X + X^\top$, $X \oplus Y = \begin{bmatrix} X & \mathsf{O}_{n,q} \\ \mathsf{O}_{p,m} & Y \end{bmatrix}$ for $X \in \mathbb{C}^{n\times m}$, $Y \in \mathbb{C}^{p\times q}$, and $\mathsf{diag}_{i=1}^{\nu} X_i = X_1 \oplus \cdots \oplus X_\nu$; $\otimes$ denotes the Kronecker product. The shorthand $\begin{bmatrix} A & B \\ * & C \end{bmatrix} = \begin{bmatrix} A & B \\ B^\top & C \end{bmatrix}$ is used for symmetric matrices, along with $[*]YX = X^\top Y X$ and $X^\top Y[*] = X^\top Y X$. Column and row concatenations are denoted $\mathsf{\textbf{Col}}_{i=1}^{n} X_i = [X_i]_{i=1}^{n} := \left[X_1^\top \cdots X_n^\top\right]^\top$ and $\mathsf{\textbf{Row}}_{i=1}^{n} X_i = [\![X_i]\!]_{i=1}^{n} = \left[X_1 \cdots X_n\right]$, respectively. The notation $\widetilde{\forall}\, x \in \mathcal{X}$ means $x \in \mathcal{X}$ holds **almost everywhere** w.r.t. the Lebesgue measure; $\mathsf{O}_{n,m}$ is the $n \times m$ zero matrix (abbreviated $\mathsf{O}_n$ when $n = m$), and $\mathbf{0}_n$ is the $n \times 1$ zero vector. The precedence of matrix operations is: *matrix (scalar) multiplications* $> \otimes > \oplus = \mathsf{diag} > +$. Empty matrices $[\,]$ follow the definition in [187, See I.7] and the computational rules in the MATLAB© platform: $I_0 = [\,]_{0,0} = [\,]_{0,0}^{-1}$, $\mathsf{O}_{0,m} = [\,]_{0,m}$, $\mathsf{O}_{n,0} = [\,]_{n,0}$, $[X_i]_{i=1}^{0} = [\,]_{0,m}$, and $[\![X_i]\!]_{i=1}^{0} = [\,]_{n,0}$ for $X_i \in \mathbb{R}^{n\times m}$, and $[\,]_{n,0} \times [\,]_{0,m} = \mathsf{O}_{n,m}$.

**Acronyms:**

DDs: Distributed Delays

KSD: Kronecker-Seuret Decomposition

KFs: Krasovskiĭ Functionals

FDEs: Functional Differential Equations

LTDSs: Linear Time-Delay Systems

SA: Spectral Abscissa

SDP: Semidefinite Programming

WDLIFs: Weakly Differentiable and Linearly Independent Functions

BMIs: Bilinear Matrix Inequalities

SRFs: Supply Rate Functions

DSFC: Dissipative State Feedback Control

CLSs: Closed-Loop Systems

## 2. Motivating Example and the KSD Concept

### 2.1. Motivating Example

We begin by using a specific example to show the main challenges in addressing controller synthesis problems for LTDSs, in both the time and frequency domains. Consider

$$\dot{\boldsymbol{x}}(t) = A_1 \boldsymbol{x}(t) + A_2 \boldsymbol{x}(t - r) + \int_{-r}^{0} A_3(\tau) \boldsymbol{x}(t + \tau) \mathsf{d}\tau, \quad (2)$$

where $A_1, A_2 \in \mathbb{R}^{n\times n}$, $A_3(\cdot) \in \mathcal{L}^2\left([-r, 0]; \mathbb{R}^{n\times n}\right)$ with delay value $r > 0$. In the frequency domain, the spectrum of LTDS (2) is $\mathbb{A} = \left\{s \in \mathbb{C} : p(s) = 0\right\}$ where

$$p(s) = \det\left(sI_n - A_1 - A_2 \mathsf{e}^{-rs} - \int_{-r}^{0} A_3(\tau) \mathsf{e}^{\tau s} \mathsf{d}\tau\right). \quad (3)$$

The set $\mathbb{A}$ contains the characteristic roots of LTDS (2), which are generally countably infinite in number due to the presence of transcendental terms $\mathsf{e}^{-rs}$ and $\int_{-r}^{0} A_3(\tau)\mathsf{e}^{\tau s}\mathsf{d}\tau$ in (3). This underlies why TDSs are generally infinite-dimensional in contrast to systems obeying ODEs.

Of widespread interest in the field of automatic control is to determine the stability of complex systems without explicitly deriving their solutions in closed form. From analyzing the spectrum set $\mathbb{A}$ with (3), the equilibrium of (2) is exponentially stable if and only if the SA of (2) satisfies $\max_{\lambda\in\mathbb{A}} \Re(\lambda) < 0$. This criterion can be verified by computing a finite number of roots thanks to the notable properties that there can be at most [3, 80, 188] a finite number of $\lambda \in \mathbb{A}$ with $\Re(\lambda) > 0$. The pseudo-spectral techniques in [188, 189] represent the latest methods that can compute the SA of (2) by discretizing either the solution operator or the infinitesimal generator [190] in $\mathbb{C}$, where $\int_{-r}^{0} A_3(\tau)\mathsf{e}^{\tau s}\mathsf{d}\tau$ is discretized via ordinary quadratures [191]. Without considering computational complexities, these methods could consistently determine the stability of (2) for a given $r > 0$ if the discretization index is sufficiently large, irrespective of the complexities of the matrix kernels of $A_3(\cdot)$.

Unfortunately, the techniques in [188, 189] cannot resolve the controller design problem for (2) where the spectrum relates to some parameters $K$ of the controller. It is a common issue among frequency domain methods [156] that solutions to stability-analysis problems do not extend readily to controller synthesis. This challenge is particularly pronounced in cases wherein the design involves performance objectives such as $\mathcal{H}^\infty/\mathcal{H}^2$ minimization [85]. A viable approach to find a stabilizing $K$ when $A_3(\tau) \equiv \mathsf{O}_n$ is to compute a locally optimal solution [80, 157] for $\min_K \max_{\lambda\in\mathbb{A}} \Re(\lambda)$, where the SA is pushed to the left in $\mathbb{C}$ as far as possible. However, the optimization problem is inherently non-smooth due to the properties of SA, requiring the use of specialized algorithms such as the BFGS-SQP method [192, 193], gradient sampling [194, 195], bundle trust-region [89, 196] to perform computations. To the best of our knowledge, no broadly applicable stabilization method along this spectral-abscissa route is currently available for the case $A_3(\tau) \not\equiv \mathsf{O}_n$, largely due to the complex-ity of $A_3(\cdot)$. Lastly, any numerical method in the frequency domain involving discretization or rational approximations [197] would face significant computational burden when a large delay $r > 0$ is discretized with an extensive number of points.

In the time domain, various established methods [57, 198] could, in theory, address the stabilization of (2) by converting the original problem into the task of solving complex equations. However, as discussed in Section 1, numerically solving these equations is generally challenging [62, 96]. In this article, we focus on using the KF approach for (2) to construct tractable controller synthesis conditions solvable by standard SDP numerical solvers. The main idea of the KF approach involves constructing a functional of the form

$$V(\boldsymbol{x}(t+\cdot)) = [*]\begin{bmatrix} P_1 & P_2 \\ * & P_3 \end{bmatrix}\begin{bmatrix} \boldsymbol{x}(t) \\ \int_{-r}^{0} F(\tau)\boldsymbol{x}(t+\tau)\mathsf{d}\tau \end{bmatrix} + \int_{-r}^{0} \boldsymbol{x}^{\top}(t+\tau)\left[Q + (\tau + r)R\right]\boldsymbol{x}(t+\tau)\mathsf{d}\tau, \quad (4)$$

where $F(\cdot)$ is a matrix-valued function with appropriate dimension which has absolutely continuous kernels, and matrices $Q, R \in \mathbb{S}^n$ and $P_1, P_2, P_3$ are matrices of appropriate dimensions determined by $F(\tau)$. Given (2) and (4) and the common procedures in constructing KFs, the main obstacles to formulating tractable conditions include:

1. $\int_{-r}^{0} A_3(\tau)\boldsymbol{x}(t+\tau)\mathsf{d}\tau$ is inevitably present in $\dot{V}$, given $\dot{\boldsymbol{x}}(t)$ in (2) and the structure of $V(\boldsymbol{x}(t+\cdot))$ in (4).

2. Computing $\dot{V}(\boldsymbol{x}(t+\cdot))$ inevitably involves the derivative of $\int_{-r}^{0} F(\tau)\boldsymbol{x}(t+\tau)\mathsf{d}\tau$ with respect to $t$, namely,

   $$\frac{\mathsf{d}}{\mathsf{d}t}\int_{-r}^{0} F(\tau)\boldsymbol{x}(t+\tau)\mathsf{d}\tau = \frac{\mathsf{d}}{\mathsf{d}t}\int_{t-r}^{t} F(\tau - t)\boldsymbol{x}(\tau)\mathsf{d}\tau = F(0)\boldsymbol{x}(t) - F(-r)\boldsymbol{x}(t-r) - \int_{-r}^{0}\frac{\mathsf{d}F(\tau)}{\mathsf{d}\tau}\boldsymbol{x}(t+\tau)\mathsf{d}\tau \quad (5)$$

3. For $\dot{V}(\boldsymbol{x}(t+\cdot))$ to be expressible in quadratic-like forms that enable the construction of SDP constraints, there must exist certain relationships between $A_3(\tau)$, $F(\tau)$ and $F'(\tau)$ that can address all the kernel functions simultaneously.

4. The infinite dimension of $A_3(\cdot)$, $F(\cdot)$ and $F'(\cdot)$ needs to be addressed by some mechanisms, since directly including them in the synthesis conditions would necessitate solving SDP problems with infinite dimension.

5. Users should always have sufficient degrees of freedom to choose different $F(\cdot)$ that can increase the generality of the $V(\boldsymbol{x}(t+\cdot))$, which in turn may improve the feasibility of the resulting SDP conditions. At the same time, it must be ensured that all kernels in $A_3(\tau)$, $F(\tau)$, and $F'(\tau)$ are handled within a unified tractable representation.

### 2.2. *KSD for Matrix Valued-Functions*

To address the challenges arising from the infinite dimensionality of matrix-valued functions, we introduced the concept of KSD in [164] to design state observers for general linear delay systems. Since several DDs can be found in a delay system, KSD should possess the capacity to handle multiple matrix-valued functions simultaneously. The following proposition presents the first ingredient of the concept of KSD.

**Proposition 1.** Assume matrix-valued functions $\mathcal{F}_{i,j}(\cdot) \in \mathscr{L}^2(\mathcal{I}_i\,;\mathbb{R}^{n_j\times m_j})$, $\forall i \in \mathbb{N}_\nu$, $\forall j \in \mathbb{N}_\hbar$ with $\nu, \hbar, n_j, m_j \in \mathbb{N}$, where $\mathcal{I}_i \subset \mathbb{R}$ and the Lebesgue measure of $\mathcal{I}_i$ is non-zero. The above membership conditions hold if and only if there exist $\boldsymbol{f}_i(\cdot) \in \mathscr{W}^{1,2}(\mathcal{I}_i\,;\mathbb{R}^{d_i})$, $\boldsymbol{\varphi}_i(\cdot) \in \mathscr{L}^2(\mathcal{I}_i\,;\mathbb{R}^{\delta_i})$, $\boldsymbol{\phi}_i(\cdot) \in \mathscr{L}^2(\mathcal{I}_i\,;\mathbb{R}^{\mu_i})$ and $M_i \in \mathbb{R}^{d_i\times\varkappa_i}$, $\gimel_{i,j} \in \mathbb{R}^{n_j\times\kappa_i m_j}$, such that

$$\mathcal{F}_{i,j}(\tau) = \gimel_{i,j}\left(\boldsymbol{g}_i(\tau)\otimes I_{m_j}\right), \quad (6a)$$

$$\frac{\mathsf{d}\boldsymbol{f}_i(\tau)}{\mathsf{d}\tau} = M_i\boldsymbol{h}_i(\tau), \qquad \boldsymbol{h}_i(\tau) = \begin{bmatrix}\boldsymbol{\varphi}_i(\tau)\\ \boldsymbol{f}_i(\tau)\end{bmatrix} \in \mathbb{R}^{\varkappa_i} \quad (6b)$$

$$\mathsf{G}_i = \int_{\mathcal{I}_i}\boldsymbol{g}_i(\tau)\boldsymbol{g}_i^{\top}(\tau)\mathsf{d}\tau \succ 0, \quad \boldsymbol{g}_i(\tau) = \begin{bmatrix}\boldsymbol{\phi}_i(\tau)\\ \boldsymbol{h}_i(\tau)\end{bmatrix} \in \mathbb{R}^{\kappa_i} \quad (6c)$$

for all $\tau \in \mathcal{I}_i$ and $i \in \mathbb{N}_\nu$ and $j \in \mathbb{N}_\hbar$, where $\nu, \hbar \in \mathbb{N}$ and $\kappa_i = d_i + \delta_i + \mu_i$, $\varkappa_i = d_i + \delta_i$ with indices $d_i \in \mathbb{N}$ and $\delta_i, \mu_i \in \mathbb{N}\cup\{0\}$ for all $i \in \mathbb{N}_\nu$. The derivatives in (6b) are weak derivatives [168].

*Proof.* It is obvious that $\mathcal{F}_{i,j}(\cdot) \in \mathscr{L}^2(\mathcal{I}_i\,;\mathbb{R}^{n_j\times m_j})$ is implied by (6a)–(6c) for all $i \in \mathbb{N}_\nu$ and $j \in \mathbb{N}_\hbar$, given the definitions of $\boldsymbol{\varphi}_i(\cdot), \boldsymbol{f}_i(\cdot), \boldsymbol{\phi}_i(\cdot)$ and the fact that $\mathscr{W}^{1,2}(\mathcal{I}_i\,;\mathbb{R}^{d_i}) \subset \mathscr{L}^2(\mathcal{I}_i\,;\mathbb{R}^{d_i})$. Hence, sufficiency is proven.

We now proceed to establish necessity. Namely, we show that $\mathcal{F}_{i,j}(\cdot) \in \mathscr{L}^2(\mathcal{I}_i\,;\mathbb{R}^{n_j\times m_j})$ implies the existence of the parameters in Proposition 1 satisfying (6a)–(6c). Let $i \in \mathbb{N}_\nu$, given any vector-valued function $\boldsymbol{f}_i(\cdot) \in \mathscr{W}^{1,2}(\mathcal{I}_i\,;\mathbb{R}^{d_i})$ satisfying $\int_{\mathcal{I}_i}\boldsymbol{f}_i(\tau)\boldsymbol{f}_i^{\top}(\tau)\mathsf{d}\tau \succ 0$, we can always select appropriate vector-valued functions $\boldsymbol{\phi}_i(\cdot) \in \mathscr{L}^2(\mathcal{I}_i\,;\mathbb{R}^{\mu_i})$ and $\boldsymbol{\varphi}_i(\cdot) \in \mathscr{L}^2(\mathcal{I}_i\,;\mathbb{R}^{\delta_i})$ with some constant matrices $M_i \in \mathbb{R}^{d_i\times\varkappa_i}$ to satisfy (6b)–(6c), because of $\boldsymbol{f}_i'(\cdot) \in \mathscr{L}^2(\mathcal{I}_i\,;\mathbb{R}^{d_i})$. This is true because we can always stack $\boldsymbol{\varphi}_i(\cdot)$ and $\boldsymbol{\phi}_i(\cdot)$ with new functions vertically that satisfy (6c), as long as they are linearly independent in a Lebesgue sense. We note that $\int_{\mathcal{I}_i}\boldsymbol{f}_i(\tau)\boldsymbol{f}_i^{\top}(\tau)\mathsf{d}\tau \succ 0$ is implied by (6c) which also implies that the functions in $\boldsymbol{g}_i(\cdot)$ in (6c) are linearly independent [199, Th. 7.2.10] in a Lebesgue sense over $\mathcal{I}_i$. Note that also $\boldsymbol{\varphi}_i(\tau)$ or $\boldsymbol{\phi}_i(\tau)$ can be an empty matrix $[\,]_{0\times 1}$.

Since $\kappa_i = \dim\left(\boldsymbol{g}_i(\tau)\right) = \dim\begin{bmatrix}\boldsymbol{\phi}_i^{\top}(\tau) & \boldsymbol{\varphi}_i^{\top}(\tau) & \boldsymbol{f}_i^{\top}(\tau)\end{bmatrix}^{\top}$ in (6b)–(6c) can be increased without fixed bounds, which means that we can stack $\boldsymbol{g}_i(\cdot)$ with an arbitrary number of new functions as long as (6c) is satisfied, we can always find appropriate $\widehat{F}_{i,j,k} \in \mathbb{R}^{n_j\times m_j}$ and $\boldsymbol{g}_i(\tau) = \left[g_{i,k}(\tau)\right]_{k=1}^{\kappa_i}$ satisfying (6c) such that

$$\mathcal{F}_{i,j}(\tau) = \sum_{k=1}^{\kappa_i}\widehat{F}_{i,j,k}\,g_{i,k}(\tau) = \left[\!\left[\widehat{F}_{i,j,k}\right]\!\right]_{k=1}^{\kappa_i}\left(\boldsymbol{g}_i(\tau)\otimes I_{m_j}\right), \quad (7)$$

for all $\tau \in \mathcal{I}_i$, $i \in \mathbb{N}_\nu$, $j \in \mathbb{N}_\hbar$ with $\kappa_i = d_i + \delta_i + \mu_i$, where functions $\boldsymbol{\varphi}_i(\cdot) \in \mathscr{L}^2(\mathcal{I}_i;\mathbb{R}^{\delta_i})$, $\boldsymbol{f}_i(\cdot) \in \mathscr{W}^{1,2}(\mathcal{I}_i\,;\mathbb{R}^{d_i})$ and $\boldsymbol{\phi}_i(\cdot) \in \mathscr{L}^2(\mathcal{I}_i\,;\mathbb{R}^{\mu_i})$ satisfy (6b)–(6c) for some $M_i \in \mathbb{R}^{d_i\times\varkappa_i}$. By setting $\gimel_{i,j} = \left[\!\left[\widehat{F}_{i,j,k}\right]\!\right]_{k=1}^{\kappa_i}$, (7) becomes (6a). ∎

**Remark 1.** Since $\int_{\mathcal{I}_i} p_i(\tau)q_i(\tau)\mathrm{d}\tau$ constitutes inner products, $\mathsf{G}_i$ in (6c) is the Gram matrix [199, Th. 7.2.10] of the functions in $\boldsymbol{g}_i(\cdot) = [g_{ij}(\cdot)]_{j=1}^{\kappa_i}$. Moreover, $\mathsf{G}_i \succ 0$ in (6c) indicates [199, Th. 7.2.10] that $g_{i,j}(\cdot)$, $j = 1\cdots\kappa_i$ are linearly independent in a Lebesgue sense over $\mathcal{I}_i$. An example of $\mathsf{G}_i$ is $\int_{-1}^{0}[\tau^i]_{i=0}^{1}\llbracket\tau\rrbracket_{i=0}^{1}\mathrm{d}\tau = \int_{-1}^{0}\left[\begin{smallmatrix}1 & \tau\\ * & \tau^2\end{smallmatrix}\right]\mathrm{d}\tau = \left[\begin{smallmatrix}1 & -1/2\\ * & 1/3\end{smallmatrix}\right] \succ 0$.

**Remark 2.** Proposition 1 is both necessary and sufficient for representing $\mathcal{F}_{i,j}(\cdot) \in \mathcal{L}^2\left(\mathcal{X};\mathbb{R}^{n_j\times m_j}\right)$. As shown in the proof of Proposition 1, any $\boldsymbol{f}_i(\cdot) \in \mathcal{W}^{1,2}(\mathcal{X};\mathbb{R}^{d_i})$ can satisfy the conditions in (6b), even if none of the functions in $\boldsymbol{f}_i(\cdot)$ are included in $\mathcal{F}_{i,j}(\cdot)$. This is because an arbitrary finite number of new functions, satisfying (6c), can be added to both $\boldsymbol{f}(\cdot)$ and $\boldsymbol{\varphi}_i(\cdot)$ such that (6b) holds for some matrix $M_i$. Hence, the decomposition in (6a)–(6c) is always achievable for $\mathcal{F}_{i,j}(\cdot)$.

Proposition 1 is the first element of KSD, which partitions $\boldsymbol{g}_i(\cdot) = \begin{bmatrix}\boldsymbol{\phi}_i^\top(\tau) & \boldsymbol{h}_i^\top(\tau)\end{bmatrix}^\top = \begin{bmatrix}\boldsymbol{\phi}_i^\top(\tau) & \boldsymbol{\varphi}_i^\top(\tau) & \boldsymbol{f}_i^\top(\tau)\end{bmatrix}^\top$ into three distinct components. Such configuration is justified by the next step of KSD, where $\boldsymbol{\phi}_i(\cdot)$ is addressed differently from $\boldsymbol{\varphi}_i(\cdot)$ and $\boldsymbol{f}_i(\cdot)$. Specifically, we assume each $\boldsymbol{\phi}_i(\cdot)$ is linearly approximated by $\boldsymbol{h}_i(\cdot) = \begin{bmatrix}\boldsymbol{\varphi}_i^\top(\cdot) & \boldsymbol{f}_i^\top(\cdot)\end{bmatrix}^\top$ as

$$\forall i \in \mathbb{N}_\nu,\ \forall\tau\in\mathcal{I}_i,\ \boldsymbol{\phi}_i(\tau) = \Gamma_i\,\mathfrak{H}_i^{-1}\boldsymbol{h}_i(\tau) + \boldsymbol{\varepsilon}_i(\tau) \tag{8}$$

where we define $\Gamma_i = \int_{\mathcal{I}_i}\boldsymbol{\phi}_i(\tau)\boldsymbol{h}_i^\top(\tau)\mathrm{d}\tau \in \mathbb{R}^{\mu_i\times\varkappa_i}$ and $\mathfrak{H}_i = \int_{\mathcal{I}_i}\boldsymbol{h}_i(\tau)\boldsymbol{h}_i^\top(\tau)\mathrm{d}\tau$. Note that $\mathfrak{H}_i \succ 0$ is implied by (6c). We also define $\boldsymbol{\varepsilon}_i(\tau) = \boldsymbol{\phi}_i(\tau) - \Gamma_i\,\mathfrak{H}_i^{-1}\boldsymbol{h}_i(\tau)$ as the approximation errors, and

$$\begin{aligned}\mathbb{S}^{\mu_i} \ni \mathfrak{E}_i &= \int_{\mathcal{I}_i}\boldsymbol{\varepsilon}_i(\tau)\boldsymbol{\varepsilon}_i^\top(\tau)\mathrm{d}\tau\\ &= \int_{\mathcal{I}_i}\boldsymbol{\phi}_i(\tau)\boldsymbol{\phi}_i^\top(\tau)\mathrm{d}\tau - \Gamma_i\,\mathfrak{H}_i^{-1}\Gamma_i^\top \succ 0\end{aligned} \tag{9}$$

is utilized to measure these errors, where $\mathfrak{E}_i \succ 0$ is ensured by the property of Schur complement of $\mathfrak{H}_i$ in $\mathsf{G}_i \succ 0$. Now by (6b) and the approximation in (8), it follows that

$$\begin{aligned}\boldsymbol{g}_i(\tau) &= \begin{bmatrix}\boldsymbol{\phi}_i(\tau)\\ \boldsymbol{h}_i(\tau)\end{bmatrix} = \begin{bmatrix}\Gamma_i\,\mathfrak{H}_i^{-1}\boldsymbol{h}_i(\tau)\\ \boldsymbol{h}_i(\tau)\end{bmatrix} + \begin{bmatrix}\boldsymbol{\varepsilon}_i(\tau)\\ \boldsymbol{0}_{\varkappa_i}\end{bmatrix}\\ &= \widehat{\Gamma}_i\boldsymbol{h}_i(\tau) + \widetilde{I}_i\boldsymbol{\varepsilon}_i(\tau),\end{aligned} \tag{10}$$

$$\widehat{\Gamma}_i = \begin{bmatrix}\Gamma_i\,\mathfrak{H}_i^{-1}\\ I_{\varkappa_i}\end{bmatrix} \in \mathbb{R}^{\kappa_i\times\varkappa_i},\quad \widetilde{I}_i = \begin{bmatrix}I_{\mu_i}\\ \mathsf{O}_{\varkappa_i,\mu_i}\end{bmatrix} \in \mathbb{R}^{\kappa_i\times\mu_i}$$

by using (6b) and the relation in (8) with matrices $\Gamma_i, \mathfrak{H}_i$ in (8). By replacing $\boldsymbol{g}_i(\cdot)$ in (6a)–(6c) with $\widehat{\Gamma}_i\boldsymbol{h}_i(\tau) + \widetilde{I}_i\boldsymbol{\varepsilon}_i(\tau)$ in (10), the presentation of the KSD concept for the matrix-valued functions $\mathcal{F}_{i,j}(\cdot) \in \mathcal{L}^2(\mathcal{I}_i;\mathbb{R}^{n_j\times m_j})$ is completed.

**Remark 3.** The decomposition for matrix-valued functions in [152] is a special case of Proposition 1 with $\nu = \hbar = 1$, $\phi_1(\cdot) = \varphi_1(\cdot) = [\,]_{0\times 1}$ and $\boldsymbol{f}_1(\cdot)$ satisfying

$$\tfrac{\mathrm{d}\boldsymbol{f}_1(\tau)}{\mathrm{d}\tau} = N\boldsymbol{f}_1(\tau)\quad\text{for some}\quad N \in \mathbb{R}^{d_1\times d_1}, \tag{11}$$

where no approximations take place. Additionally, the approximation scheme in [162] is a special case of our KSD with $\varphi_1(\cdot) = [\,]_{0\times 1}$ and $\boldsymbol{\phi}_1(\cdot) \in \mathcal{L}^2(\mathcal{X};\mathbb{R}^{\mu_1})$, where the approximator function $\boldsymbol{f}_1(\cdot)$ must satisfy (11). The differential equation in (11) implies that the functions in the $\boldsymbol{f}_1(\cdot)$ of [152, 162] must be the solutions to linear homogeneous differential equations with constant coefficients, which are far more restrictive than that of $\boldsymbol{f}_1(\cdot) \in \mathcal{W}^{1,2}(\mathcal{X};\mathbb{R}^{d_1})$ in Proposition 1. Moreover, [200, Proposition 1] is a special case of Proposition 1 with $\nu = 2$, $\hbar = 1$, $\phi_i(\cdot) = [\,]_{0\times 1}$, $\boldsymbol{\varphi}_i(\cdot) \in \mathcal{L}^2(\mathcal{I}_i;\mathbb{R}^{\delta_i})$ and $\boldsymbol{f}_i(\cdot) \in \mathcal{C}^1(\mathcal{I}_i;\mathbb{R}^{d_i}) \subset \mathcal{W}^{1,2}(\mathcal{I}_i;\mathbb{R}^{d_i})$, where no approximations are employed. Finally, the approximation framework in [144] is also a special case of our KSD with $\nu = 1$, $\varphi_1(\cdot) = [\,]_{0\times 1}$, $\phi_1(\cdot) \in \mathcal{C}(\mathcal{X};\mathbb{R})$, where $\boldsymbol{f}_1(\cdot)$ contains Legendre polynomials [201, 22.2.10] only.

The idea of Proposition 1 is not difficult to grasp. For any matrix-valued function $A(\tau) \in \mathbb{R}^{n\times n}$, we can always write

$$A(\tau) = \sum_{j=1}^{\kappa_1} A_j g_j(\tau) = \widehat{A}\left(\boldsymbol{g}(\tau)\otimes I_n\right),\quad \forall\tau\in\mathcal{X} \tag{12}$$

as a combination of some matrices $A_j \in \mathbb{R}^{n\times n}$ and linearly independent functions $\{g_j(\cdot)\}_{j=1}^{\kappa_1}$ with $\boldsymbol{g}(\cdot) = [g_j(\cdot)]_{j=1}^{\kappa_1}$, based on the properties of matrix-valued functions and the Kronecker product [199]. Moreover, the decomposition in (12) is constructible for multiple matrix-value functions $\mathcal{F}_{i,j}(\cdot) \in \mathcal{L}^2(\mathcal{I}_i;\mathbb{R}^{n_j\times m_j})$ using many functions $\boldsymbol{g}_i(\cdot)$ on each set $\mathcal{I}_i$ for every $j \in \mathbb{N}_\hbar$ simultaneously. The advantage of the decomposition in (12) is that $\widehat{A}$ has finite dimensions and the infinite dimension of $A(\cdot)$ is transferred to $\boldsymbol{g}(\cdot)$. What sets apart our Proposition 1, however, are the conditions in (6b)–(6c). Indeed, the relations in (6b) are formulated to facilitate the construction of KFs like (4) with $F(\tau) = \boldsymbol{f}(\tau)\otimes I_n$, so that (5) can be equivalently expressed by

$$\begin{aligned}\frac{\mathrm{d}}{\mathrm{d}t}\int_{-r}^{0}F(\tau)\boldsymbol{x}(t+\tau)\mathrm{d}\tau &= \frac{\mathrm{d}}{\mathrm{d}t}\int_{t-r}^{t}F(\tau-t)\boldsymbol{x}(\tau)\mathrm{d}\tau\\ &= F(0)\boldsymbol{x}(t) - F(-r)\boldsymbol{x}(t-r)\\ &\quad - (M\otimes I_n)\int_{-r}^{0}(\boldsymbol{h}(\tau)\otimes I_n)\boldsymbol{x}(t+\tau)\mathrm{d}\tau.\end{aligned} \tag{13}$$

Suppose $n_j = m_j$ in Proposition 1 with $\nu = \hbar = 1$. Then an application of the KSD in Proposition 1 with (13) can denote $A(\tau) = \widehat{A}\left(\boldsymbol{g}(\tau)\otimes I_m\right)$, $F(\tau) = \boldsymbol{f}(\tau)\otimes I_n$ and $F'(\tau) = M\boldsymbol{h}(\tau)\otimes I_n$ simultaneously using $\boldsymbol{g}(\tau)\otimes I_n$, given the relations in (6b) and $\boldsymbol{f}(\tau) = \begin{bmatrix}\mathsf{O}_{d_1,\delta_1} & I_{d_1}\end{bmatrix}\boldsymbol{h}(\tau)$. Meanwhile, the choice of $\boldsymbol{g}(\cdot)$ is not unique, and we can stack the column of $\boldsymbol{g}(\tau)$ with any number of new functions that satisfy (6b)–(6c). As a result, the machinery of the KSD concept is designed to address the theoretical obstacles we have outlined below (4) in the construction of (4) for (2).

**Remark 4.** The designation "KSD" is coined in light of the mathematical constructs in Proposition 1 and (8)–(10). The acronym begins with "Kronecker" to reflect the use of the Kronecker product in the decompositions in (6a). The structure of $A(\tau) = \widehat{A}\left(\boldsymbol{g}(\tau)\otimes I_n\right)$ can be found in various subjects such as sum-of-squares optimization[1] [202] and the construction of complete-type KFs [203]. Additionally, the idea of approximating the kernel functions was initially explored by A. Seuret and his colleagues in [144] in analyzing the stability of linear systems with a distributed delay, where the continuous kernel function of a scalar distributed delay is approximated by Legendre polynomials. In recognition of this significant contribution, we have utilized "Seuret" in KSD to indicate the use of function approximations.

## 3. Problem Formulation via KSD

In this section, we introduce the open-loop model and formulate the control-design problem using Proposition 1.

### 3.1. Open-Loop System Equations

Consider a general linear time-delay system

$$
\begin{aligned}
\dot{\boldsymbol{x}}(t) &= \mathcal{A}\mathbf{x}_t(\cdot) + \mathcal{B}_1\mathbf{u}_t(\cdot) + D_1\boldsymbol{w}(t), \quad \widetilde{\forall} t \geq t_0 \\
\boldsymbol{z}(t) &= \mathcal{C}\mathbf{x}_t(\cdot) + \mathcal{B}_2\mathbf{u}_t(\cdot) + D_2\boldsymbol{w}(t) \\
\mathbf{x}_t(\theta) &= \boldsymbol{x}(t+\theta), \quad \mathbf{u}_t(\theta) = \boldsymbol{u}(t+\theta), \quad \forall \theta \in \mathcal{J}
\end{aligned}
\tag{14}
$$

denoted by a linear FDE in an extended sense [3, page 58] for almost all $t \geq t_0 \in \mathbb{R}$ with respect to the Lebesgue measure, where $\mathcal{J} = [-r, 0]$ with known delay $r > 0$, operators $\mathcal{A} \in \mathcal{C}\left(\mathcal{C}\left(\mathcal{J};\mathbb{R}^n\right);\mathbb{R}^n\right)$, $\mathcal{B}_1 \in \mathcal{C}\left(\mathcal{C}\left(\mathcal{J};\mathbb{R}^p\right);\mathbb{R}^n\right)$, $\mathcal{B}_2 \in \mathcal{C}\left(\mathcal{C}\left(\mathcal{J};\mathbb{R}^p\right);\mathbb{R}^m\right)$ and $\mathcal{C} \in \mathcal{C}\left(\mathcal{C}\left(\mathcal{J};\mathbb{R}^n\right);\mathbb{R}^m\right)$ are linear, and $D_1 \in \mathbb{R}^{n\times q}$, $D_2 \in \mathbb{R}^{m\times q}$ are constant matrices. Moreover, $\boldsymbol{x}(t) \in \mathbb{R}^n$ is the solution to the FDE in the extended sense which is absolutely continuous for $t \geq t_0$, $\mathbf{x}_t(\cdot) \in \mathcal{C}\left(\mathcal{J};\mathbb{R}^n\right)$ is the state of the FDE which satisfies $\mathbf{x}_t(\cdot) \in \mathcal{W}^{1,2}\left(\mathcal{J};\mathbb{R}^n\right)$ for all $t \geq t_0 + r$, $\boldsymbol{u}(t) \in \mathbb{R}^p$ is the control input, $\boldsymbol{w}(\cdot) \in \mathcal{L}^2\left(\mathbb{R}_{\geq t_0};\mathbb{R}^q\right)$ stands for disturbances, and $\boldsymbol{z}(t) \in \mathbb{R}^m$ is the regulated output with dimension indices $n, m, p, q \in \mathbb{N}$.

**Remark 5.** Any linear delay structure can be represented by continuous linear operators as in (14) thanks to the Riesz representation theorem. For the general theory of linear FDEs, see [3, Chapters 6-7].

The FDE in (14) can be reformulated equivalently [3, Chapter 7] as an abstract differential equation on Hilbert (Banach) space using the theory of $\mathcal{C}_0$ semigroup [52, 190]. Note, however, that the infinitesimal generator of the $\mathcal{C}_0$ semigroup representation is not identical to operator $\mathcal{A}$ in (14), since the operators have different domains and ranges.

The linear operators in (14) are abstract and difficult to deal with directly. Thanks to the Riesz representation theorem [49], these operators can always be equivalently represented by

$$
\begin{aligned}
\mathcal{A}\mathbf{x}_t(\cdot) &= \int_{-r}^{0} \mathsf{d}[\mathsf{A}(\tau)]\,\mathbf{x}_t(\tau) = \int_{-r}^{0} \mathsf{d}[\mathsf{A}(\tau)]\,\boldsymbol{x}(t+\tau), \\
\mathcal{C}\mathbf{x}_t(\cdot) &= \int_{-r}^{0} \mathsf{d}[\mathsf{C}(\tau)]\,\mathbf{x}_t(\tau) = \int_{-r}^{0} \mathsf{d}[\mathsf{C}(\tau)]\,\boldsymbol{x}(t+\tau), \\
\mathcal{B}_i\mathbf{u}_t(\cdot) &= \int_{-r}^{0} \mathsf{d}[\mathsf{B}_i(\tau)]\,\mathbf{u}_t(\tau) = \int_{-r}^{0} \mathsf{d}[\mathsf{B}_i(\tau)]\,\boldsymbol{u}(t+\tau),
\end{aligned}
\tag{15}
$$

by the use of Lebesgue-Stieltjes integrals, where $\mathsf{A}(\tau) \in \mathbb{R}^{n\times n}$, $\mathsf{C}(\tau) \in \mathbb{R}^{m\times n}$, $\mathsf{B}_1(\tau) \in \mathbb{R}^{n\times p}$ and $\mathsf{B}_2(\tau) \in \mathbb{R}^{m\times p}$ are matrix-valued functions of bounded variation [49] that are right-continuous over $\mathcal{J}$. By (15), the abstract operators in (14) are denoted by concrete integral expressions. Yet, the Lebesgue-Stieltjes integrals in (15) may still be too rich to be addressed analytically. A reasonable simplification is to let

$$
\begin{aligned}
\mathsf{A}(\tau) &= \sum_{i=0}^{\nu} A_i \mathit{1}(\tau + r_i) \;+ \int_{-r}^{\tau} \widetilde{A}(\theta)\mathsf{d}\theta \\
\mathsf{B}_1(\tau) &= \sum_{i=0}^{\nu} B_i \mathit{1}(\tau + r_i) \;+ \int_{-r}^{\tau} \widetilde{B}(\theta)\mathsf{d}\theta \\
\mathsf{B}_2(\tau) &= \sum_{i=0}^{\nu} \mathfrak{B}_i \mathit{1}(\tau + r_i) \;+ \int_{-r}^{\tau} \widetilde{\mathfrak{B}}(\theta)\mathsf{d}\theta \\
\mathsf{C}(\tau) &= \sum_{i=0}^{\nu} C_i \mathit{1}(\tau + r_i) \;+ \int_{-r}^{\tau} \widetilde{C}(\theta)\mathsf{d}\theta
\end{aligned}
\tag{16}
$$

for all $\tau \in \mathcal{J}$, where $\mathit{1}(\tau)$ is the step function, $r = r_\nu > \cdots > r_2 > r_1 > r_0 = 0$ with $\nu \in \mathbb{N}$, and $A_i \in \mathbb{R}^{n\times n}$, $C_i \in \mathbb{R}^{m\times n}$, $B_i \in \mathbb{R}^{n\times p}$, $\mathfrak{B}_i \in \mathbb{R}^{m\times p}$ are constant matrices, and

$$
\begin{aligned}
&\widetilde{A}(\cdot) \in \mathcal{L}^2(\mathcal{J};\mathbb{R}^{n\times n}), \quad \widetilde{C}(\cdot) \in \mathcal{L}^2(\mathcal{J};\mathbb{R}^{m\times n}), \\
&\widetilde{B}(\cdot) \in \mathcal{L}^2(\mathcal{J};\mathbb{R}^{n\times p}), \quad \widetilde{\mathfrak{B}}(\cdot) \in \mathcal{L}^2(\mathcal{J};\mathbb{R}^{m\times p})
\end{aligned}
\tag{17}
$$

are the integral kernels. With (15)–(17), system (14) becomes

$$
\begin{aligned}
\dot{\boldsymbol{x}}(t) =& \sum_{i=0}^{\nu} A_i\boldsymbol{x}(t-r_i) + \int_{-r}^{0} \widetilde{A}(\tau)\boldsymbol{x}(t+\tau)\mathsf{d}\tau \\
&+ \sum_{i=0}^{\nu} B_i\boldsymbol{u}(t-r_i) + \int_{-r}^{0} \widetilde{B}(\tau)\boldsymbol{u}(t+\tau)\mathsf{d}\tau + D_1\boldsymbol{w}(t), \\
\boldsymbol{z}(t) =& \sum_{i=0}^{\nu} C_i\boldsymbol{x}(t-r_i) + \int_{-r}^{0} \widetilde{C}(\tau)\boldsymbol{x}(t+\tau)\mathsf{d}\tau \\
&+ \sum_{i=0}^{\nu} \mathfrak{B}_i\boldsymbol{u}(t-r_i) + \int_{-r}^{0} \widetilde{\mathfrak{B}}(\tau)\boldsymbol{u}(t+\tau)\mathsf{d}\tau + D_2\boldsymbol{w}(t), \\
\mathbf{x}_{t_0}(\theta) =& \boldsymbol{x}(t_0+\theta) = \boldsymbol{\psi}(\theta), \quad \forall \theta \in \mathcal{J}
\end{aligned}
\tag{18}
$$

by the properties of Lebesgue-Stieltjes integrals, where $\boldsymbol{\psi}(\cdot) \in \mathcal{C}(\mathcal{J};\mathbb{R}^n)$ is the initial condition. Now (18) contains the terms of pointwise and distributed delays explicitly, and we are ready to use (18) as the foundation of our subsequent analysis.

**Remark 6.** The simplifications in (16) are not restrictive. With the exception of some pathological instances such as the Can-tor function [204], (16) is adequate to denote most linear delay structures encountered in practical systems.

[1] In the cases wherein $\boldsymbol{g}(\tau)$ contains polynomials of $\tau$

**Remark 7.** A variety of cybernetic systems with general DDs can be modeled by (18) such as the characterization of event-triggered mechanisms [36], networked control systems [161], and chemical reaction networks [205, Eq. (30)], etc. Also, distributed delays have widespread applications in the SIR model of epidemics [206].

The DD integrals in (18) can always be decomposed as

$$
\begin{aligned}
\int_{-r}^{0} \widetilde{A}(\tau)\boldsymbol{x}(t+\tau)\mathsf{d}\tau &= \sum_{i=1}^{\nu}\int_{\mathcal{I}_i} \widetilde{A}_i(\tau)\boldsymbol{x}(t+\tau)\mathsf{d}\tau \\
\int_{-r}^{0} \widetilde{C}(\tau)\boldsymbol{x}(t+\tau)\mathsf{d}\tau &= \sum_{i=1}^{\nu}\int_{\mathcal{I}_i} \widetilde{C}_i(\tau)\boldsymbol{x}(t+\tau)\mathsf{d}\tau \\
\int_{-r}^{0} \widetilde{B}(\tau)\boldsymbol{u}(t+\tau)\mathsf{d}\tau &= \sum_{i=1}^{\nu}\int_{\mathcal{I}_i} \widetilde{B}_i(\tau)\boldsymbol{u}(t+\tau)\mathsf{d}\tau \\
\int_{-r}^{0} \widetilde{\mathfrak{B}}(\tau)\boldsymbol{u}(t+\tau)\mathsf{d}\tau &= \sum_{i=1}^{\nu}\int_{\mathcal{I}_i} \widetilde{\mathfrak{B}}_i(\tau)\boldsymbol{u}(t+\tau)\mathsf{d}\tau
\end{aligned} \tag{19}
$$

with $\mathcal{I}_i = [-r_i, -r_{i-1}]$ and matrix-valued functions

$$
\begin{aligned}
&\widetilde{A}_i(\cdot) \in \mathscr{L}^2(\mathcal{I}_i\,;\mathbb{R}^{n\times n}), \quad \widetilde{C}_i(\cdot) \in \mathscr{L}^2(\mathcal{I}_i\,;\mathbb{R}^{m\times n}), \\
&\widetilde{B}_i(\cdot) \in \mathscr{L}^2(\mathcal{I}_i\,;\mathbb{R}^{n\times p}), \quad \widetilde{\mathfrak{B}}_i(\cdot) \in \mathscr{L}^2(\mathcal{I}_i\,;\mathbb{R}^{m\times p})
\end{aligned} \tag{20}
$$

for all $i \in \mathbb{N}_\nu$ by the properties of integrals. This partition of $[-r, 0]$ into $\mathcal{I}_i$ in (19)–(20) is introduced for employing the KF approach in the next section.

### 3.2. *Formulating Closed-Loop System Equation via KSD*

Having introduced the KSD concept in Proposition 1 and the open-loop system in Section 3.1, we formulate the dissipative state feedback control (DSFC) problem for (18) using a static full state feedback controller, where Lemma 3 in Appendix A will be extensively employed in this paper.

First, applying the KSD concept in Proposition 1 to all matrix-valued functions in (20) with indices $\hbar = 4$, $n_1 = n_2 = n$, $n_3 = n_4 = m$, $m_1 = m_3 = n$, $m_2 = m_4 = p$, and letting $\mathfrak{F}_{i,1}(\tau) = \widetilde{A}_i(\tau)$, $\mathfrak{F}_{i,2}(\tau) = \widetilde{B}_i(\tau)$, $\mathfrak{F}_{i,3}(\tau) = \widetilde{C}_i(\tau)$, $\mathfrak{F}_{i,4}(\tau) = \widetilde{\mathfrak{B}}_i(\tau)$ gives us the decompositions

$$
\begin{aligned}
&\widetilde{A}_i(\tau) = \widehat{A}_i\left(\boldsymbol{g}_i(\tau)\otimes I_n\right), \quad \widetilde{B}_i(\tau) = \widehat{B}_i\left(\boldsymbol{g}_i(\tau)\otimes I_p\right) \\
&\widetilde{C}_i(\tau) = \widehat{C}_i\left(\boldsymbol{g}_i(\tau)\otimes I_n\right), \quad \widetilde{\mathfrak{B}}_i(\tau) = \widehat{\mathfrak{B}}_i\left(\boldsymbol{g}_i(\tau)\otimes I_p\right)
\end{aligned} \tag{21}
$$

where we defined constant matrices $\gimel_{i,1} = \widehat{A}_i \in \mathbb{R}^{n\times\kappa_i n}$, $\gimel_{i,2} = \widehat{B}_i \in \mathbb{R}^{n\times\kappa_i p}$, $\gimel_{i,3} = \widehat{C}_i \in \mathbb{R}^{m\times\kappa_i n}$, $\gimel_{i,4} = \widehat{\mathfrak{B}}_i \in \mathbb{R}^{m\times\kappa_i p}$, and $\boldsymbol{g}_i(\tau) = \begin{bmatrix}\boldsymbol{\phi}_i^\top(\tau) & \boldsymbol{\varphi}_i^\top(\tau) & \boldsymbol{f}_i^\top(\tau)\end{bmatrix}^\top$ contains appropriate functions satisfying (6b)–(6c).

Suppose all components of state $\boldsymbol{x}(t)$ can be measured for feedback. We apply static controller $\boldsymbol{u}(t) = K\boldsymbol{x}(t)$, $K \in \mathbb{R}^{p\times n}$ to (18), which gives the CLS

$$
\begin{aligned}
\dot{\boldsymbol{x}}(t) &= D_1\boldsymbol{w}(t) + \sum_{i=0}^{\nu}\left(A_i + B_iK\right)\boldsymbol{x}(t-r_i) \\
&\quad + \sum_{i=1}^{\nu}\int_{\mathcal{I}_i}\left[\widetilde{A}_i(\tau) + \widetilde{B}_i(\tau)K\right]\boldsymbol{x}(t+\tau)\mathsf{d}\tau, \\
\boldsymbol{z}(t) &= D_2\boldsymbol{w}(t) + \sum_{i=0}^{\nu}\left(C_i + \mathfrak{B}_iK\right)\boldsymbol{x}(t-r_i) \\
&\quad + \sum_{i=1}^{\nu}\int_{\mathcal{I}_i}\left[\widetilde{C}_i(\tau) + \widetilde{\mathfrak{B}}_i(\tau)K\right]\boldsymbol{x}(t+\tau)\mathsf{d}\tau, \\
\boldsymbol{x}_{t_0}(\theta) &= \boldsymbol{x}(t_0+\theta) = \boldsymbol{\psi}(\theta), \quad \forall\theta\in\mathscr{J}
\end{aligned} \tag{22}
$$

in light of (19).

Given the KSD representation of DDs in (21), both $\widetilde{B}_i(\tau)K$ and $\widetilde{\mathfrak{B}}_i(\tau)K$ in (21) and (22) involve the term $\left(\boldsymbol{g}_i(\tau)\otimes I_p\right)K$ which satisfies identity

$$
\begin{aligned}
&\left(\boldsymbol{g}_i(\tau)\otimes I_p\right)K = \left(\boldsymbol{g}_i(\tau)\otimes I_p\right)\left(I_1\otimes K\right) \\
&\qquad = I_{\kappa_i}\boldsymbol{g}_i(\tau)\otimes KI_n = \left(I_{\kappa_i}\otimes K\right)\left(\boldsymbol{g}_i(\tau)\otimes I_n\right).
\end{aligned} \tag{23}
$$

Now let $\boldsymbol{\chi}(t,\theta) = [\boldsymbol{x}(t + \acute{r}_i\theta - r_{i-1})]_{i=1}^{\nu} \in \mathbb{R}^{\nu n}$ with $\theta \in [-1, 0]$ and $\acute{r}_i = r_i - r_{i-1}$, inspired by the state variable $z(t,\tau)$ in [207]. By replacing all DDs in (22) with the decompositions in (21) and then employing (23) and Lemma 3 with $\boldsymbol{\chi}(t,\theta)$, we find that for all $i \in \mathbb{N}_\nu$

$$
\begin{aligned}
\widetilde{A}_i(\tau) + \widetilde{B}_i(\tau)K &= \left[\widehat{A}_i + \widehat{B}_i\left(I_{\kappa_i}\otimes K\right)\right]G_i(\tau), \\
\widetilde{C}_i(\tau) + \widetilde{\mathfrak{B}}_i(\tau)K &= \left[\widehat{C}_i + \widehat{\mathfrak{B}}_i\left(I_{\kappa_i}\otimes K\right)\right]G_i(\tau),
\end{aligned} \tag{24a}
$$

$$
\sum_{i=1}^{\nu}(A_i + B_iK)\boldsymbol{x}(t-r_i) = \left[\!\left[(A_i + B_iK)\right]\!\right]_{i=1}^{\nu}\boldsymbol{\chi}(t,-1), \tag{24b}
$$

$$
\sum_{i=1}^{\nu}(C_i + \mathfrak{B}_iK)\boldsymbol{x}(t-r_i) = \left[\!\left[(C_i + \mathfrak{B}_iK)\right]\!\right]_{i=1}^{\nu}\boldsymbol{\chi}(t,-1), \tag{24c}
$$

where $G_i(\tau) = \left(\boldsymbol{g}_i(\tau)\otimes I_n\right)$ with $\boldsymbol{g}_i(\cdot)$ in (6c). With (24b)–(24a), we can reformulate the CLS in (22) as

$$
\begin{aligned}
&\dot{\boldsymbol{x}}(t) = \left(A_0 + B_0K\right)\boldsymbol{x}(t) + \left[\!\left[(A_i + B_iK)\right]\!\right]_{i=1}^{\nu}\boldsymbol{\chi}(t,-1) \\
&\;+ \sum_{i=1}^{\nu}\int_{\mathcal{I}_i}\left[\widehat{A}_i + \widehat{B}_i\left(I_{\kappa_i}\otimes K\right)\right]G_i(\tau)\boldsymbol{x}(t+\tau)\mathsf{d}\tau + D_1\boldsymbol{w}(t), \\
&\boldsymbol{z}(t) = \left(C_0 + \mathfrak{B}_0K\right)\boldsymbol{x}(t) + \left[\!\left[(C_i + \mathfrak{B}_iK)\right]\!\right]_{i=1}^{\nu}\boldsymbol{\chi}(t,-1) \\
&\;+ \sum_{i=1}^{\nu}\int_{\mathcal{I}_i}\left[\widehat{C}_i + \widehat{\mathfrak{B}}_i\left(I_{\kappa_i}\otimes K\right)\right]G_i(\tau)\boldsymbol{x}(t+\tau)\mathsf{d}\tau + D_2\boldsymbol{w}(t), \\
&\boldsymbol{x}_{t_0}(\theta) = \boldsymbol{x}(t_0+\theta) = \boldsymbol{\psi}(\theta), \quad \forall\theta\in\mathscr{J}
\end{aligned} \tag{25}
$$

where all DDs in (22) are expressed via $\boldsymbol{g}_i(\cdot)$.

Given the expressions of $\boldsymbol{g}_i(\tau)$ in (10) in terms of $\boldsymbol{h}_i(\tau)$ and $\boldsymbol{\varepsilon}_i(\tau)$ in light of the approximation scheme in (8)–(10), the identities in (23) can be further expanded as

$$
\begin{aligned}
&\left(\boldsymbol{g}_i(\tau)\otimes I_p\right)K = \left(I_{\kappa_i}\otimes K\right)\left(\boldsymbol{g}_i(\tau)\otimes I_n\right) \\
&= \left(I_{\kappa_i}\otimes K\right)\left[\left(\widehat{\Gamma}_i\boldsymbol{h}_i(\tau) + \widetilde{I}_i\boldsymbol{\varepsilon}_i(\tau)\right)\otimes I_n\right]
\end{aligned}
$$

$$
\begin{aligned}
&= \left(I_{\kappa_i} \otimes K\right)\left(\widehat{\Gamma}_i \otimes I_n\right) H_i(\tau) + \left(I_{\kappa_i} \otimes K\right)\left(\widetilde{I}_i \otimes I_n\right) E_i(\tau) \\
&= \left(\widehat{\Gamma}_i \otimes K\right) H_i(\tau) + \left(\widetilde{I}_i \otimes K\right) E_i(\tau) \qquad (26)
\end{aligned}
$$

with

$$
\begin{aligned}
\boldsymbol{g}_i(\tau) \otimes I_n &= \left(\widehat{\Gamma}_i \boldsymbol{h}_i(\tau) + \widetilde{I}_i \boldsymbol{\varepsilon}_i(\tau)\right) \otimes I_n \\
&= \left(\widehat{\Gamma}_i \otimes I_n\right) H_i(\tau) + \left(\widetilde{I}_i \otimes I_n\right) E_i(\tau) \qquad (27)
\end{aligned}
$$

by (A.1), where $H_i(\tau) = \boldsymbol{h}_i(\tau) \otimes I_n$ and $E_i(\tau) = \boldsymbol{\varepsilon}_i(\tau) \otimes I_n$. By (26)–(27) and (A.1) with matrices $\mathfrak{H}_i \succ 0$, $\mathfrak{E}_i \succ 0$ in (8) and (9), we can rewrite all the DD integral matrices in (25) as

$$
\begin{aligned}
&\left[\widehat{A}_i + \widehat{B}_i \left(I_{\kappa_i} \otimes K\right)\right] \left(\boldsymbol{g}_i(\tau) \otimes I_n\right) \\
&= \left[\widehat{A}_i \left(\widehat{\Gamma}_i \otimes I_n\right) + \widehat{B}_i \left(\widehat{\Gamma}_i \otimes K\right)\right] H_i(\tau) \\
&\quad + \left[\widehat{A}_i \left(\widetilde{I}_i \otimes I_n\right) + \widehat{B}_i \left(\widetilde{I}_i \otimes K\right)\right] E_i(\tau) \\
&= \left[\widehat{A}_i \left(T_i \otimes I_n\right) + \widehat{B}_i \left(T_i \otimes K\right)\right] \left[\mathfrak{L}_{\boldsymbol{h}_i}^{-1} \boldsymbol{h}_i(\tau) \otimes I_n\right] \\
&+ \left[\widehat{A}_i \left(\widetilde{T}_i \otimes I_n\right) + \widehat{B}_i \left(\widetilde{T}_i \otimes K\right)\right] \left[\mathfrak{L}_{\boldsymbol{\varepsilon}_i}^{-1} \boldsymbol{\varepsilon}_i(\tau) \otimes I_n\right], \qquad (28)
\end{aligned}
$$

$$
\begin{aligned}
&\forall i \in \mathbb{N}_\nu, \left[\widehat{C}_i + \widehat{\mathfrak{B}}_i \left(I_{\kappa_i} \otimes K\right)\right] \left(\boldsymbol{g}_i(\tau) \otimes I_n\right) \\
&= \left[\widehat{C}_i \left(T_i \otimes I_n\right) + \widehat{\mathfrak{B}}_i \left(T_i \otimes K\right)\right] \left[\mathfrak{L}_{\boldsymbol{h}_i}^{-1} \boldsymbol{h}_i(\tau) \otimes I_n\right] \qquad (29) \\
&\quad + \left[\widehat{C}_i \left(\widetilde{T}_i \otimes I_n\right) + \widehat{\mathfrak{B}}_i \left(\widetilde{T}_i \otimes K\right)\right] \left[\mathfrak{L}_{\boldsymbol{\varepsilon}_i}^{-1} \boldsymbol{\varepsilon}_i(\tau) \otimes I_n\right]
\end{aligned}
$$

with matrices $T_i = \widehat{\Gamma}_i \mathfrak{L}_{\boldsymbol{h}_i} = \begin{bmatrix} \Gamma_i \mathfrak{L}_{\boldsymbol{h}_i}^{\top,-1} \\ \mathfrak{L}_{\boldsymbol{h}_i} \end{bmatrix}$ and $\widetilde{T}_i = \widetilde{I}_i \mathfrak{L}_{\boldsymbol{\varepsilon}_i}$ for all $i \in \mathbb{N}_\nu$, where $\mathfrak{L}_{\boldsymbol{h}_i} \mathfrak{L}_{\boldsymbol{h}_i}^\top = \mathfrak{H}_i$ and $\mathfrak{L}_{\boldsymbol{\varepsilon}_i} \mathfrak{L}_{\boldsymbol{\varepsilon}_i}^\top = \mathfrak{E}_i$ are the lower triangular Cholesky decompositions. Now applying (28)–(29) with the properties in (A.3b) to the DD integrals in (25) further produces

$$
\begin{aligned}
&\sum_{i=1}^{\nu} \int_{\mathcal{I}_i} \left[\widehat{A}_i + \widehat{B}_i \left(I_{\kappa_i} \otimes K\right)\right] \left(\boldsymbol{g}_i(\tau) \otimes I_n\right) \boldsymbol{x}(t+\tau) \mathsf{d}\tau \\
&= \left[\!\left[\widehat{A}_i \left(T_i \otimes I_n\right) + \widehat{B}_i \left(T_i \otimes K\right)\right]\!\right]_{i=1}^{\nu} \boldsymbol{\xi}(t) \\
&\qquad + \left[\!\left[\widehat{A}_i \left(\widetilde{T}_i \otimes I_n\right) + \widehat{B}_i \left(\widetilde{T}_i \otimes K\right)\right]\!\right]_{i=1}^{\nu} \boldsymbol{e}(t), \qquad (30a)
\end{aligned}
$$

$$
\begin{aligned}
&\sum_{i=1}^{\nu} \int_{\mathcal{I}_i} \left[\widehat{C}_i + \widehat{\mathfrak{B}}_i \left(I_{\kappa_i} \otimes K\right)\right] \left(\boldsymbol{g}_i(\tau) \otimes I_n\right) \boldsymbol{x}(t+\tau) \mathsf{d}\tau \\
&= \left[\!\left[\widehat{C}_i \left(T_i \otimes I_n\right) + \widehat{\mathfrak{B}}_i \left(T_i \otimes K\right)\right]\!\right]_{i=1}^{\nu} \boldsymbol{\xi}(t) \\
&\qquad + \left[\!\left[\widehat{C}_i \left(\widetilde{T}_i \otimes I_n\right) + \widehat{\mathfrak{B}}_i \left(\widetilde{T}_i \otimes K\right)\right]\!\right]_{i=1}^{\nu} \boldsymbol{e}(t), \qquad (30b)
\end{aligned}
$$

where

$$
\begin{aligned}
\boldsymbol{\xi}(t) &= \left[\int_{\mathcal{I}_i} \left(\mathfrak{L}_{\boldsymbol{h}_i}^{-1} \boldsymbol{h}_i(\tau) \otimes I_n\right) \boldsymbol{x}(t+\tau) \mathsf{d}\tau\right]_{i=1}^{\nu}, \\
\boldsymbol{e}(t) &= \left[\int_{\mathcal{I}_i} \left(\mathfrak{L}_{\boldsymbol{\varepsilon}_i}^{-1} \boldsymbol{\varepsilon}_i(\tau) \otimes I_n\right) \boldsymbol{x}(t+\tau) \mathsf{d}\tau\right]_{i=1}^{\nu}.
\end{aligned} \qquad (31)
$$

Finally, by utilizing the identities in (30a)–(31) with the properties in (A.1)–(A.3d) on (25), we obtain the final CLS

$$
\begin{aligned}
\dot{\boldsymbol{x}}(t) &= \left(\mathbf{A} + \mathbf{B}_1 \left[\left(I_\beta \otimes K\right) \oplus \mathsf{O}_q\right]\right) \boldsymbol{\vartheta}(t), && \widetilde{\forall} t \geq t_0 \\
\boldsymbol{z}(t) &= \left(\mathbf{C} + \mathbf{B}_2 \left[\left(I_\beta \otimes K\right) \oplus \mathsf{O}_q\right]\right) \boldsymbol{\vartheta}(t), && \qquad (32) \\
\mathbf{x}_{t_0}(\theta) &= \boldsymbol{x}(t_0 + \theta) = \boldsymbol{\psi}(\theta), && \forall \theta \in \mathcal{J}
\end{aligned}
$$

with $t_0$ and $\boldsymbol{\psi}(\cdot)$ in (18), where $\beta = 1+\nu+\kappa$ with $\kappa = \sum_{i=1}^{\nu} \kappa_i$ and $\kappa_i = d_i + \delta_i + \mu_i$ in Proposition 1, and

$$
\begin{aligned}
\mathbf{A} &= \left[\left[\!\left[A_i\right]\!\right]_{i=0}^{\nu} \quad \left[\!\left[\widehat{A}_i \left(T_i \otimes I_n\right)\right]\!\right]_{i=1}^{\nu} \cdots \right. \\
&\qquad \left. \cdots \left[\!\left[\widehat{A}_i \left(\widetilde{T}_i \otimes I_n\right)\right]\!\right]_{i=1}^{\nu} \quad D_1\right], \qquad (33a) \\
\mathbf{B}_1 &= \left[\left[\!\left[B_i\right]\!\right]_{i=0}^{\nu} \left[\!\left[\widehat{B}_i \left(T_i \otimes I_p\right)\right]\!\right]_{i=1}^{\nu} \cdots \right. \\
&\qquad \left. \cdots \left[\!\left[\widehat{B}_i \left(\widetilde{T}_i \otimes I_p\right)\right]\!\right]_{i=1}^{\nu} \quad \mathsf{O}_{n,q}\right], \qquad (33b) \\
\mathbf{C} &= \left[\left[\!\left[C_i\right]\!\right]_{i=0}^{\nu} \quad \left[\!\left[\widehat{C}_i \left(T_i \otimes I_n\right)\right]\!\right]_{i=1}^{\nu} \cdots \right. \\
&\qquad \left. \cdots \left[\!\left[\widehat{C}_i \left(\widetilde{T}_i \otimes I_n\right)\right]\!\right]_{i=1}^{\nu} \quad D_2\right], \qquad (33c) \\
\mathbf{B}_2 &= \left[\left[\!\left[\mathfrak{B}_i\right]\!\right]_{i=0}^{\nu} \quad \left[\!\left[\widehat{\mathfrak{B}}_i \left(T_i \otimes I_p\right)\right]\!\right]_{i=1}^{\nu} \cdots \right. \\
&\qquad \left. \cdots \left[\!\left[\widehat{\mathfrak{B}}_i \left(\widetilde{T}_i \otimes I_p\right)\right]\!\right]_{i=1}^{\nu} \quad \mathsf{O}_{m,q}\right], \qquad (33d) \\
\boldsymbol{\omega}(t) &= \left[\boldsymbol{x}^\top(t) \quad \boldsymbol{\chi}^\top(t,-1) \quad \boldsymbol{\xi}^\top(t)\right]^\top, \qquad (33e) \\
\boldsymbol{\vartheta}(t) &= \left[\boldsymbol{\omega}^\top(t) \quad \boldsymbol{e}^\top(t) \quad \boldsymbol{w}^\top(t)\right]^\top. \qquad (33f)
\end{aligned}
$$

Note that $\boldsymbol{\vartheta}(t)$ is structured so that (25) can be expressed as the products of the matrices in (33a)–(33d) and $\boldsymbol{\vartheta}(t)$.

## 4. Dissipative State Feedback Controller Design

The following criteria are presented to characterize the stability and dissipativity of system (32) with an SRF.

### *4.1. Criteria for Stability and Dissipativity*

**Lemma 1.** *The origin of the CLS in (32) with $\boldsymbol{w}(t) \equiv \mathbf{0}_q$ is exponentially stable if there exist $\epsilon_1, \epsilon_2, \epsilon_3 > 0$ and a differentiable functional $\mathsf{v} : \mathcal{C}(\mathcal{J}; \mathbb{R}^n) \to \mathbb{R}$ such that*

$$
\epsilon_1 \left\|\boldsymbol{\psi}(0)\right\|^2 \leq \mathsf{v}(\boldsymbol{\psi}(\cdot)) \leq \epsilon_2 \left\|\boldsymbol{\psi}(\cdot)\right\|_\infty^2, \qquad (34a)
$$

$$
\widetilde{\forall} t \geq t_0, \ \tfrac{\mathsf{d}}{\mathsf{d}t} \mathsf{v}(\mathbf{x}_t(\cdot)) \leq -\epsilon_3 \left\|\boldsymbol{x}(t)\right\|^2 \qquad (34b)
$$

*for any $\boldsymbol{\psi}(\cdot) \in \mathcal{C}(\mathcal{J}; \mathbb{R}^n)$ in (32), where $\left\|\boldsymbol{\psi}(\cdot)\right\|_\infty^2 = \sup_{-r \leq \tau \leq 0} \left\|\boldsymbol{\psi}(\tau)\right\|^2$. Furthermore, notation $\mathbf{x}_t(\cdot)$ is defined by the relation $\forall t \geq t_0, \forall \theta \in \mathcal{J}, \mathbf{x}_t(\theta) = \boldsymbol{x}(t+\theta)$ in which $\boldsymbol{x}(\cdot)$ satisfies the FDE in (32) with $\boldsymbol{w}(t) \equiv \mathbf{0}_q$.*

*Proof.* See [164, Lemma 2]. ■

**Remark 8.** The FDE in (32) is defined in the extended sense [3, page 58] w.r.t the Lebesgue measures. Thus, Lemma 1 is **not** a special case of the classical Krasovskiĭ stability criteria [3, Section 5.2] that were specifically developed for traditional FDEs with classical derivatives and piecewise continuous right-hand sides.

**Definition 1.** System (32) with an SRF $\mathsf{s}(\boldsymbol{z}(t), \boldsymbol{w}(t))$ is dissipative if there exists $\mathsf{v}(\bullet) \in \mathcal{C}^1\big(\mathcal{C}(\mathscr{J}; \mathbb{R}^n_{\geq 0}); \mathbb{R}\big)$ such that

$$\tfrac{\mathsf{d}}{\mathsf{d}t}\mathsf{v}(\boldsymbol{\mathsf{x}}_t(\cdot)) - \mathsf{s}(\boldsymbol{z}(t), \boldsymbol{w}(t)) \leq 0, \quad \breve{\forall} t \geq t_0 \tag{35}$$

with $t_0 \in \mathbb{R}$, $\boldsymbol{z}(t)$ and $\boldsymbol{w}(t)$ in (32). Moreover, $\boldsymbol{\mathsf{x}}_t(\cdot)$ in (35) is defined by the relation $\forall t \geq t_0$, $\forall \theta \in \mathscr{J}$, $\boldsymbol{\mathsf{x}}_t(\theta) = \boldsymbol{x}(t+\theta)$, where $\boldsymbol{x}(t)$ satisfies the equations in (32).

If the dissipative inequality in (35) holds, it follows that

$$\forall t \geq t_0, \quad \mathsf{v}(\boldsymbol{\mathsf{x}}_t(\cdot)) - \mathsf{v}(\boldsymbol{\mathsf{x}}_{t_0}(\cdot)) \leq \int_{t_0}^{t} \mathsf{s}(\boldsymbol{z}(\tau), \boldsymbol{w}(\tau))\mathsf{d}\tau, \tag{36}$$

which is the original definition of dissipativity in [208], since $\dot{\mathsf{v}}(\boldsymbol{\mathsf{x}}_t(\cdot))$ is integrable.

To characterize and enforce dissipativity for (32), we employ a quadratic supply rate function

$$\begin{gathered}\mathsf{s}(\boldsymbol{z}(t), \boldsymbol{w}(t)) = \begin{bmatrix} \boldsymbol{z}(t) \\ \boldsymbol{w}(t) \end{bmatrix}^\top \begin{bmatrix} \widetilde{J}^\top J_1^{-1} \widetilde{J} & J_2 \\ * & J_3 \end{bmatrix} \begin{bmatrix} \boldsymbol{z}(t) \\ \boldsymbol{w}(t) \end{bmatrix}, \\ \widetilde{J}^\top J_1^{-1} \widetilde{J} \preceq 0, \qquad \widetilde{J} \in \mathbb{R}^{m\times m}, \ \mathbb{S}^m \ni J_1^{-1} \prec 0, \\ J_2 \in \mathbb{R}^{m\times q}, \quad J_3 \in \mathbb{S}^q \end{gathered} \tag{37}$$

based on [209, Eq. (9)] with minor modifications (signs, matrices). The SRF in (37) can describe multiple performance criteria:

- $\mathcal{L}^2$ gain performance: $J_1 = -\gamma I_m$, $\widetilde{J} = I_m$, $J_2 = \mathsf{O}_{m,q}$, $J_3 = \gamma I_q$ with $\gamma > 0$
- Strict Passivity: $J_1 \prec 0$, $\widetilde{J} = \mathsf{O}_m$, $J_2 = I_m$, $J_3 = \mathsf{O}_m$
- Other sector constraints enumerated in [210, Table 1].

### 4.2. Main Results on DSFC Design via SDP

The principal theorem on the dissipative state feedback control problem is proposed as follows, where its proof is placed in Appendix C to ensure coherent presentations.

**Theorem 1.** *Let $\boldsymbol{g}_i(\cdot), M_i$ and $\Gamma_i, \mathfrak{H}_i, \mathfrak{E}_i$ in (8)–(10) and $\widehat{A}_i, \widehat{C}_i, \widehat{B}_i, \widehat{\mathfrak{B}}_i$ in (21) be given for utilizing the KSD in Proposition 1 with (21). Then the CLS in (32) with SRF (37) is dissipative, and the trivial solution to (32) with $\boldsymbol{w}(t) \equiv \mathbf{0}_q$ is exponentially stable if there exist $K \in \mathbb{R}^{p\times n}$, $P_1 \in \mathbb{S}^n$, $P_2 \in \mathbb{R}^{n\times dn}$, $P_3 \in \mathbb{S}^{dn}$ and $Q_i, R_i \in \mathbb{S}^n$, $i \in \mathbb{N}_\nu$ such that*

$$\begin{bmatrix} P_1 & P_2 \\ * & P_3 \end{bmatrix} + \left[\mathsf{O}_n \oplus \left(\mathop{\mathsf{diag}}_{i=1}^{\nu} I_{d_i} \otimes Q_i\right)\right] \succ 0 \tag{38a}$$

$$\mathbf{Q} = \mathop{\mathsf{diag}}_{i=1}^{\nu} Q_i \succ 0, \quad \mathbf{R} = \mathop{\mathsf{diag}}_{i=1}^{\nu} R_i \succ 0 \tag{38b}$$

$$\begin{bmatrix} \boldsymbol{\Psi} & \boldsymbol{\Sigma}^\top \widetilde{J}^\top \\ * & J_1 \end{bmatrix} = \mathsf{Sy}\left[\mathbf{P}^\top \boldsymbol{\Pi}\right] + \boldsymbol{\Phi}(P_2, P_3, \mathbf{Q}, \mathbf{R}, \boldsymbol{\Sigma}) \prec 0 \tag{38c}$$

*where $\boldsymbol{\Sigma} = \mathbf{C} + \mathbf{B}_2[(I_\beta \otimes K) \oplus \mathsf{O}_q]$ with $\mathbf{C}, \mathbf{B}_2$ in (33c)–(33d), and*

$$\boldsymbol{\Psi} = \mathsf{Sy}\left(S^\top \begin{bmatrix} P_1 & P_2 \\ * & P_3 \end{bmatrix} \begin{bmatrix} \boldsymbol{\Omega} \\ \mathbf{M} \otimes I_n \quad \mathsf{O}_{dn,(\mu n+q)} \end{bmatrix} - \begin{bmatrix} \mathsf{O}_{(\beta n),m} \\ J_2^\top \end{bmatrix} \boldsymbol{\Sigma}\right) + \Xi(\mathbf{Q}, \mathbf{R}), \tag{39a}$$

$$S = \begin{bmatrix} I_n & \mathsf{O}_{n,\nu n} & \mathsf{O}_{n,\varkappa n} & \mathsf{O}_{n,\mu n} & \mathsf{O}_{n,q} \\ \mathsf{O}_{dn} & \mathsf{O}_{dn,\nu n} & \widehat{I} & \mathsf{O}_{dn,\mu n} & \mathsf{O}_{dn,q} \end{bmatrix}, \tag{39b}$$

$$\begin{aligned}\Xi = &\left[(\mathbf{Q} + \mathbf{R}\Lambda) \oplus \mathsf{O}_n \oplus \mathsf{O}_{\kappa n} \oplus \mathsf{O}_q\right] - \\ &\left[\mathsf{O}_n \oplus \mathbf{Q} \oplus \left(\mathop{\mathsf{diag}}_{i=1}^{\nu} I_{\varkappa_i} \otimes R_i\right) \oplus \left(\mathop{\mathsf{diag}}_{i=1}^{\nu} I_{\mu_i} \otimes R_i\right) \oplus J_3\right],\end{aligned} \tag{39c}$$

$$\widehat{I} = \left(\mathop{\mathsf{diag}}_{i=1}^{\nu} \mathfrak{L}_{\boldsymbol{f}_i}^{-1} \left[\mathsf{O}_{d_i,\delta_i} \quad I_{d_i}\right] \mathfrak{L}_{\boldsymbol{h}_i}\right) \otimes I_n, \tag{39d}$$

$$\Lambda = \mathop{\mathsf{diag}}_{i=1}^{\nu} \acute{r}_i I_n, \quad \acute{r}_i = r_i - r_{i-1}, \tag{39e}$$

$$\begin{aligned}\mathbf{M} = &\left[\mathop{\mathsf{diag}}_{i=1}^{\nu} \mathfrak{L}_{\boldsymbol{f}_i}^{-1} \boldsymbol{f}_i(-r_{i-1}) \quad \mathbf{0}_d \quad \mathsf{O}_{d,\varkappa}\right] \\ &- \left[\mathbf{0}_d \quad \mathop{\mathsf{diag}}_{i=1}^{\nu} \mathfrak{L}_{\boldsymbol{f}_i}^{-1} \boldsymbol{f}_i(-r_i) \quad \mathop{\mathsf{diag}}_{i=1}^{\nu} \mathfrak{L}_{\boldsymbol{f}_i}^{-1} M_i \mathfrak{L}_{\boldsymbol{h}_i}\right]\end{aligned} \tag{39f}$$

*with $d = \sum_{i=1}^{\nu} d_i$, $\varkappa = \sum_{i=1}^{\nu} \varkappa_i = \sum_{i=1}^{\nu} d_i + \delta_i$, $\mu = \sum_{i=1}^{\nu} \mu_i$ and $M_i$ in Proposition 1. Moreover, $\boldsymbol{\Omega} := \mathbf{A} + \mathbf{B}_1[(I_\beta \otimes K) \oplus \mathsf{O}_q]$ with $\mathbf{A}, \mathbf{B}_1$ in (33a)–(33b), and $\mathfrak{L}_{\boldsymbol{f}_i}$ is the Cholesky decomposition of $\mathfrak{F}_i = \int_{\mathcal{I}_i} \boldsymbol{f}_i(\tau)\boldsymbol{f}_i^\top(\tau)\mathsf{d}\tau = \mathfrak{L}_{\boldsymbol{f}_i}\mathfrak{L}_{\boldsymbol{f}_i}^\top \succ 0$, and*

$$\mathbf{P} = \left[P_1 \quad \mathsf{O}_{n,\nu n} \quad P_2\widehat{I} \quad \mathsf{O}_{n,(\mu n+q+m)}\right], \quad \boldsymbol{\Pi} = \left[\boldsymbol{\Omega} \quad \mathsf{O}_{n,m}\right], \tag{40a}$$

$$\begin{aligned}\boldsymbol{\Phi} = \mathsf{Sy}\Bigg(&\begin{bmatrix} P_2 \\ \mathsf{O}_{\nu n,dn} \\ \widehat{I}^\top P_3 \\ \mathsf{O}_{(\mu n+q+m),dn} \end{bmatrix} \left[\mathbf{M} \otimes I_n \quad \mathsf{O}_{dn,(\mu n+q+m)}\right] \\ &+ \begin{bmatrix} \mathsf{O}_{(\beta n),m} \\ -J_2^\top \\ \widetilde{J} \end{bmatrix} \left[\boldsymbol{\Sigma} \quad \mathsf{O}_m\right]\Bigg) + \left[\Xi(\mathbf{Q}, \mathbf{R}) \oplus J_1\right].\end{aligned} \tag{40b}$$

*Finally, the number of unknowns is $(0.5d^2 + d + \nu + 0.5)n^2 + (0.5d + 0.5 + \nu + p)n \in \mathcal{O}(d^2n^2)$.*

As (32) is equivalent to (25), the existence of numerical solutions to (38a)–(38c) always ensures stability and dissipativity for the original CLS (25), not an approximated system. This stands in contrast to the methods in [211, 212], where numerical solutions may not strictly imply the stability of the original TDSs but some approximated systems. Finally, a qualitative analysis of the conservatism of Theorem 1 is provided in Appendix D.

**Remark 9.** By the statements in Proposition 1, functions $\boldsymbol{g}_i(\cdot)$ must be chosen so as to be able to represent all the DD-matrix kernels in (20), while adding more functions satisfying (6b)–(6c) to $\boldsymbol{g}_i(\cdot)$ is always applicable since $\mathsf{dim}(\boldsymbol{g}_i(\tau))$ is unbounded. At the same time, the proposed framework grants complete flexibility in determining how the distributed delay kernels of $\widetilde{A}_i(\cdot), \widetilde{C}_i(\cdot)$ are positioned within $\boldsymbol{g}_i(\tau) = \left[\boldsymbol{\phi}_i^\top(\tau) \quad \boldsymbol{\varphi}_i^\top(\tau) \quad \boldsymbol{f}_i^\top(\tau)\right]^\top$. For instance, all kernels in (20) can be included by $\boldsymbol{f}_i(\cdot)$ if they are all weakly differentiable, where no approximations of functions are needed. On the other hand, some kernels in (20) can be placed in $\boldsymbol{\varphi}_i(\cdot)$ or $\boldsymbol{\phi}_i(\cdot)$ if they are not weakly differentiable. In fact, we can allocate all the kernels in (20) to $\boldsymbol{\phi}_i(\cdot)$ that are approximated by $\boldsymbol{f}_i(\cdot)$ and $\boldsymbol{\varphi}_i(\cdot)$, where $\boldsymbol{f}_i(\cdot)$ is completely independent of the kernels in (20). This property is significant since one can utilize any $\boldsymbol{f}_i(\cdot) \in \mathscr{W}^{1,2}(\mathcal{I}_i\,; \mathbb{R}^{d_i})$ for (6b) and (21) and (C.1b) to increase the feasibility of the synthesis condition. The computational complexity of Theorem 1 is of $\mathcal{O}(d^2 n^2)$, depending on $d_i = \mathsf{dim}(\boldsymbol{f}_i(\tau))$. Thus, the choice of $\boldsymbol{g}_i(\cdot)$ should balance both feasibility and computational complexity. We refer readers to Appendix D for more information on the potential strategies to reduce the conservatism of the SDP constraints in Theorem 1.

The inequality in (38c) is bilinear (nonconvex) due to the products between $K$ and $P_1, P_2$. The subsequent theorem, whose proof is detailed in Appendix E, addresses this issue by decoupling the BMI in (38c) using the Projection Lemma [169–171]. The core strategy of this type of methods [213–216] is to construct convex SDP constraints via the introduction of slack variables, while preserving the structural integrity of $P_2 \in \mathbb{R}^{n \times dn}$.

**Lemma 2** (Projection Lemma)**.** *[169–171] Given $n; p; q \in \mathbb{N}$, $\Pi \in \mathbb{S}^n, P \in \mathbb{R}^{q \times n}, Q \in \mathbb{R}^{p \times n}$, there exists $\Theta \in \mathbb{R}^{p \times q}$ such that the following two assertions are equivalent :*

$$\Pi + P^\top \Theta^\top Q + Q^\top \Theta P \prec 0, \tag{41a}$$

$$P_\perp^\top \Pi P_\perp \prec 0 \quad \& \quad Q_\perp^\top \Pi Q_\perp \prec 0, \tag{41b}$$

*where the columns of $P_\perp$, $Q_\perp$ contain bases of the null space of $P$, $Q$, respectively, satisfying $PP_\perp = \mathsf{O}$ and $QQ_\perp = \mathsf{O}$.*

**Theorem 2.** *Given scalars $\{\alpha_i\}_{i=1}^{\beta} \subset \mathbb{R}$ and $\boldsymbol{g}_i(\cdot), M_i$ and $\Gamma_i, \mathfrak{H}_i, \mathfrak{E}_i$ in (8)–(10) and $\widehat{A}_i, \widehat{C}_i, \widehat{B}_i, \widehat{\mathfrak{B}}_i$ in (21) for utilizing Proposition 1 with (21), the CLS in (32) with SRF (37) is dissipative and the origin of (32) with $\boldsymbol{w}(t) \equiv \mathbf{0}_q$ is globally exponentially stable if there exist $X, \acute{P}_1 \in \mathbb{S}^n$, $\acute{P}_2 \in \mathbb{R}^{n \times dn}$, $\acute{P}_3 \in \mathbb{S}^{dn}$ and $\acute{Q}_i, \acute{R}_i \in \mathbb{S}^n$, and $V \in \mathbb{R}^{p \times n}$ such that*

$$\begin{bmatrix} \acute{P}_1 & \acute{P}_2 \\ * & \acute{P}_3 \end{bmatrix} + \left[\mathsf{O}_n \oplus \left(\operatorname*{diag}_{i=1}^{\nu} I_{d_i} \otimes \acute{Q}_i\right)\right] \succ 0, \tag{42a}$$

$$\acute{\mathbf{Q}} = \operatorname*{diag}_{i=1}^{\nu} \acute{Q}_i \succ 0, \quad \acute{\mathbf{R}} = \operatorname*{diag}_{i=1}^{\nu} \acute{R}_i \succ 0, \tag{42b}$$

$$\mathsf{Sy}\left(\begin{bmatrix} I_n \\ [\alpha_i I_n]_{i=1}^{\beta} \\ \mathsf{O}_{(q+m),n} \end{bmatrix} \begin{bmatrix} -X & \acute{\mathbf{\Pi}} \end{bmatrix}\right) + \begin{bmatrix} \mathsf{O}_n & \acute{\mathbf{P}} \\ * & \acute{\mathbf{\Phi}} \end{bmatrix} \prec 0, \tag{42c}$$

*where* $\acute{\mathbf{P}} = \begin{bmatrix} \acute{P}_1 & \mathsf{O}_{n,\nu n} & \acute{P}_2 \widehat{I} & \mathsf{O}_{n,(\mu n+q+m)} \end{bmatrix}$ *and*

$$\acute{\mathbf{\Pi}} = \begin{bmatrix} \mathbf{A}\left[(I_\beta \otimes X) \oplus I_q\right] + \mathbf{B}_1\left[(I_\beta \otimes V) \oplus \mathsf{O}_q\right] & \mathsf{O}_{n,m} \end{bmatrix}$$

*with $\widehat{I}$ in (39d), and $\acute{\mathbf{\Phi}} := \mathbf{\Phi}\big(\acute{P}_2, \acute{P}_3, \acute{\mathbf{Q}}, \acute{\mathbf{R}}, \acute{\mathbf{\Sigma}}\big)$ and $\acute{\mathbf{\Sigma}} = \mathbf{C}\left[(I_\beta \otimes X) \oplus I_q\right] + \mathbf{B}_2\left[(I_\beta \otimes V) \oplus \mathsf{O}_q\right]$ with $\mathbf{A}, \mathbf{B}_1, \mathbf{B}_2, \mathbf{C}$ in (33a)–(33d) and $\mathbf{\Phi}$ in the form of (38c). The controller gain $K$ is calculated via $K = VX^{-1}$, where the non-singularity of $X$ is guaranteed by the $(1,1)$ block of (42c). The total number of unknowns is $(0.5d^2 + d + \nu + 1)n^2 + (0.5d + 1 + \nu + p)n \in \mathcal{O}(d^2 n^2)$.*

**Remark 10.** While the step (E.5) in the proof of Theorem 2 may introduce conservatism, $\acute{P}_1, \acute{P}_2$ in (42a) retain the same degrees of freedom as of $P_1, P_2$ in (38a). Consequently, employing Lemma 2 at step (E.4) preserves the structural integrity of the matrix parameters in (C.1a), resulting in less conservatism than that would be caused by directly letting $P_2 = [\![ b_i P_1 ]\!]_{i=1}^{d}$ with known scalars $\{b_i\}_{i=1}^{d} \subset \mathbb{R}$ to convexify the BMI in (38c). Note, the introduction of slack variable $W$ does not enhance the feasibility of Theorem 2 relative to Theorem 1.

**Remark 11.** For scalars $\{\alpha_i\}_{i=1}^{\beta} \subset \mathbb{R}$, we can assume $\alpha_i = 0$ for $i = 2, \ldots, \beta$ with an adjustable parameter $\alpha_1 \in \mathbb{R} \setminus \{0\}$. It is crucial to stress that $\alpha_1 \neq 0$ is necessary, as $\alpha_1 = 0$ renders the $A_0$ related-diagonal-block of $\acute{\mathbf{\Theta}}$ in (E.9) infeasible.

### *4.3. Inner Convex Approximation of BMI*

Although (42a)–(42c) are convex, the simplification in (E.5) may introduce conservatism compared to the original synthesis condition in Theorem 1. Thus, there is a preference for methodologies capable of solving (38c) directly. There are various numerical approaches available that can compute local solutions for general non-convex SDPs [217]. Here, we propose an offline sequential convex SDP algorithm based on the inner convex approximation strategy developed by [218, 219]. Our algorithm ensures monotonic convergence to a KKT point based on the convergence analysis by [218, Th. 3.2] for problems subject to general bilinear matrix inequalities. Each iteration in our algorithm is formulated as a convex SDP program, and the initial point of the whole algorithm can be provided by a feasible solution to the convex SDP problem in Theorem 2, effectively combining the advantages of both Theorem 1 and Theorem 2.

First of all, notice that the inequality in (38c) is nonconvex whereas (38a)–(38b) remain convex even if $K \neq \mathsf{O}_{p \times n}$. Now we can reformulate the inequality in (38c) as

$$\mathsf{Sy}\left[\mathbf{P}^\top \mathbf{\Pi}\right] + \mathbf{\Phi} = \mathsf{Sy}\left(\mathbf{P}^\top \mathbf{B}\left[(I_\beta \otimes K) \oplus \mathsf{O}_{q+m}\right]\right) + \widehat{\mathbf{\Phi}} \prec 0 \tag{43}$$

with $\mathbf{B} := \begin{bmatrix} \mathbf{B}_1 & \mathsf{O}_{n,m} \end{bmatrix}$ and $\widehat{\mathbf{\Phi}} := \mathsf{Sy}\left(\mathbf{P}^\top \begin{bmatrix} \mathbf{A} & \mathsf{O}_{n,m} \end{bmatrix}\right) + \mathbf{\Phi}$, where $\mathbf{P}$ is given by (40a), and $\mathbf{A}$ and $\mathbf{B}_1$ are given by (33a)–(33b). Note that $\widehat{\mathbf{\Phi}}$ is convex w.r.t the unknowns inside. By the conclusions of [218, Example 3], we see that function $\Delta\big(\bullet, \widetilde{\mathbf{G}}, \bullet, \widetilde{\mathbf{N}}\big)$, which is defined as

$$\mathbb{S}^\ell \ni \Delta\big(\mathbf{G}, \widetilde{\mathbf{G}}, \mathbf{N}, \widetilde{\mathbf{N}}\big) := [*]\left[Z \oplus (I_n - Z)\right]^{-1} \begin{bmatrix} \mathbf{G} - \widetilde{\mathbf{G}} \\ \mathbf{N} - \widetilde{\mathbf{N}} \end{bmatrix} + \mathsf{Sy}\left(\widetilde{\mathbf{G}}^\top \mathbf{N} + \mathbf{G}^\top \widetilde{\mathbf{N}} - \widetilde{\mathbf{G}}^\top \widetilde{\mathbf{N}}\right) + \mathbf{T} \tag{44a}$$

with $Z \oplus (I_n - Z) \succ 0$ satisfying

$$\forall \mathbf{G}; \widetilde{\mathbf{G}} \in \mathbb{R}^{n \times \ell},\ \forall \mathbf{N}; \widetilde{\mathbf{N}} \in \mathbb{R}^{n \times \ell},\quad \mathbf{T} + \mathsf{Sy}\left(\mathbf{G}^\top \mathbf{N}\right) = \Delta\big(\mathbf{G}, \mathbf{G}, \mathbf{N}, \mathbf{N}\big) \preceq \Delta\big(\mathbf{G}, \widetilde{\mathbf{G}}, \mathbf{N}, \widetilde{\mathbf{N}}\big), \tag{44b}$$

is a PSD-convex mapping [220, Eq. (6)] and a PSD-convex overestimate [218, Example 3] of the matrix-valued function

$\Delta(\mathbf{G},\mathbf{G},\mathbf{N},\mathbf{N}) = \mathbf{T} + \mathsf{Sy}\left[\mathbf{G}^\top\mathbf{N}\right]$ w.r.t the parameterization

$$\mathbf{vec}(\widetilde{\mathbf{G}}) = \mathbf{vec}(\mathbf{G}) \quad \& \quad \mathbf{vec}(\widetilde{\mathbf{N}}) = \mathbf{vec}(\mathbf{N}).$$

Now let $\mathbf{T} = \widehat{\mathbf{\Phi}}$, $\mathbf{G} = \mathbf{P}$ with $\widetilde{P}_1 \in \mathbb{S}^n$, $\widetilde{P}_2 \in \mathbb{R}^{n\times dn}$ and

$$\begin{aligned}
&\widetilde{\mathbf{G}} = \widetilde{\mathbf{P}} = \begin{bmatrix}\widetilde{P}_1 & \mathsf{O}_{n,\nu n} & \widetilde{P}_2\widehat{I} & \mathsf{O}_{n,(\mu n+q+m)}\end{bmatrix},\\
&\mathbf{N} = \mathbf{B}\mathbf{K}, \quad \mathbf{K} = \left(I_\beta \otimes K\right) \oplus \mathsf{O}_{q+m},\\
&\widetilde{\mathbf{N}} = \mathbf{B}\widetilde{\mathbf{K}}, \quad \widetilde{\mathbf{K}} = \left(I_\beta \otimes \widetilde{K}\right) \oplus \mathsf{O}_{q+m}, \quad \widetilde{K} \in \mathbb{R}^{p\times n}
\end{aligned} \tag{45}$$

for (44a) with $\ell = \beta n + q + m$ and $Z \oplus (I_n - Z) \succ 0$, where $\widehat{\mathbf{\Phi}}, \mathbf{P}, \mathbf{B}$ and $K$ are given by (43). Then we can conclude

$$\begin{aligned}
&\mathfrak{S}\left(\mathbf{P},\mathbf{P},K,K\right) = \widehat{\mathbf{\Phi}} + \mathsf{Sy}\left[\mathbf{P}^\top\mathbf{B}\left[(I_\beta\otimes K)\oplus\mathsf{O}_{q+m}\right]\right]\\
&\preceq \mathfrak{S}\left(\mathbf{P},\widetilde{\mathbf{P}},K,\widetilde{K}\right) := \widehat{\mathbf{\Phi}} + \mathsf{Sy}\left(\widetilde{\mathbf{P}}^\top\mathbf{N} + \mathbf{P}^\top\widetilde{\mathbf{N}} - \widetilde{\mathbf{P}}^\top\widetilde{\mathbf{N}}\right)\\
&+ \begin{bmatrix}\mathbf{P}^\top - \widetilde{\mathbf{P}}^\top & \mathbf{N}^\top - \widetilde{\mathbf{N}}^\top\end{bmatrix}\left[Z \oplus (I_n - Z)\right]^{-1}[*]
\end{aligned} \tag{46}$$

by (44b), where $\mathfrak{S}\left(\mathbf{P},\mathbf{P},K,K\right) = \mathsf{Sy}\left[\mathbf{P}^\top\mathbf{\Pi}\right] + \mathbf{\Phi}$ as in (43) that is implied by $\mathfrak{S}\left(\mathbf{P},\widetilde{\mathbf{P}},K,\widetilde{K}\right) \prec 0$ with (46). Thus, it is easy to see that $\mathfrak{S}\left(\mathbf{P},\widetilde{\mathbf{P}},K,\widetilde{K}\right) \prec 0$ if and only if

$$\begin{bmatrix}\widehat{\mathbf{\Phi}} + \mathsf{Sy}\left(\widetilde{\mathbf{P}}^\top\mathbf{N} + \mathbf{P}^\top\widetilde{\mathbf{N}} - \widetilde{\mathbf{P}}^\top\widetilde{\mathbf{N}}\right) & \mathbf{P}^\top - \widetilde{\mathbf{P}}^\top & \mathbf{N}^\top - \widetilde{\mathbf{N}}^\top\\ * & -Z & \mathsf{O}_n\\ * & * & Z - I_n\end{bmatrix} \prec 0 \tag{47}$$

with $\mathbf{N}, \widetilde{\mathbf{N}}$ in (45) by an application of the Schur complement with $Z \oplus (I_n - Z) \succ 0$. For this reason, (43) is inferred from (47), where (47) can then be computed by convex SDP solvers if $\widetilde{\mathbf{P}}$ and $\widetilde{K}$ are known.

By compiling all the preceding steps according to the expositions in [218], Algorithm 1 is established, where $f(\mathbf{x})$ is the convex objective function with $\mathbf{x}$ comprising all the decision variables in (47), and $\|\mathbf{x}\|_\infty := \max_i |x_i|$, and $\|X\|_F := \sqrt{\mathrm{tr}(X^\top X)}$ is the Frobenius norm. Scalars $\rho_1, \rho_2$ and $\varepsilon$ are given constants for regularization and regulating error tolerance, respectively.

## 5. Solutions to a Variant DSFC Problem

### 5.1. Designing Controllers with Delays via KSD in the Absence of Input Delays

If delays are absent from the control input $\boldsymbol{u}(\cdot)$ in (18), then we can adjust Theorems 1 and 2 and Algorithm 1 to address a different DSFC problem where the feedback can contain any number of pointwise and distributed delays.

Specifically, let $B_i = \widetilde{B}_i(\tau) = \mathsf{O}_{n,p}$, $\mathfrak{B}_i = \widetilde{\mathfrak{B}}_i(\tau) = \mathsf{O}_{m,p}$ for all $i \in \mathbb{N}_\nu$ in (18) where input delays are not present. Our goal is to construct a controller

$$\boldsymbol{u}(t) = \sum_{i=0}^{\nu} K_i\boldsymbol{x}(t - r_i) + \sum_{i=1}^{\nu}\int_{\mathcal{I}_i}\widetilde{K}_i(\tau)\boldsymbol{x}(t+\tau)\mathrm{d}\tau \tag{48}$$

to stabilize the above system with the SRF in (37), where $K_0, K_i \in \mathbb{R}^{p\times n}$ and $\widetilde{K}_i(\cdot) \in \mathcal{L}^2(\mathcal{I}_i;\mathbb{R}^{p\times n})$ for all $i \in \mathbb{N}_\nu$. This class of systems is of critical importance for the research

---

**Algorithm 1:** An iterative solution to the non-convex SDP denoted by (38a)–(38c) in Theorem 1

---

**begin**

  **choose** $\{\alpha_i\}_{i=1}^{\beta} \subset \mathbb{R}$, **solve** (42a)–(42c) **return** $K$,

  **solve** (38a)–(38c) with $K$ **return** $\mathbf{Y} = \begin{bmatrix}P_1 & P_2\end{bmatrix}$,

  **solve** (38a)–(38c) again with the above $\mathbf{Y}$ **return** $K$,

  **update** $\widetilde{\mathbf{Y}} = \begin{bmatrix}\widetilde{P}_1 & \widetilde{P}_2\end{bmatrix} \longleftarrow \mathbf{Y}$, $\widetilde{K} \longleftarrow K$,

  **solve** $\min_{\mathbf{x}} f(\mathbf{x}) + \rho_1\left\|\mathbf{Y} - \widetilde{\mathbf{Y}}\right\|_F^2 + \rho_2\left\|K - \widetilde{K}\right\|_F^2$
   subject to (38a)–(38b), (47) with (45) and the parameters in Theorem 1, **return** $\mathbf{Y}$ and $K$

  **while** $\dfrac{\left\|\mathbf{vec}\begin{bmatrix}\mathbf{Y} & K^\top\end{bmatrix} - \mathbf{vec}\begin{bmatrix}\widetilde{\mathbf{Y}} & \widetilde{K}^\top\end{bmatrix}\right\|_\infty}{\left\|\mathbf{vec}\begin{bmatrix}\widetilde{\mathbf{Y}} & \widetilde{K}^\top\end{bmatrix}\right\|_\infty + 1} \geq \varepsilon$ **do**

    **update** $\mathbb{R}^{n\times(1+d)n} \ni \widetilde{\mathbf{Y}} \longleftarrow \mathbf{Y}$, $\widetilde{K} \longleftarrow K$,

    **solve** $\min_{\mathbf{x}} f(\mathbf{x}) + \rho_1\left\|\mathbf{Y} - \widetilde{\mathbf{Y}}\right\|_F^2 + \rho_2\left\|K - \widetilde{K}\right\|_F^2$
     subject to (38a)–(38b), (47) with (45) and the parameters in Theorem 1, **return** $\mathbf{Y}$ and $K$

  **end**

**end**

---

of LTDSs, as early results by [77] have suggested that a controller with DDs can be necessary to stabilize unstable LTDSs in the absence of input delays.

Now apply Proposition 1 with $\hbar = 3$, $n_1 = n$, $n_2 = m$, $n_3 = p$, $m_1 = m_2 = m_3 = n$, and $\mathcal{F}_{i,1}(\tau) = \widetilde{A}_i(\tau)$, $\mathcal{F}_{i,2}(\tau) = \widetilde{C}_i(\tau)$, and $\mathcal{F}_{i,3}(\tau) = \widetilde{K}_i(\tau)$ in (48) and take into account the matrix kernels in (21). Then we obtain the decompositions in (21) together with

$$\widetilde{K}_i(\tau) = \mathcal{K}_i\left(\boldsymbol{g}_i(\tau)\otimes I_n\right), \quad \mathcal{K}_i \in \mathbb{R}^{p\times\kappa_i n} \tag{49}$$

for all $i \in \mathbb{N}_\nu$ and $\forall\tau \in \mathcal{I}_i$.

Now by adding (49) to the statements of Proposition 1 and taking into account the spirit of the derivations in (8)–(31), we can derive the resulting CLS

$$\begin{aligned}
&\dot{\boldsymbol{x}}(t) = \left(\mathbf{A} + B_0\boldsymbol{K}\right)\boldsymbol{\vartheta}(t), \quad \boldsymbol{z}(t) = \left(\mathbf{C} + \mathfrak{B}_0\boldsymbol{K}\right)\boldsymbol{\vartheta}(t),\\
&\quad\forall\theta \in \mathcal{J}, \quad \boldsymbol{x}(t_0+\theta) = \boldsymbol{\psi}(\theta), \quad \boldsymbol{\psi}(\cdot) \in \mathcal{C}(\mathcal{J};\mathbb{R}^n),\\
&\boldsymbol{K} = \Big[\,[\![K_i]\!]_{i=0}^{\nu} \quad [\![\mathcal{K}_i\left(T_i\otimes I_n\right)]\!]_{i=1}^{\nu}\cdots\\
&\qquad\qquad\cdots\left[\!\!\left[\mathcal{K}_i\left(\widetilde{T}_i\otimes I_n\right)\right]\!\!\right]_{i=1}^{\nu} \quad \mathsf{O}_{p,q}\Big]
\end{aligned} \tag{50}$$

with $\mathbf{A}, \mathbf{C}, \boldsymbol{\vartheta}(t)$ in (33a)–(33f) and $T_i; \widetilde{T}_i$ in (29) and the controller parameters in (48)–(49).

### 5.2. Main Results on DSFC Design in the Absence of Input Delays via SDP

Since the LTDS in (50) is also denoted using $\boldsymbol{\vartheta}(t)$ in (33f), Theorems 1 and 2 can be modified to accommodate the CLS in (50) accordingly, which results in the following two theorems.

**Theorem 3.** *Let $\boldsymbol{g}_i(\cdot), M_i$ and $\Gamma_i, \mathfrak{H}_i, \mathfrak{E}_i$ in (8)–(10), (49) and $\widehat{A}_i, \widehat{C}_i, \widehat{B}_i, \widehat{\mathfrak{B}}_i$ in (21) be given for utilizing the KSD in Proposition 1 with (21). Then the CLS in (50) with SRF (37) is dissipative and the trivial solution to (50) with $\boldsymbol{w}(t) \equiv \mathbf{0}_q$ is exponentially stable if there exist $P_1 \in \mathbb{S}^n$, $P_2 \in \mathbb{R}^{n \times dn}$, $P_3 \in \mathbb{S}^{dn}$ with $d = \sum_{i=1}^{\nu} d_i$, and $Q_i, R_i \in \mathbb{S}^n$, $K_0, K_i \in \mathbb{R}^{p \times n}$, $\mathcal{K}_i \in \mathbb{R}^{p \times \kappa_i n}$, $i \in \mathbb{N}_\nu$ such that (38a)–(38c) hold with $\boldsymbol{\Omega} = \mathbf{A} + B_0 \boldsymbol{K}$, $\boldsymbol{\Sigma} = \mathbf{C} + \mathfrak{B}_0 \boldsymbol{K}$ and $\mathbf{A}, \mathbf{C}$ in (33a)–(33c), $\boldsymbol{K}$ in (50) and $B_0, \mathfrak{B}_0$ in (18). Finally, the number of unknowns is $(0.5d^2 + d + \nu + 0.5)n^2 + (0.5d + 0.5 + \nu + \beta p)n$.*

*Proof.* The theorem is proven immediately by using $\boldsymbol{\Omega} = \mathbf{A} + B_0 \boldsymbol{K}$ and $\boldsymbol{\Sigma} = \mathbf{C} + \mathfrak{B}_0 \boldsymbol{K}$ in (38c). ■

**Theorem 4.** *Let $\boldsymbol{g}_i(\cdot), M_i$ and $\Gamma_i, \mathfrak{H}_i, \mathfrak{E}_i$ in (8)–(10) for the KSD in Section 3.2 with (49) be given with known scalars $\{\alpha_i\}_{i=1}^{\beta}$. Then the CLS in (50) with SRF (37) is dissipative and the trivial solution to (50) with $\boldsymbol{w}(t) \equiv \mathbf{0}_q$ is exponentially stable if there exist $\acute{P}_1, X \in \mathbb{S}^n, \acute{P}_2 \in \mathbb{R}^{n \times dn}, \acute{P}_3 \in \mathbb{S}^{dn}$ and $\acute{Q}_i, \acute{R}_i \in \mathbb{S}^n$ and $V_0, V_i \in \mathbb{R}^{p \times n}$, $\mathcal{V}_i \in \mathbb{R}^{p \times \kappa_i n}$, $i \in \mathbb{N}_\nu$ such that (42a)–(42c) hold with*

$$\acute{\boldsymbol{\Pi}} = \left[\mathbf{A}\left[(I_\beta \otimes X) \oplus I_q\right] + B_0 \boldsymbol{V} \quad \mathsf{O}_{n,m}\right],$$
$$\acute{\boldsymbol{\Sigma}} = \mathbf{C}\left[(I_\beta \otimes X) \oplus I_q\right] + \mathfrak{B}_0 \boldsymbol{V},$$

*where $\mathbf{A}, \mathbf{C}$ are given by (33a)–(33c) with $B_0, \mathfrak{B}_0$ in (18),*

$$\text{and } \boldsymbol{V} = \Big[ \llbracket V_i \rrbracket_{i=0}^{\nu} \quad \llbracket \mathcal{V}_i (T_i \otimes I_n) \rrbracket_{i=1}^{\nu} \cdots \\ \cdots \left\llbracket \mathcal{V}_i \left(\widetilde{T}_i \otimes I_n\right) \right\rrbracket_{i=1}^{\nu} \quad \mathsf{O}_{p,q} \Big] \tag{51}$$

*with $T_i, \widetilde{T}_i$ in (29). Moreover, the controller gains are calculated via $K_0 = V_0 X^{-1}$ and $K_i = V_i X^{-1}$ and $\mathcal{K}_i = \mathcal{V}_i \left(I_{\kappa_i} \otimes X^{-1}\right)$ for $i \in \mathbb{N}_\nu$. Finally, the number of unknowns is $(0.5d^2 + d + \nu + 1)n^2 + (0.5d + \nu + 1 + \beta p)n$.*

*Proof.* The proof is achieved by adjusting the proof of Theorem 2 at step (E.11), which should be modified as

$$\acute{\boldsymbol{\Pi}} = \left[\mathbf{A}\left[(I_\beta \otimes X) \oplus I_q\right] + B_0 \boldsymbol{K}\left[(I_\beta \otimes X) \oplus I_q\right] \quad \mathsf{O}_{n,m}\right] \\ = \left[\mathbf{A}\left[(I_\beta \otimes X) \oplus I_q\right] + B_0 \boldsymbol{V} \quad \mathsf{O}_{n,m}\right]$$

with $\boldsymbol{V}$ in (51), where $V_0 = K_0 X$, $V_i = K_i X$ and $\mathcal{V}_i = \mathcal{K}_i (I_{\kappa_i} \otimes X)$ for all $i \in \mathbb{N}_\nu$. Also, we note that equality $\boldsymbol{K}\left[(I_\beta \otimes X) \oplus I_q\right] = \boldsymbol{V}$ with $\boldsymbol{K}, \boldsymbol{V}$ in (50)–(51) can be established by (A.1). ■

The BMI in Theorem 3 can also be computed by a modified version of Algorithm 1, which gives Algorithm 2. This can be done by using the substitutions $\mathbf{N} = B_0 \left[\boldsymbol{K} \quad \mathsf{O}_{p,m}\right]$, $\widetilde{\mathbf{N}} = B_0 \left[\widetilde{\boldsymbol{K}} \quad \mathsf{O}_{p,m}\right]$, $K \leftarrow \boldsymbol{\mathfrak{K}}$, $\widetilde{K} \leftarrow \widetilde{\boldsymbol{\mathfrak{K}}}$ for (47), with the parameters in Theorem 3, and $\boldsymbol{K}$ in (50) and $\widetilde{K}_0, \widetilde{K}_i \in \mathbb{R}^{p \times n}$, $\widetilde{\mathcal{K}}_i \in \mathbb{R}^{p \times \kappa_i n}$ with

$$\widetilde{\boldsymbol{K}} = \Big[ \left\llbracket \widetilde{K}_i \right\rrbracket_{i=0}^{\nu} \quad \left\llbracket \widetilde{\mathcal{K}}_i (T_i \otimes I_n) \right\rrbracket_{i=1}^{\nu} \cdots \\ \cdots \left\llbracket \widetilde{\mathcal{K}}_i \left(\widetilde{T}_i \otimes I_n\right) \right\rrbracket_{i=1}^{\nu} \quad \mathsf{O}_{p,q} \Big], \tag{52}$$
$$\boldsymbol{\mathfrak{K}} = \left[ \llbracket K_i \rrbracket_{i=0}^{\nu} \ \llbracket \mathcal{K}_i \rrbracket_{i=1}^{\nu} \right], \quad \widetilde{\boldsymbol{\mathfrak{K}}} = \left[ \left\llbracket \widetilde{K}_i \right\rrbracket_{i=0}^{\nu} \ \left\llbracket \widetilde{\mathcal{K}}_i \right\rrbracket_{i=1}^{\nu} \right].$$

---

**Algorithm 2:** An iterative solution to the non-convex SDP in Theorem 3

---

**begin**

  **choose** $\{\alpha_i\}_{i=1}^{\beta} \subset \mathbb{R}$, **utilize** Theorem 4 **compute** $\boldsymbol{\mathfrak{K}}$,
  **utilize** Theorem 3 with $\boldsymbol{\mathfrak{K}}$ **compute** $\mathbf{Y} = \left[P_1 \quad P_2\right]$
  **utilize** Theorem 3 with $\mathbf{Y}$ **compute** $\boldsymbol{\mathfrak{K}}$ again

  **update** $\widetilde{\mathbf{Y}} = \left[\widetilde{P}_1 \quad \widetilde{P}_2\right] \leftarrow \mathbf{Y}$, $\widetilde{\boldsymbol{\mathfrak{K}}} \longleftarrow \boldsymbol{\mathfrak{K}}$

  **solve** $\min_{\mathbf{x}} f(\mathbf{x}) + \rho_1 \left\|\mathbf{Y} - \widetilde{\mathbf{Y}}\right\|_F^2 + \rho_2 \left\|\boldsymbol{\mathfrak{K}} - \widetilde{\boldsymbol{\mathfrak{K}}}\right\|_F^2$ subject to (38a)–(38b), (47) with (45) and $\mathbf{N} = B_0 \left[\boldsymbol{K} \quad \mathsf{O}_{p,m}\right]$, $\widetilde{\mathbf{N}} = B_0 \left[\widetilde{\boldsymbol{K}} \quad \mathsf{O}_{p,m}\right]$, with the parameters in (52), **return** $\mathbf{Y}$ and $\boldsymbol{\mathfrak{K}}$

  **while** $\dfrac{\left\| \begin{bmatrix} \mathbf{vec}(\mathbf{Y}) \\ \mathbf{vec}(\boldsymbol{\mathfrak{K}}) \end{bmatrix} - \begin{bmatrix} \mathbf{vec}(\widetilde{\mathbf{Y}}) \\ \mathbf{vec}(\widetilde{\boldsymbol{\mathfrak{K}}}) \end{bmatrix} \right\|_\infty}{\left\| \begin{bmatrix} \mathbf{vec}(\widetilde{\mathbf{Y}}) \\ \mathbf{vec}(\widetilde{\boldsymbol{\mathfrak{K}}}) \end{bmatrix} \right\|_\infty + 1} \geq \varepsilon$ **do**

    **update** $\widetilde{\mathbf{Y}} \longleftarrow \mathbf{Y}$, $\mathbb{R}^{p \times \beta n} \ni \widetilde{\boldsymbol{\mathfrak{K}}} \longleftarrow \boldsymbol{\mathfrak{K}}$

    **solve** again the SDP problem in the previous step involving (52), **return** $\mathbf{Y}$ and $\boldsymbol{\mathfrak{K}}$

  **end**

**end**

---

Finally, the diagrams below offer an intuitive explanation of how the proposed theorems and iterative algorithms relate to one another, illustrating their combined application as a unified solution to the DSFC problem for the time-delay systems in (18) and (50).

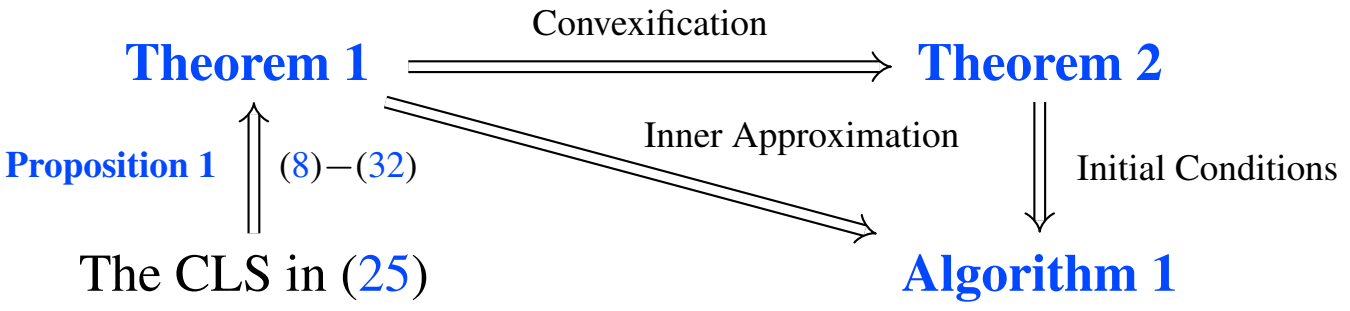

Schematic Diagram of DSFC design for the CLS in (32).

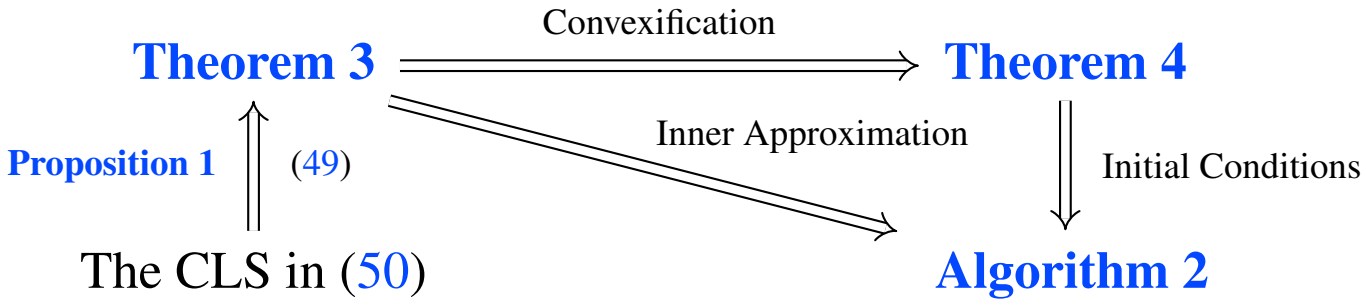

Schematic Diagram of DSFC design for the CLS in (50).

**Remark 12.** From a mathematical perspective, it should be noted that we have not applied the controller in (48) to the sys-tem in (18). The presence of input delays in (18) will certainly lead to the introduction of composite integrals in the form of $\int_{-r}^{0}\int_{-r}^{0} B(\theta)K(\tau)\boldsymbol{x}(t+\tau+\theta)\mathrm{d}\tau\mathrm{d}\theta$ if (48) is adopted for a feedback. Handling such a scenario is far from straightforward and is currently deemed infeasible when addressing a dissipative synthesis problem.

## 6. Numerical Experiments and Simulations

We conducted numerical experiments with our proposed methodologies to demonstrate their effectiveness and the major advantages of utilizing the Carathéodory framework in modeling FDEs where some parameters are subject to discontinuities and glitches. These tests involved using all the components of the proposed method to tackle the DSFC problem for general LTDSs. All computations were performed using MATLAB© with Yalmip [221] as the optimization parser, and Mosek [167], SDPT3 [222] were employed as the numerical solvers for SDP problems.

### *6.1. DSFC Design of an LTDS with Multiple Delays*

Consider a system in the form of (18) with delays $r_1 = 1$, $r_2 = 1.7$ and state space matrices in (53) with $n = m = 2$, $p = q = 1$, where $\varsigma(t) = \mathbf{0}_2$ holds almost everywhere w.r.t the Lebesgue measure. Function $\varsigma(t)$ could represent glitches or other anomalies in the input gain matrix $\widetilde{B}_1(\tau)$, which could not be characterized if we were to use the traditional derivative for $\dot{\boldsymbol{x}}(t)$ in (18). As our methodology is formulated for FDEs in the extended sense [3, page 58], $\varsigma(t)$ is effectively treated as zeros. In fact, $\varsigma(t)$ can be seamlessly incorporated into any matrix in (53) and be regarded as $\mathbf{0}_2$ by our proposed framework. This serves as an excellent example to demonstrate the benefits of utilizing the Carathéodory framework in the modeling of TDSs. Employing the spectral method from [188], we find that the nominal system with $\boldsymbol{w}(t) \equiv \mathbf{0}_q$ is unstable. Moreover, we employ the $\mathscr{L}^2$ gain

$$\gamma > 0,\ J_1 = -\gamma I_2,\ \widetilde{J} = I_2,\ J_2 = \mathbf{0}_2,\ J_3 = \gamma \tag{54}$$

as the performance objective for the SRF in (37) with $\gamma$ to be minimized.

The parameters in (53) are selected with sufficient levels of complexity to illustrate the strength of the proposed method. Our approach can handle practical examples such as those noted in Remark 7, whose DD-matrix kernels have simpler structures than those in (53). To the best of our knowledge, no existing methods can effectively solve the DSFC problem for LTDS (18) with the parameters in (53), due to the complexity of the DDs and the presence of multiple non-commensurate delays with non-Hurwitz $A_0$.

Assuming that all the system's states are measurable, our goal is to determine the controller gain of $\boldsymbol{u}(t) = K\boldsymbol{x}(t)$ to stabilize open-loop system (18), while minimizing the $\mathscr{L}^2$ gain. By examining the DD kernels in (53), let

$$\boldsymbol{\phi}_1(\tau) = \begin{bmatrix}\exp(\sin 20\tau)\\ \exp(\cos 20\tau)\end{bmatrix},\quad \boldsymbol{\phi}_2(\tau) = \begin{bmatrix}\exp(\sin 18\tau)\\ \exp(\cos 18\tau)\end{bmatrix} \tag{55a}$$

$$\varphi_1(\tau) = \frac{1}{\sin^2 1.2\tau + 1},\qquad \varphi_2(\tau) = \frac{1}{\cos^2 0.7\tau + 1} \tag{55b}$$

$$\boldsymbol{f}_1(\tau) = \begin{bmatrix}[\tau^i]_{i=0}^{\sigma_1}\\ [\sin 20i\tau]_{i=1}^{\lambda_1}\\ [\cos 20i\tau]_{i=1}^{\lambda_1}\end{bmatrix},\quad \boldsymbol{f}_2(\tau) = \begin{bmatrix}[\tau^i]_{i=0}^{\sigma_2}\\ [\sin 18i\tau]_{i=1}^{\lambda_2}\\ [\cos 18i\tau]_{i=1}^{\lambda_2}\end{bmatrix} \tag{55c}$$

for (21)–(6c) with $d_i = 2\lambda_i + \sigma_i + 1$, $\mu_i = 2$, $\delta_i = 1$ and

$$\begin{aligned} M_1 &= \begin{bmatrix}\mathbf{0}_{d_1} & \begin{bmatrix}\mathbf{0}_{\sigma_1}^\top & 0\\ \operatorname*{diag}_{i=1}^{\sigma_1} i & \mathbf{0}_{\sigma_1}\end{bmatrix} \oplus \begin{bmatrix}\mathsf{O}_{\lambda_1} & \operatorname*{diag}_{i=1}^{\lambda_1} 20i\\ -\operatorname*{diag}_{i=1}^{\lambda_1} 20i & \mathsf{O}_{\lambda_1}\end{bmatrix}\end{bmatrix}\\ M_2 &= \begin{bmatrix}\mathbf{0}_{d_2} & \begin{bmatrix}\mathbf{0}_{\sigma_2}^\top & 0\\ \operatorname*{diag}_{i=1}^{\sigma_2} i & \mathbf{0}_{\sigma_2}\end{bmatrix} \oplus \begin{bmatrix}\mathsf{O}_{\lambda_2} & \operatorname*{diag}_{i=1}^{\lambda_2} 18i\\ -\operatorname*{diag}_{i=1}^{\lambda_2} 18i & \mathsf{O}_{\lambda_2}\end{bmatrix}\end{bmatrix}\end{aligned} \tag{55d}$$

for the relations in (6b). As a result, we can construct

$$\begin{aligned} \widehat{A}_1 &= \begin{bmatrix} 0 & 0.8 & 0 & -0.3 & 0 & 0 & 0.1 & 0 & \mathbf{0}_{2\sigma_1}^\top & 3 & 0 & \mathbf{0}_{4\lambda_1-2}^\top\\ 0 & 0 & 0 & 0 & 1 & 0 & 0.3 & 0 & \mathbf{0}_{2\sigma_1}^\top & 0 & 3 & \mathbf{0}_{4\lambda_1-2}^\top \end{bmatrix}\\ \widehat{A}_2 &= \begin{bmatrix} 0 & 0 & 0 & 0.3 & 0 & -1 & 0 & 0 & \mathbf{0}_{2\sigma_2+2\lambda_2}^\top & -10 & 0 & \mathbf{0}_{2\lambda_2-2}^\top\\ 0.1 & 0 & 0 & 0 & 0 & 0 & 0 & 0.2 & \mathbf{0}_{2\sigma_2+2\lambda_2}^\top & 0 & -10 & \mathbf{0}_{2\lambda_2-2}^\top \end{bmatrix}\\ \widehat{B}_1 &= \begin{bmatrix} 0 & 0 & -0.01 & 0.1 & 0.01 & \mathbf{0}_{\sigma_1-1+2\lambda_1}^\top\\ 0 & 0 & 0.02 & 0 & 0.1 & \mathbf{0}_{\sigma_1-1+2\lambda_1}^\top \end{bmatrix}\\ \widehat{B}_2 &= \begin{bmatrix} 0.01 & 0.2 & 0.01 & \mathbf{0}_{\sigma_2+1+2\lambda_2}^\top\\ 0.02 & 0.1 & 0 & \mathbf{0}_{\sigma_2+1+2\lambda_2}^\top \end{bmatrix}\\ \widehat{C}_1 &= \begin{bmatrix} 0 & 0 & 0 & 0 & 0 & 1 & 0.7 & -0.2 & \mathbf{0}_{2\sigma_1}^\top & 0 & 0 & \mathbf{0}_{2\lambda_1-2}^\top & 1 & 0 & \mathbf{0}_{2\lambda_1-2}^\top\\ -0.5 & 0 & 0 & 0 & 0 & 0 & 0.4 & 0.8 & \mathbf{0}_{2\sigma_1}^\top & 0 & -1 & \mathbf{0}_{2\lambda_1-2}^\top & 0 & 0 & \mathbf{0}_{2\lambda_1-2}^\top \end{bmatrix}\\ \widehat{C}_2 &= \begin{bmatrix} 0 & 0 & 0 & 1 & 0 & 0 & 0.2 & 0.3 & \mathbf{0}_{2\sigma_2}^\top & 1 & 0 & \mathbf{0}_{4\lambda_2-2}^\top\\ 0 & 0 & 0 & 0 & 0 & -1 & 0 & 0.1 & \mathbf{0}_{2\sigma_2}^\top & 0 & 0 & \mathbf{0}_{4\lambda_2-2}^\top \end{bmatrix}\\ \widehat{\mathfrak{B}}_1 &= \begin{bmatrix} 0.1 & 0 & -0.1 & 0 & 0.01 & \mathbf{0}_{\sigma_1+2\lambda_1-1}^\top\\ 0.2 & 0 & 0 & 0 & 0 & \mathbf{0}_{\sigma_1+2\lambda_1-1}^\top \end{bmatrix}\\ \widehat{\mathfrak{B}}_2 &= \begin{bmatrix} 0.01 & 0.2 & 0.1 & 0 & 0 & \mathbf{0}_{\sigma_2+2\lambda_2-1}^\top\\ 0.02 & 0 & 0.2 & 0 & 0 & \mathbf{0}_{\sigma_2+2\lambda_2-1}^\top \end{bmatrix} \end{aligned} \tag{56}$$

to satisfy the conditions in (21)–(6c).

**Remark 13.** The functions in (55a)–(55c) were selected for several reasons. First, $\boldsymbol{\phi}_i(\cdot)$ in (55a) can be easily approximated by $\boldsymbol{f}_i(\cdot)$ in (55c), and some functions in $\boldsymbol{f}_i(\cdot)$ are included by the DD integral matrices specified in (53). On the other hand, the functions in $\varphi_i(\cdot)$ in (55b) are directly factorized, as they are hard to approximate by $\boldsymbol{f}_i(\cdot)$ in (55c) with $\sigma_i, \lambda_i$ of reasonable values due to their irregular oscillation behaviors. Finally, the combination of trigonometric functions and polynomials for $\boldsymbol{f}_i(\cdot)$ in (55c) would ensure the feasibility of the synthesis condition, as the structure of $\boldsymbol{f}_i(\cdot)$ directly determines the generality of KF (C.1a). Hence, the choice in (55a)–(55c) demonstrates all the components of the proposed KSD concept in conjunction with the SDP framework, which balances feasibility and the implied computational complexity $\mathcal{O}(d^2n^2)$ determined by

$$
A_0 = \begin{bmatrix} -2 & 0 \\ 2 & 0.01 \end{bmatrix}, \; A_1 = \begin{bmatrix} -1 & 0.1 \\ 0.2 & 0 \end{bmatrix}, \; A_2 = \begin{bmatrix} -0.1 & 0 \\ 0 & -0.2 \end{bmatrix}, \quad B_0 = \begin{bmatrix} 0 \\ 1 \end{bmatrix}, \quad B_1 = \begin{bmatrix} 0.01 \\ 0.1 \end{bmatrix}, \quad B_2 = -\begin{bmatrix} 0.1 \\ 0.1 \end{bmatrix}
$$

$$
\widetilde{A}_1(\tau) = \begin{bmatrix} 0.1 + 3\sin(20\tau) & 0.8\exp(\sin 20\tau) - 0.3\exp(\cos 20\tau) \\ 0.3 + \dfrac{1}{\sin^2(1.2\tau) + 1.0} & 3\sin(20\tau) \end{bmatrix}, \widetilde{A}_2(\tau) = \begin{bmatrix} -10\cos(18\tau) & 0.3\exp(\cos 18\tau) - \dfrac{1}{\cos^2 0.7\tau + 1} \\ 0.1\exp(\sin 18\tau) & 0.2 - 10\cos(18\tau) \end{bmatrix},
$$

$$
\widetilde{B}_1(\tau) = \begin{bmatrix} 0.01\tau - \dfrac{0.01}{\sin^2(1.2\tau) + 1} + 0.1 \\ 0.1\tau + \dfrac{0.02}{\sin^2(1.2\tau) + 1} \end{bmatrix} + \varsigma(t), \quad \widetilde{B}_2(\tau) = \begin{bmatrix} 0.2\exp(\cos 18\tau) + 0.01\exp(\sin 18\tau) + \dfrac{0.01}{\cos^2(0.7\tau) + 1} \\ 0.1\exp(\cos 18\tau) + 0.02\exp(\sin(18\tau)) \end{bmatrix} \tag{53}
$$

$$
C_0 = \begin{bmatrix} -0.1 & 0.2 \\ 0 & 0.1 \end{bmatrix}, \; C_1 = \begin{bmatrix} -0.1 & 0 \\ 0 & 0.2 \end{bmatrix}, \; C_2 = \begin{bmatrix} 0 & 0.1 \\ -0.1 & 0 \end{bmatrix}, \; \mathfrak{B}_1 = \begin{bmatrix} 0.01 \\ 0.01 \end{bmatrix}, \; \mathfrak{B}_2 = -\begin{bmatrix} 0.01 \\ 0.1 \end{bmatrix}, \; D_1 = \begin{bmatrix} 0.2 \\ 0.3 \end{bmatrix},
$$

$$
\widetilde{C}_1(\tau) = \begin{bmatrix} 0.7 + \cos(20\tau) & \frac{1}{\sin^2 1.2\tau + 1} - 0.2 \\ 0.4 - 0.5\exp(\sin 20\tau) & 0.8 - \sin(20\tau) \end{bmatrix}, \quad \widetilde{C}_2(\tau) = \begin{bmatrix} 0.2 + \sin(18\tau) & 0.3 + \exp(\cos 18\tau) \\ 0 & 0.1 - \frac{1}{\cos^2 0.7\tau + 1} \end{bmatrix}, \quad D_2 = \begin{bmatrix} 0.12 \\ 0.1 \end{bmatrix}
$$

$$
\widetilde{\mathfrak{B}}_1(\tau) = \begin{bmatrix} 0.01\tau + 0.1\exp(\sin 20\tau) - \dfrac{0.1}{\sin^2(1.2\tau) + 1} \\ 0.2\exp(\sin 20\tau) \end{bmatrix}, \quad \widetilde{\mathfrak{B}}_2(\tau) = \begin{bmatrix} 0.2\exp(\cos 18\tau) + 0.01\exp(\sin 18\tau) + \dfrac{0.1}{\cos^2(0.7\tau) + 1} \\ 0.02\exp(\sin 18\tau) + \dfrac{0.2}{\cos^2(0.7\tau) + 1} \end{bmatrix}
$$

$\mathsf{dim}(\boldsymbol{f}_i(\tau)) = d_i$.

To compute $K$, apply Theorem 2 to (32) with $\sigma_1 = \sigma_2 = \lambda_1 = \lambda_2 = 1$ and $\alpha_i = 0$, $i = 2, \ldots, \beta$, $\alpha_1 = 5$ and the parameters in (53)–(56), where $\Gamma_i$, $\mathfrak{F}_i$, $\mathfrak{H}_i$, $\mathfrak{E}_i$ in (8)–(10) are computed via the MATLAB© function `vpaintegral` that performs numerical integrations with high-level variable precision. Our SDP program yields numerical results $K = -\begin{bmatrix} 1.3794 & 1.8668 \end{bmatrix}$ with $\min\gamma = 0.8986$. This $K$ is then used for initializing Algorithm 1. After running Algorithm 1 for the same system with different numbers of iterations (NoIs) for the same system with different sets of $\sigma_i, \lambda_i$, the numerical results are listed in Tables 1 and 2, where the spectral abscissae of the resulting CLSs with $\boldsymbol{w}(t) \equiv \mathbf{0}_q$ were calculated by the method in [188]. Our results bring out the fact that adding more WDLIFs (larger $\lambda_1, \lambda_2$) to $\boldsymbol{f}_i(\cdot)$ may increase the feasibility of the synthesis conditions, leading to smaller $\min\gamma$. Moreover, it confirms that using Algorithm 1 can produce controller gains of significantly better performance than those produced by employing Theorem 2 alone ($\min\gamma = 0.8986$). Hence, this illustrates the contribution of Algorithm 1.

**Table 1:** Controller gains $K$ with $\min\gamma$ produced with $\sigma_1 = \sigma_2 = \lambda_1 = \lambda_2 = 1$, $d_1 = d_2 = 4$

| $K$ | $\begin{bmatrix} -1.5456 \\ -1.9359 \end{bmatrix}^\top$ | $\begin{bmatrix} -1.5365 \\ -1.9539 \end{bmatrix}^\top$ | $\begin{bmatrix} -1.5180 \\ -1.9696 \end{bmatrix}^\top$ | $\begin{bmatrix} -1.5033 \\ -1.9815 \end{bmatrix}^\top$ |
|---|---|---|---|---|
| $\min\gamma$ | 0.6573 | 0.6542 | 0.6523 | 0.6509 |
| SA | $-0.7223$ | $-0.7214$ | $-0.7224$ | $-0.7233$ |
| NoIs | 5 | 10 | 15 | 20 |

**Table 2:** Controller gains $K$ with $\min\gamma$ produced with $\sigma_1 = \sigma_2 = 1$, $\lambda_1 = \lambda_2 = 2$, $d_1 = d_2 = 6$.

| $K$ | $\begin{bmatrix} -1.5538 \\ -1.9566 \end{bmatrix}^\top$ | $\begin{bmatrix} -1.5848 \\ -1.9638 \end{bmatrix}^\top$ | $\begin{bmatrix} -1.5870 \\ -1.9721 \end{bmatrix}^\top$ | $\begin{bmatrix} -1.5810 \\ -1.9805 \end{bmatrix}^\top$ |
|---|---|---|---|---|
| $\min\gamma$ | 0.6443 | 0.6398 | 0.6376 | 0.6361 |
| SA | $-0.7179$ | $-0.7113$ | $-0.7099$ | $-0.7099$ |
| NoIs | 5 | 10 | 15 | 20 |

For numerical simulations, we consider the CLS in (32) with the parameters in (53) and controller gain $K = \begin{bmatrix} -1.5810 & -1.9805 \end{bmatrix}$ in Table 2 that guarantees $\min\gamma = 0.6361$. Let $t_0 = 0$ and $\boldsymbol{\psi}(\tau) = \begin{bmatrix} 5 & 3 \end{bmatrix}^\top$, $\tau \in [-r_2, 0]$ as the initial condition, and $w(t) = 5\sin 3\pi t(\mathit{1}(t) - \mathit{1}(t-10))$ as the disturbance where $\mathit{1}(t)$ is the Heaviside step function. For noise and glitch signal $\varsigma(t)$, we employ the Band-Limited White Noise block in Simulink to generate a white noise signal $n(t)$ for $\varsigma(t) = n(t)\begin{bmatrix} 1 & 1 \end{bmatrix}^\top$ in (53), where `Sample time` $=$ $0.002s$ and the default values for `Seed` and `Noise power` were adopted. Since $n(t)$ can only be realized as a discrete sequence $n(t) = n(kT)$ within a numerical simulation environment, function $\varsigma(t) = \varsigma(kT)$ has only a **finite** number of nonzero values, which satisfies $\widetilde{\forall} t \geq 0$, $\varsigma(t) = 0$ and is in line with the definitions in (53). Finally, the computations were performed via the ODE solver `ode8` with $0.002s$ as the sampling time.

The results in Figures 1–3 clearly demonstrate that the goal of dissipative stabilization for the CLS is achieved by the resulting controller, as we can observe the gradual suppression of noise effects and their final extinction from the plots of $\boldsymbol{x}(t), u(t), \boldsymbol{z}(t)$, albeit noise $n(t)$ is always nonzero. We can also see a substantial attenuation of the effect of disturbance $w(t) = 5 \sin 3\pi t(\mathit{1}(t) - \mathit{1}(t-10))$ in terms of the oscillation amplitudes over $t \in [5, 10]$. This underscores both the advantages of the Carathéodory framework in addressing discontinuities and glitches in the system parameters, and the guaranteed $\mathcal{L}^2$ performance design of our robust controller.

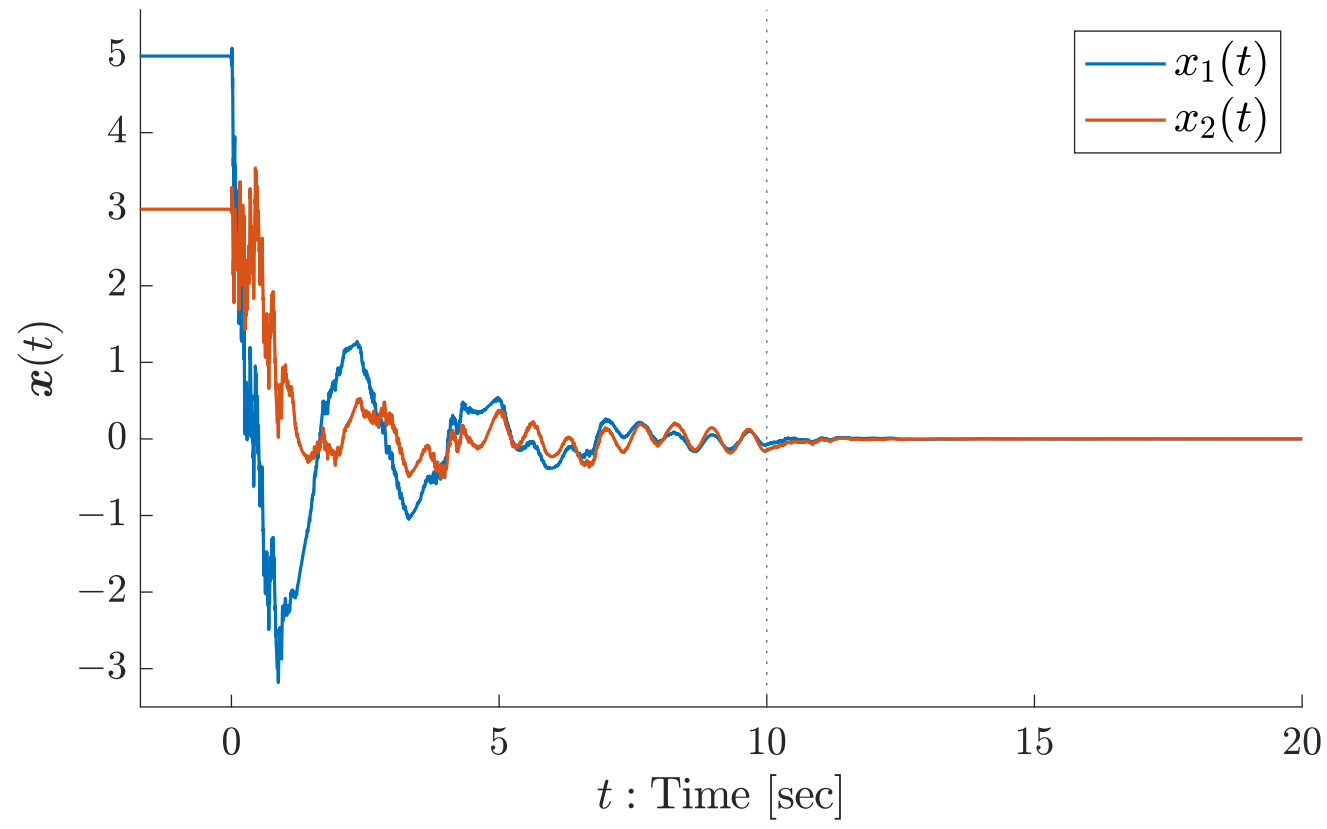

**Figure 1:** Plot of $\boldsymbol{x}(t)$ with $K$ in Table 2 ensuring $\min \gamma = 0.6361$

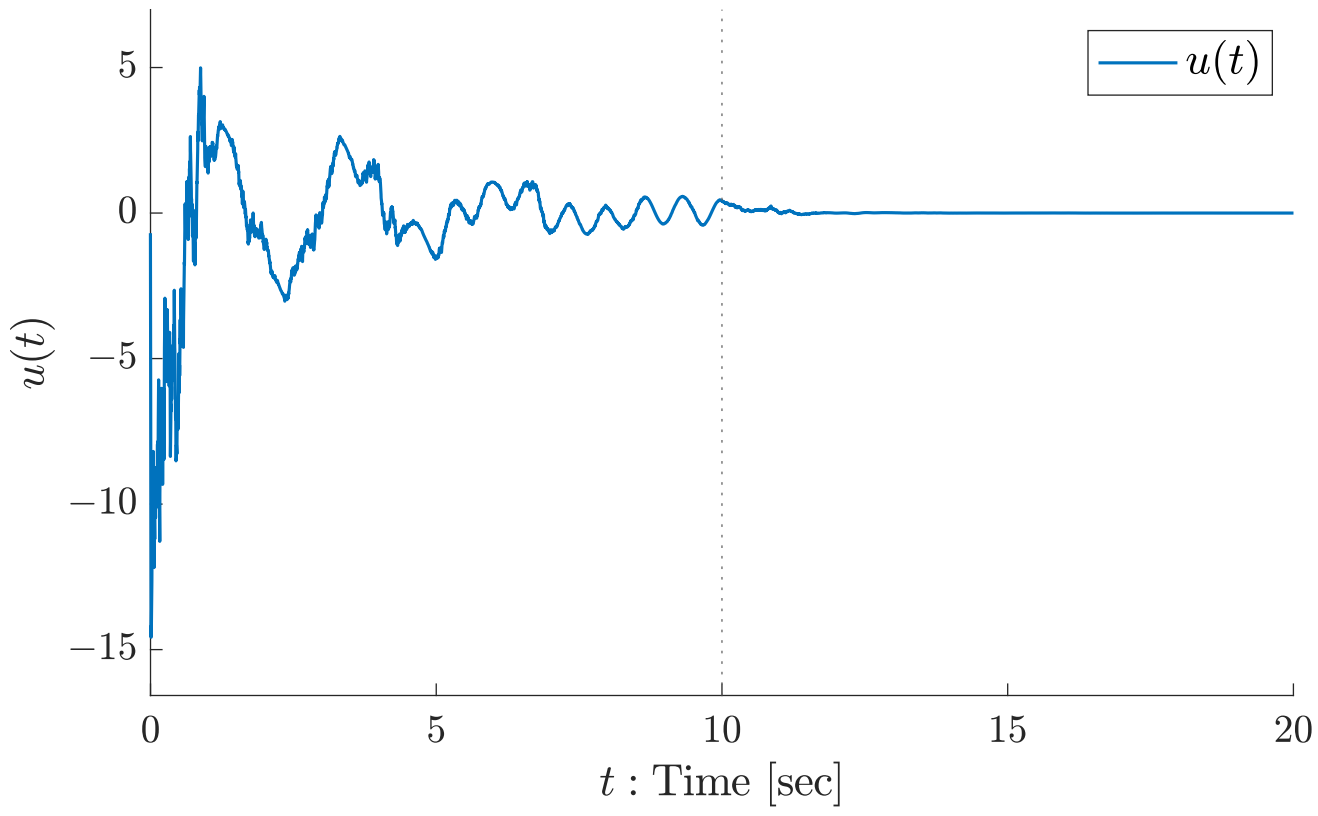

**Figure 2:** Plot of $u(t) = K\boldsymbol{x}(t)$ ensuring $\min \gamma = 0.6361$

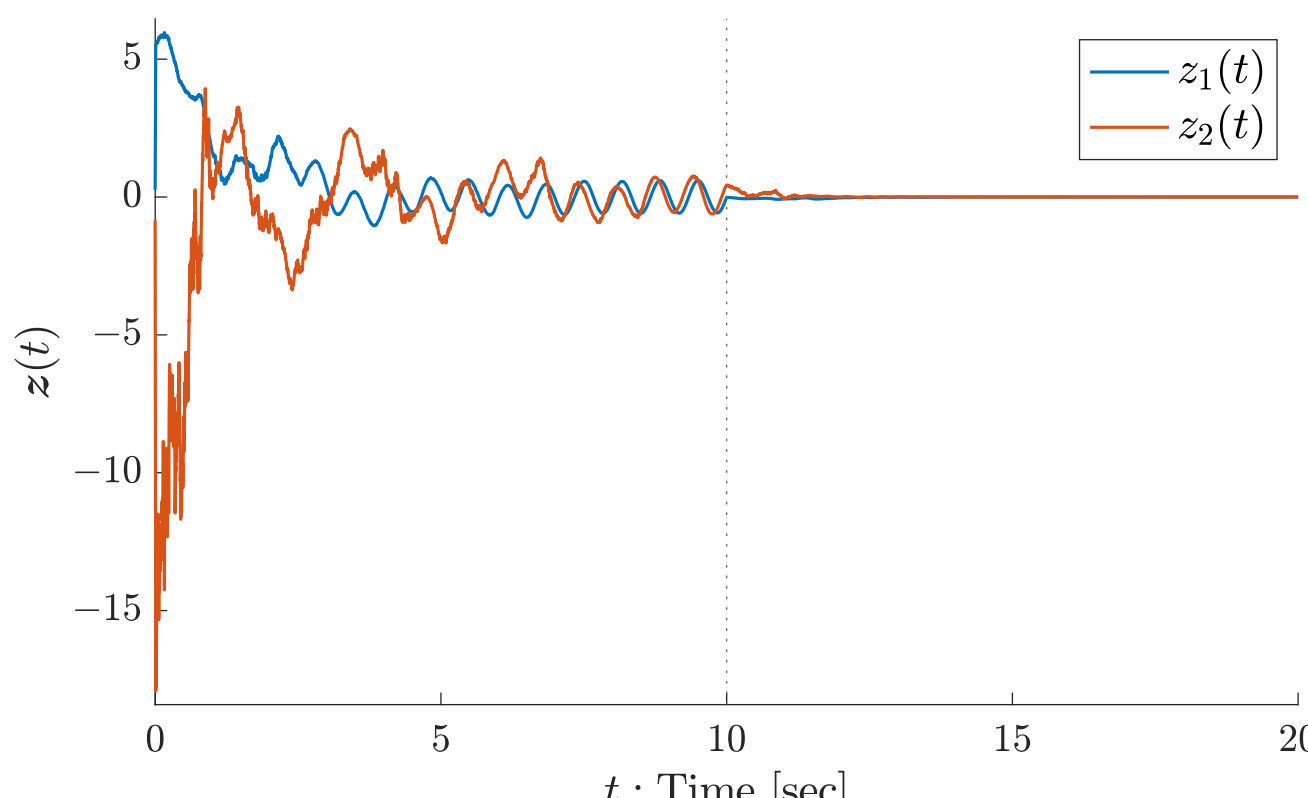

**Figure 3:** Plot of $\boldsymbol{z}(t)$ of (32) ensuring $\min \gamma = 0.6361$

### 6.2. DSFC Design of an LTDS with Controller Delays

In this subsection, we present the testing results of a numerical example that can demonstrate the advantages of employing controllers with general delay structures in the absence of input delays in (18).

Consider the open-loop system in (18) with the same parameters as in Section 6.1 except for the input delay matrices $B_0 = \begin{bmatrix} 0 & 1 \end{bmatrix}^\top + \varsigma(t)$, $B_i = \widetilde{B}_i(\tau) = \mathsf{O}_{n,p}$ and $\mathfrak{B}_i = \widetilde{\mathfrak{B}}_i(\tau) = \mathsf{O}_{m,p}$, $\forall i \in \mathbb{N}_\nu$, where $\varsigma(t)$ follows the same definition in Section 6.1. As such, we can employ the controller in (48)–(49) to achieve dissipative stabilization while minimizing $\gamma$ in (54).

The steps of computing controller gains here are entirely identical to those in the previous subsection, except that we applied Algorithm 2 to the CLS in (50) in conjunction with Theorems 3 and 4. Specifically, we applied the KSD concept with (49) to the DD kernel matrices in (53) and (48) using the same functions in (55a)–(55c) for $\boldsymbol{g}_i(\cdot)$ with $M_i$ in (55d). This results in the same parameters $\widehat{A}_i, \widehat{C}_i$ as in (56), while $K_i, \mathcal{K}_i$ in (49) are the controller gains to be computed.

Assume $\alpha_i = 0, i = 2, \cdots, \beta$, and $\alpha_1 = 5$, respectively, for the utilization of Theorem 4. We apply the entire package consisting of Theorems 3 and 4 and Algorithm 2 to the CLS in (50) to compute the controller parameters in (49). This produced the numerical results listed in Tables 3 and 4. To further demonstrate the advantages of the general delay structures in (49) in the absence of input delays, we also computed $\min \gamma$ for the case of $\boldsymbol{u}(t) = K_0 \boldsymbol{x}(t)$ with $K_i = \widetilde{K}_i(\tau) = \mathsf{O}_{p,n}, \forall i \in \mathbb{N}_\nu$, which corresponds to a static state controller without memories. This yielded the numerical results in Tables 5 and 6, where the values of $\min \gamma$ can be compared to those in Tables 3 and 4 under the same $\sigma_1, \sigma_2, \lambda_1, \lambda_2$.

It is evident that for both cases in Tables 3–4 and Tables 5–6, adding more WDLIFs to $\boldsymbol{f}_i(\cdot)$ leads to smaller $\min \gamma$ similar to the results in Tables 1 and 2. However, comparing $\min \gamma$ in Table 3 with Table 5, and Table 4 with Table 6 shows that using controllers with general delay structures can significantly improve the $\mathcal{L}^2$ gain performance $\min \gamma$ over static controller $\boldsymbol{u}(t) = K_0 \boldsymbol{x}(t)$. This justifies the use of the delay structures in (48) in the absence of input delays, despite the additional resources required for the implementation of delay structures in the controller.

Since the CLS in (50) is a retarded FDE, its nominal stability is guaranteed by [223, Th. 2] when the DDs in (49) are implemented numerically via standard quadratures [191], as long as the numerical accuracy reaches a certain degree. This property ensures that the DDs in (49) can be implemented numerically for real-world applications.

**Table 3:** $\min \gamma$ computed by Algorithm 2 with $\sigma_1 = \sigma_2 = \lambda_1 = \lambda_2 = 1, d_1 = d_2 = 4$

| $\min \gamma$ | 0.5242 | 0.524 | 0.5238 | 0.5237 |
|---|---|---|---|---|
| Spectral Abscissa | −0.6983 | −0.6979 | −0.6976 | −0.6989 |
| Number of Iterations | 5 | 10 | 15 | 20 |

**Table 4:** $\min \gamma$ computed by Algorithm 2 with $\sigma_1 = \sigma_2 = 1$, $\lambda_1 = \lambda_2 = 2$, $d_1 = d_2 = 6$

| $\min \gamma$ | 0.5236 | 0.5234 | 0.5232 | 0.523 |
|---|---|---|---|---|
| Spectral Abscissa | −0.7039 | −0.6915 | −0.6899 | −0.6907 |
| Number of Iterations | 5 | 10 | 15 | 20 |

**Table 5:** $\min \gamma$ computed by Algorithm 2 with $\sigma_1 = \sigma_2 = \lambda_1 = \lambda_2 = 1$, $d_1 = d_2 = 4$, $K_i = \widetilde{K}_i(\tau) = \mathsf{O}_{p,n}$

| $\min \gamma$ | 0.5785 | 0.5760 | 0.5736 | 0.5714 |
|---|---|---|---|---|
| Spectral Abscissa | −0.7259 | −0.7319 | −0.7358 | −0.7385 |
| Number of Iterations | 5 | 10 | 15 | 20 |

For numerical simulation, we use the parameters

$$K_0 = \begin{bmatrix} -9.1247 & -22.9729 \end{bmatrix}, \quad K_1 = \begin{bmatrix} 0.2429 & -0.117 \end{bmatrix}$$
$$K_2 = \begin{bmatrix} 0.0972 & 0.1237 \end{bmatrix}$$
$$\mathcal{K}_1 = [\, -0.2358 \ -0.3997 \ -0.3585 \ 0.1859 \ -0.4752$$
$$\cdots -0.8404 \ 4.0466 \ -2.6485 \ -0.9318 \ -3.3031 \ -1.2266$$
$$\cdots -2.4909 \ 0.1010 \ -0.5549 \ 0.8075 \ 1.4130 \ 0.0620 \ 0.148 \,]$$
$$\mathcal{K}_2 = [\, 0.0276 \ -0.0464 \ 0.2009 \ -0.3101 \ -0.8961 \ -0.9072$$
$$\cdots 6.3285 \ -1.1248 \ 5.0129 \ -1.1893 \ 1.1018 \ 1.5506$$
$$\cdots 0.3628 \ -0.0725 \ 4.3602 \ 6.5996 \ 0.0282 \ 0.6417 \,]$$

corresponding to $\min \gamma = 0.523$ in Table 4 for the controller in (49). The rest of the simulation setting is identical to that of Section 6.1 including the configuration of noise signal $\varsigma(t)$. The numerical results of simulations are presented in Figures 4–6.

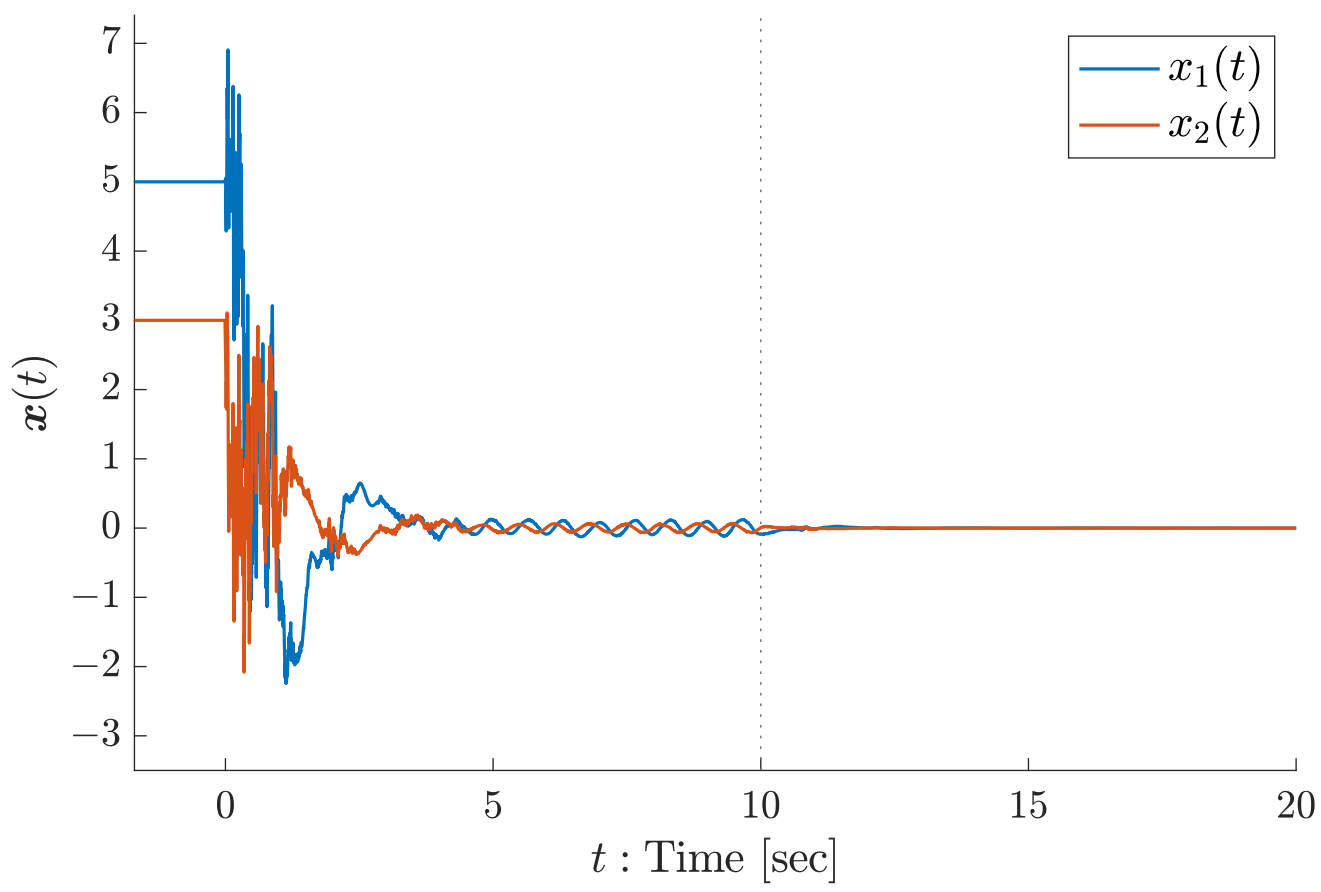

**Figure 4:** Plot of $\boldsymbol{x}(t)$ of (50) with parameters (57)

## 7. Conclusion and Remarks

We have proposed an SDP framework based on the KSD concept developed by [164] that can effectively address two different DSFC problems of the LTDS in (18). While these two problems differ due to the presence of general delays in the input and controller, the structures of the CLSs in (32) and (50) indicate that both cases can be similarly addressed by the proposed SDP methodology. The generality of the model in (18) is guaranteed as it closely mirrors the expressions of general LTDSs (14) expressed via the Lebesgue-Stieltjes integrals (16), where the FDEs are understood in the extended sense. It has been demonstrated that the KSD concept coordinates effectively with the construction of our general KF in (C.1a), overcoming the challenges posed by the infinite dimension of DD-matrix kernels while ensuring that the synthesis conditions in Theorems 1 and 2 are represented by matrix inequalities with finite dimensions. The generality and effectiveness of our framework are rooted in three key ingredients: first, the general structure of functions $\boldsymbol{f}_i(\cdot)$ in both Proposition 1 and the functional in (C.1a), in conjunction with the use of the KSD concept that introduces no conservatism in terms of kernel lifting; second, the extensive use of the general integral inequality (B.2) derived from the least-squares principle; and third, a robust package comprising two theorems and an iterative algorithm aiming to minimize the conservatism arising from the non-convexity in (38c). Numerical experiments have demonstrated that our framework can effectively compute dissipative controllers for systems with intricate delay structures, even when the kernel functions exhibit vastly different characteristics.

**Table 6:** $\min \gamma$ computed by Algorithm 2 with $\sigma_1 = \sigma_2 = 1$, $\lambda_1 = \lambda_2 = 2$, $d_1 = d_2 = 6$, $K_i = \widetilde{K}_i(\tau) = \mathsf{O}_{p,n}$

| $\min \gamma$ | 0.5723 | 0.5669 | 0.5626 | 0.559 |
|---|---|---|---|---|
| Spectral Abscissa | −0.7179 | −0.7113 | −0.7099 | −0.7099 |
| Number of Iterations | 5 | 10 | 15 | 20 |

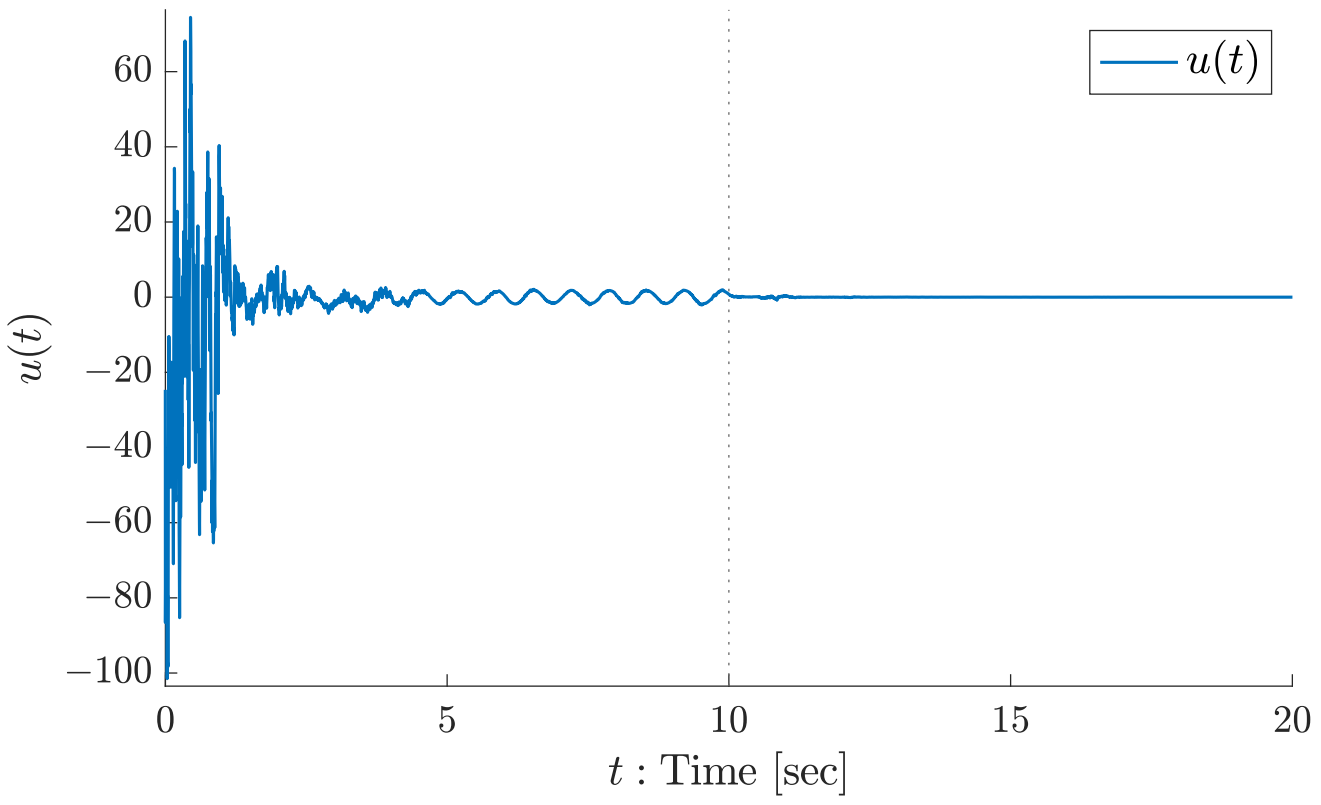

**Figure 5:** The trajectory of control action $u(t)$ of (48) with (57)

The methodologies proposed in this paper offer a promising foundation for developing new solutions for delay-related practical systems, such as neural networks [33], event-triggered systems [36] and other relevant applications [224]. The above claim is evidenced by the numerous publications [33, 36, 225–228] on engineering systems that utilized the decomposition approach developed by [152], which is a special instance of the more complex KSD concept proposed here. Moreover, the pro-posed SDP framework can be utilized as a blueprint to solve other semi-open problems in the field of LTDSs in several directions, such as dissipative dynamical output feedback control and filtering problems. In particular, the proposed KSD can effectively parameterize the matrix-valued integral kernels of the complex DDs in the predictor controllers for LTDSs with both state and input pointwise delays [104, 105, 229], as we have demonstrated in (49) its ability to represent controllers with complex delay structures. This would be of significant research interest since designing dissipative predictor controllers for LTDSs with both state and input delays remains one of the most important semi-open problems in the field of time-delay systems.

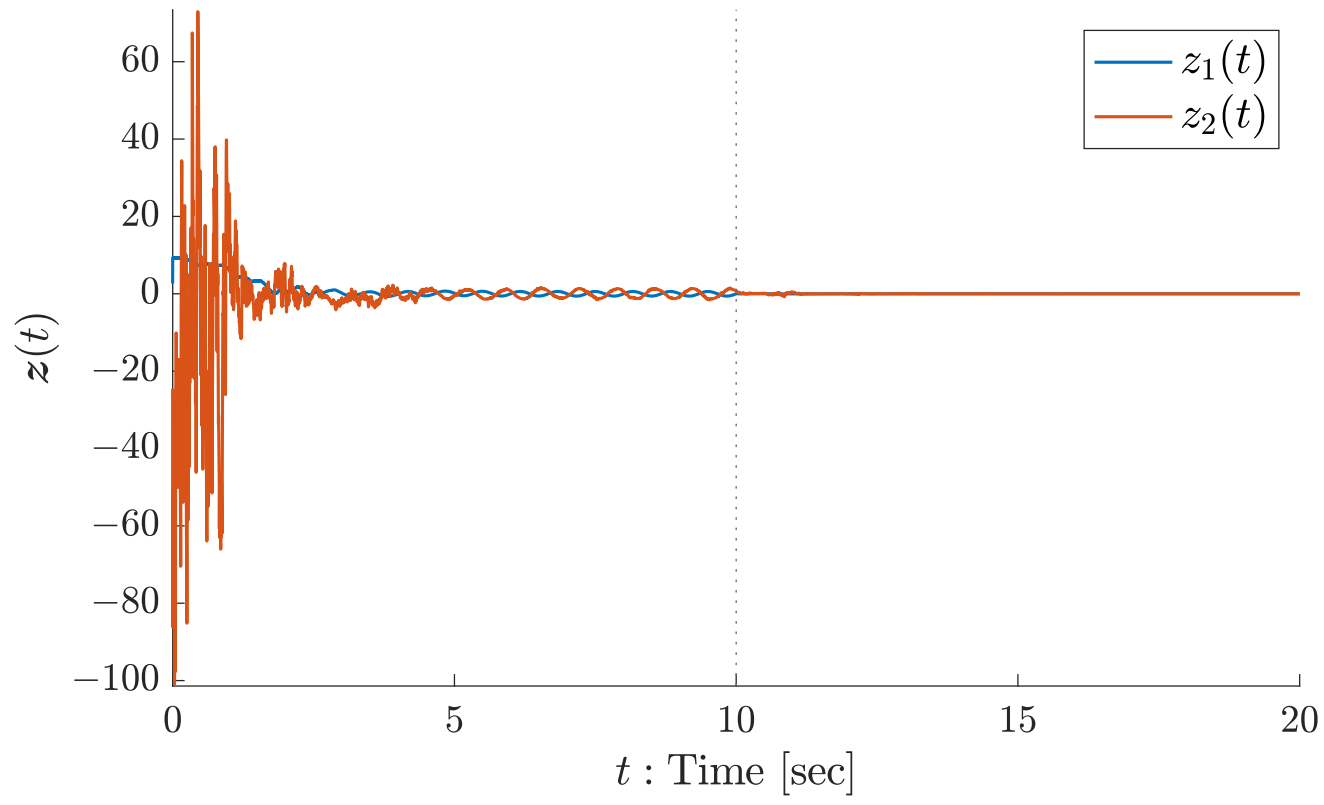

**Figure 6:** The trajectories of $\boldsymbol{z}(t)$ of (50) with parameters (57)

## Appendix A. Matrix Identities

**Lemma 3.** $\forall X \in \mathbb{R}^{n\times m},\ \forall Y \in \mathbb{R}^{m\times p},\ \forall Z \in \mathbb{R}^{q\times r}$,

$$
\begin{aligned}
&(X\otimes I_q)(Y\otimes Z) = XY\otimes Z \\
&\quad = XY\otimes ZI_r = (X\otimes Z)(Y\otimes I_r) \\
&\quad = I_n XY\otimes (ZI_r) = (I_n\otimes Z)(XY\otimes I_r).
\end{aligned} \tag{A.1}
$$

*Moreover,* $\forall X \in \mathbb{R}^{n\times m}$, *we have*

$$
\begin{bmatrix} A & B \\ C & D \end{bmatrix} \otimes X = \begin{bmatrix} A\otimes X & B\otimes X \\ C\otimes X & D\otimes X \end{bmatrix},\quad I_\nu \otimes X = \operatorname*{diag}_{i=1}^{\nu} X \tag{A.2}
$$

*for any* $A, B, C, D$ *with appropriate dimensions.*

*Finally, for any matrix* $X_i, Y_i, Z_i$ *and vector* $\mathbf{v}_i$ *with matching dimensions,*

$$
\begin{aligned}
\sum_{i=1}^{n} (X_i + Y_i Z_i)\mathbf{v}_i
&= \llbracket X_i + Y_i Z_i \rrbracket_{i=1}^{n} \left[\mathbf{v}_i\right]_{i=1}^{n} \qquad \text{(A.3a)} \\
&= \left( \llbracket X_i \rrbracket_{i=1}^{n} + \llbracket Y_i Z_i \rrbracket_{i=1}^{n} \right) \left[\mathbf{v}_i\right]_{i=1}^{n} \\
&= \llbracket X_i \rrbracket_{i=1}^{n} \left[\mathbf{v}_i\right]_{i=1}^{n} + \llbracket Y_i Z_i \rrbracket_{i=1}^{n} \left[\mathbf{v}_i\right]_{i=1}^{n} \\
&= \llbracket X_i \rrbracket_{i=1}^{n} \left[\mathbf{v}_i\right]_{i=1}^{n} + \llbracket Y_i \rrbracket_{i=1}^{n} \left[\operatorname*{diag}_{i=1}^{n} Z_i\right] \left[\mathbf{v}_i\right]_{i=1}^{n}, \qquad \text{(A.3b)}
\end{aligned}
$$

$$
\begin{aligned}
\sum_{i=1}^{n} Y_i Z_i \mathbf{v}_i &= \llbracket Y_i \rrbracket_{i=1}^{n} \left[\operatorname*{diag}_{i=1}^{n} Z_i\right] \left[\mathbf{v}_i\right]_{i=1}^{n} \\
&= \begin{bmatrix} \llbracket Y_i \rrbracket_{i=1}^{n} & \mathsf{O} \end{bmatrix} \left[ \left[\operatorname*{diag}_{i=1}^{n} Z_i\right] \oplus \mathsf{O} \right] \left[\mathbf{v}_i\right]_{i=1}^{n+1},
\end{aligned} \tag{A.3c}
$$

$$
\left[Y_i Z_i\right]_{i=1}^{n} = \left(\operatorname*{diag}_{i=1}^{n} Y_i\right) \left[Z_i\right]_{i=1}^{n}, \tag{A.3d}
$$

$$
\begin{aligned}
\operatorname*{diag}_{i=1}^{n} Y_i Z_i \otimes I_\nu &= \left( \left[\operatorname*{diag}_{i=1}^{n} Y_i\right] \left[\operatorname*{diag}_{i=1}^{n} Z_i\right] \right) \otimes I_\nu \\
&= \left[\operatorname*{diag}_{i=1}^{n} Y_i \otimes I_\nu\right] \left[\operatorname*{diag}_{i=1}^{n} Z_i \otimes I_\nu\right] = \operatorname*{diag}_{i=1}^{n} \left(Y_i Z_i \otimes I_\nu\right).
\end{aligned} \tag{A.3e}
$$

*Proof.* The proof is accomplished by utilizing the properties of the Kronecker products [199] in conjunction with the definition of the $n$-ary operators $\mathsf{diag}_{i=1}^{n}$, $\llbracket\cdot\rrbracket_{i=1}^{n}$ and $[\cdot]_{i=1}^{n}$. ■

## Appendix B. A General Integral Inequality

We define a weighted Lebesgue function space as

$$
\mathscr{L}^2_{\varpi}\left(\mathcal{K};\mathbb{R}^d\right) := \left\{ \boldsymbol{h}(\cdot) \in \mathcal{M}\left(\mathcal{K};\mathbb{R}^d\right) : \left\|\boldsymbol{h}(\cdot)\right\|_{\varpi} < +\infty \right\} \tag{B.1}
$$

with $d \in \mathbb{N}$ and $\|\boldsymbol{h}(\cdot)\|^2_{\varpi} := \int_{\mathcal{K}} \varpi(\tau)\boldsymbol{h}^\top(\tau)\boldsymbol{h}(\tau)\mathsf{d}\tau$, where $\varpi(\cdot) \in \mathcal{M}(\mathcal{K};\mathbb{R}_{\geq 0})$, and weighted function $\varpi(\cdot)$ has at most countably many zeros. Furthermore, $\mathcal{K} \subseteq \mathbb{R}\cup\{\pm\infty\}$ and its Lebesgue measure [49] is nonzero.

**Lemma 4.** *Given a sequence of measurable sets* $\mathcal{K}_i \subseteq \mathcal{K}$ *and* $\varpi(\cdot)$ *in* (B.1) *and* $\mathbb{S}^n \ni X_i \succeq 0$ *with* $i \in \mathbb{N}_\nu$ *and* $\nu \in \mathbb{N}$. *Consider functions* $\mathfrak{h}_i(\cdot) \in \mathscr{L}^2_{\varpi}\left(\mathcal{K}_i;\mathbb{R}^{l_i}\right)$ *and* $\mathfrak{q}_i(\cdot) \in \mathscr{L}^2_{\varpi}\left(\mathcal{K}_i;\mathbb{R}^{\lambda_i}\right)$ *satisfying*

$$
\forall i \in \mathbb{N}_\nu,\ \int_{\mathcal{K}_i} \varpi(\tau) \begin{bmatrix} \mathfrak{q}_i(\tau) \\ \mathfrak{h}_i(\tau) \end{bmatrix} \begin{bmatrix} \mathfrak{q}_i^\top(\tau) & \mathfrak{h}_i^\top(\tau) \end{bmatrix} \mathsf{d}\tau \succ 0
$$

*with* $l_i \in \mathbb{N}$ *and* $\lambda_i \in \mathbb{N}\cup\{0\}$ *for all* $i \in \mathbb{N}_\nu$, *which indicates* $\mathfrak{H}_i = \int_{\mathcal{K}_i} \varpi(\tau)\mathfrak{h}_i(\tau)\mathfrak{h}_i^\top(\tau)\mathsf{d}\tau \succ 0$. *Then* $\forall \boldsymbol{x}(\cdot) \in \mathscr{L}^2_{\varpi}(\mathcal{K};\mathbb{R}^n)$

$$
\begin{aligned}
&\sum_{i=1}^{\nu} \int_{\mathcal{K}_i} \varpi(\tau)\boldsymbol{x}^\top(\tau)X_i\boldsymbol{x}(\tau)\mathsf{d}\tau \\
&\geq [*]\left[\operatorname*{diag}_{i=1}^{\nu} \mathfrak{H}_i^{-1}\otimes X_i\right] \left[\int_{\mathcal{K}_i} \varpi(\tau)\left[\mathfrak{h}_i(\tau)\otimes I_n\right]\boldsymbol{x}(\tau)\mathsf{d}\tau\right]_{i=1}^{\nu} \\
&\quad + [*]\left[\operatorname*{diag}_{i=1}^{\nu} \mathfrak{E}_i^{-1}\otimes X_i\right] \left[\int_{\mathcal{K}_i} \varpi(\tau)\left[\mathfrak{e}_i(\tau)\otimes I_n\right]\boldsymbol{x}(\tau)\mathsf{d}\tau\right]_{i=1}^{\nu} \\
&\geq [*]\left[\operatorname*{diag}_{i=1}^{\nu} \mathfrak{H}_i^{-1}\otimes X_i\right] \left[\int_{\mathcal{K}_i} \varpi(\tau)\left[\mathfrak{h}_i(\tau)\otimes I_n\right]\boldsymbol{x}(\tau)\mathsf{d}\tau\right]_{i=1}^{\nu}
\end{aligned} \tag{B.2}
$$

*where function* $\mathfrak{e}_i(\tau) = \mathfrak{q}_i(\tau) - \mathfrak{A}_i\mathfrak{h}_i(\tau)$ *contains the residual errors resulting from the least-squares approximation of* $\mathfrak{q}_i(\tau)$ *by* $\mathfrak{h}_i(\tau)$, *and* $\mathfrak{A}_i = \int_{\mathcal{K}_i} \varpi(\tau)\mathfrak{q}_i(\tau)\mathfrak{h}_i^\top(\tau)\mathsf{d}\tau\mathfrak{H}_i^{-1} \in \mathbb{R}^{\lambda_i\times l_i}$ *contains the coefficients of approximation, and* $\mathfrak{E}_i =$ $\int_{\mathcal{K}_i} \varpi(\tau)\mathfrak{e}_i(\tau)\mathfrak{e}_i^\top(\tau)\mathsf{d}\tau \succ 0$ *measures the level of errors in the case of* $\lambda_i \in \mathbb{N}$. *For any case of* $\lambda_i = 0, i = 1, \cdots, \nu$, *we define* $\mathfrak{E}_i = [\,]_{0,0}$ *as an empty matrix and* $[\,]_{0,0}^{-1} = [\,]_{0,0}$.

*Proof.* By applying [162, Eq. (17)] $\nu$ times with Lemma 3, the first inequality in (B.2) is established. Note that the definition of $\mathfrak{H}_i$ defined here corresponds to the inverse of the analogous matrix defined in [162]. Furthermore, $\mathfrak{E}_i^{-1} \succ 0$ in (B.2) is ensured by [162, Eq. (18)], which also implies the second inequality in (B.2). ■

## Appendix C. Proof of Theorem 1

*Proof.* The proof is divided into two parts, and is based on the construction of the complete-type KF

$$\mathsf{v}(\mathbf{x}_t(\cdot)) = \boldsymbol{\eta}^\top(t)\begin{bmatrix} P_1 & P_2 \\ * & P_3 \end{bmatrix}\boldsymbol{\eta}(t) + \sum_{i=1}^{\nu}\int_{\mathcal{I}_i} \boldsymbol{x}^\top(t+\tau)\left[Q_i + (\tau + r_i)R_i\right]\boldsymbol{x}(t+\tau)\mathsf{d}\tau \tag{C.1a}$$

in conjunction with the use of Lemma 4 and the KSD representation in Section 3.2, where $\mathbf{x}_t(\cdot)$ follows the definition in (35), and $P_1 \in \mathbb{S}^n$, $P_2 \in \mathbb{R}^{n\times dn}$, $P_3 \in \mathbb{S}^{dn}$ and $Q_i \in \mathbb{S}^n$, $R_i \in \mathbb{S}^n$ and

$$\begin{aligned} \boldsymbol{\eta}(t) &= \begin{bmatrix} \boldsymbol{x}^\top(t) & \mathfrak{z}^\top(t) \end{bmatrix}^\top, \\ \mathfrak{z}(t) &= \left[\int_{\mathcal{I}_i} \left(\mathfrak{L}_{\boldsymbol{f}_i}^{-1}\boldsymbol{f}_i(\tau) \otimes I_n\right)\boldsymbol{x}(t+\tau)\mathsf{d}\tau\right]_{i=1}^{\nu} \end{aligned} \tag{C.1b}$$

with $\mathfrak{L}_{\boldsymbol{f}_i}^{-1}$ and $\boldsymbol{f}_i(\cdot)$ in Proposition 1, $\forall i \in \mathbb{N}_\nu$. $\mathfrak{L}_{\boldsymbol{f}_i}^{-1}$ are well defined and unique by (6c).

### Part I:

We first establish the following identities on the integrals in (C.1a) and (C.1b). Utilizing the Leibniz integral rule [230] in view of $\boldsymbol{\chi}(t,\theta) = [\boldsymbol{x}(t + \acute{r}_i\theta - r_{i-1})]_{i=1}^{\nu} \in \mathbb{R}^{\nu n}$, we find that

$$\begin{aligned} &\sum_{i=1}^{\nu}\frac{\mathsf{d}}{\mathsf{d}t}\int_{\mathcal{I}_i} \boldsymbol{x}^\top(t+\tau)\ \left[Q_i + (r_i+\tau)R_i\right]\boldsymbol{x}(t+\tau)\mathsf{d}\tau \\ &= \sum_{i=1}^{\nu}[*]\ (Q_i + \acute{r}_iR_i)\,\boldsymbol{x}(t - r_{i-1}) - \sum_{i=1}^{\nu}[*]\,Q_i\boldsymbol{x}(t-r_i) \\ &\quad - \sum_{i=1}^{\nu}\int_{\mathcal{I}_i}\boldsymbol{x}^\top(t+\tau)R_i\boldsymbol{x}(t+\tau)\mathsf{d}\tau \\ &= \boldsymbol{\chi}^\top(t,0)\,(\mathbf{Q} + \mathbf{R}\Lambda)\,\boldsymbol{\chi}(t,0) - \boldsymbol{\chi}^\top(t,-1)\,\mathbf{Q}\,\boldsymbol{\chi}(t,-1) \\ &\quad - \sum_{i=1}^{\nu}\int_{\mathcal{I}_i}\boldsymbol{x}^\top(t+\tau)R_i\,\boldsymbol{x}(t+\tau)\mathsf{d}\tau \end{aligned} \tag{C.2}$$

with $\boldsymbol{\chi}(t,-1) = [\boldsymbol{x}(t-r_i)]_{i=1}^{\nu}$, $\boldsymbol{\chi}(t,0) = [\boldsymbol{x}(t-r_{i-1})]_{i=1}^{\nu}$, where $\mathbf{Q} = \operatorname{diag}_{i=1}^{\nu} Q_i$, and $\mathbf{R} = \operatorname{diag}_{i=1}^{\nu} R_i$ with $\Lambda$ in (39e). By Lemma 3 and the relation $\boldsymbol{f}_i(\tau) = \begin{bmatrix}\mathsf{O}_{d_i,\delta_i} & I_{d_i}\end{bmatrix}\boldsymbol{h}_i(\tau)$ with $\boldsymbol{f}_i(\cdot), \boldsymbol{h}_i(\cdot)$ in (6b), function $\mathfrak{z}(t)$ in (C.1b) satisfies

$$\begin{aligned} \mathfrak{z}(t) &= \left[\int_{\mathcal{I}_i}\left[\mathfrak{L}_{\boldsymbol{f}_i}^{-1}\boldsymbol{f}_i(\tau)\otimes I_n\right]\boldsymbol{x}(t+\tau)\mathsf{d}\tau\right]_{i=1}^{\nu} \\ &= \left[\left(\operatorname*{diag}_{i=1}^{\nu}\mathfrak{L}_{\boldsymbol{f}_i}^{-1}\right)\otimes I_n\right]\left[\int_{\mathcal{I}_i}\left(\boldsymbol{f}_i(\tau)\otimes I_n\right)\boldsymbol{x}(t+\tau)\mathsf{d}\tau\right]_{i=1}^{\nu} \\ &= \left[\left(\operatorname*{diag}_{i=1}^{\nu}\mathfrak{L}_{\boldsymbol{f}_i}^{-1}\begin{bmatrix}\mathsf{O}_{d_i,\delta_i} & I_{d_i}\end{bmatrix}\right)\otimes I_n\right] \\ &\qquad\qquad \times\left[\int_{\mathcal{I}_i}\left(\boldsymbol{h}_i(\tau)\otimes I_n\right)\boldsymbol{x}(t+\tau)\mathsf{d}\tau\right]_{i=1}^{\nu} \\ &= \left[\left(\operatorname*{diag}_{i=1}^{\nu}\mathfrak{L}_{\boldsymbol{f}_i}^{-1}\begin{bmatrix}\mathsf{O}_{d_i,\delta_i} & I_{d_i}\end{bmatrix}\mathfrak{L}_{\boldsymbol{h}_i}\right)\otimes I_n\right]\boldsymbol{\xi}(t) = \widehat{I}\boldsymbol{\xi}(t) \end{aligned} \tag{C.3a}$$

with $\boldsymbol{\xi}(t)$ in (31), and its weak derivative satisfies

$$\begin{aligned} \frac{\mathsf{d}\mathfrak{z}(t)}{\mathsf{d}t} &= \frac{\mathsf{d}}{\mathsf{d}t}\left[\int_{\mathcal{I}_i}\left(\mathfrak{L}_{\boldsymbol{f}_i}^{-1}\otimes I_n\right)F_i(\tau)\boldsymbol{x}(t+\tau)\mathsf{d}\tau\right]_{i=1}^{\nu} \\ &= \left[\int_{\mathcal{I}_i}\left(\mathfrak{L}_{\boldsymbol{f}_i}^{-1}\otimes I_n\right)F_i(\tau)\frac{\partial}{\partial t}\boldsymbol{x}(t+\tau)\mathsf{d}\tau\right]_{i=1}^{\nu} \\ &= \left[\int_{\mathcal{I}_i}\left(\mathfrak{L}_{\boldsymbol{f}_i}^{-1}\otimes I_n\right)F_i(\tau)\frac{\partial}{\partial \tau}\boldsymbol{x}(t+\tau)\mathsf{d}\tau\right]_{i=1}^{\nu} \\ &= \left[\left(\operatorname*{diag}_{i=1}^{\nu}\mathfrak{L}_{\boldsymbol{f}_i}^{-1}\right)\otimes I_n\right]\left[\operatorname*{diag}_{i=1}^{\nu}F_i(-r_{i-1})\right]\boldsymbol{\chi}(t,0) \\ &\quad - \left[\left(\operatorname*{diag}_{i=1}^{\nu}\mathfrak{L}_{\boldsymbol{f}_i}^{-1}\right)\otimes I_n\right]\left[\operatorname*{diag}_{i=1}^{\nu}F_i(-r_i)\right]\boldsymbol{\chi}(t,-1) \\ &- \left[\left(\operatorname*{diag}_{i=1}^{\nu}\mathfrak{L}_{\boldsymbol{f}_i}^{-1}\right)\otimes I_n\right]\left[\int_{\mathcal{I}_i}\left(\frac{\mathsf{d}}{\mathsf{d}\tau}\boldsymbol{f}_i(\tau)\otimes I_n\right)\boldsymbol{x}(t+\tau)\mathsf{d}\tau\right]_{i=1}^{\nu} \\ &= \left[\left(\operatorname*{diag}_{i=1}^{\nu}\mathfrak{L}_{\boldsymbol{f}_i}^{-1}\boldsymbol{f}_i(-r_{i-1})\right)\otimes I_n\right]\boldsymbol{\chi}(t,0) \\ &\quad - \left[\left(\operatorname*{diag}_{i=1}^{\nu}\mathfrak{L}_{\boldsymbol{f}_i}^{-1}\boldsymbol{f}_i(-r_i)\right)\otimes I_n\right]\boldsymbol{\chi}(t,-1) \\ &\quad - \left[\left(\operatorname*{diag}_{i=1}^{\nu}\mathfrak{L}_{\boldsymbol{f}_i}^{-1}\right)\otimes I_n\right]\left[\left(\operatorname*{diag}_{i=1}^{\nu}M_i\mathfrak{L}_{\boldsymbol{h}_i}\right)\otimes I_n\right]\boldsymbol{\xi}(t) \\ &= \left[\left(\operatorname*{diag}_{i=1}^{\nu}\mathfrak{L}_{\boldsymbol{f}_i}^{-1}\boldsymbol{f}_i(-r_{i-1})\right)\otimes I_n \quad \mathsf{O}_{dn} \quad \mathsf{O}_{dn,\varkappa n}\right]\boldsymbol{\omega}(t) \\ &\quad - \left[\mathsf{O}_{dn} \quad \left(\operatorname*{diag}_{i=1}^{\nu}\mathfrak{L}_{\boldsymbol{f}_i}^{-1}\boldsymbol{f}_i(-r_i)\right)\otimes I_n \cdots\right. \\ &\qquad\qquad \left.\cdots\left(\operatorname*{diag}_{i=1}^{\nu}\mathfrak{L}_{\boldsymbol{f}_i}^{-1}M_i\mathfrak{L}_{\boldsymbol{h}_i}\right)\otimes I_n\right]\boldsymbol{\omega}(t) \\ &= \left(\mathbf{M}\otimes I_n\right)\boldsymbol{\omega}(t) \end{aligned} \tag{C.3b}$$

with $\mathfrak{H}_i$ in (8), $\boldsymbol{\omega}(t)$ in (33e) and $\mathbf{M}$ in (39f), where $F_i(\tau) = \boldsymbol{f}_i(\tau)\otimes I_n$. Specifically, the first and second equalities in (C.3b) were established by the Leibniz rule [230] for Lebesgue integrals in light of $\frac{\partial}{\partial t}\boldsymbol{x}(t+\tau) = \frac{\partial}{\partial \tau}\boldsymbol{x}(t+\tau)$ (weak derivative), and by using the integration by parts formula together with

$$\begin{aligned} &\left[\int_{\mathcal{I}_i}\left(\mathfrak{L}_{\boldsymbol{f}_i}^{-1}\frac{\mathsf{d}}{\mathsf{d}\tau}\boldsymbol{f}_i(\tau)\otimes I_n\right)\boldsymbol{x}(t+\tau)\mathsf{d}\tau\right]_{i=1}^{\nu} \\ &\qquad = \left[\left(\operatorname*{diag}_{i=1}^{\nu}\mathfrak{L}_{\boldsymbol{f}_i}^{-1}M_i\mathfrak{L}_{\boldsymbol{h}_i}\right)\otimes I_n\right]\times \\ &\left[\int_{\mathcal{I}_i}\left(\mathfrak{L}_{\boldsymbol{h}_i}^{-1}\boldsymbol{h}_i(\tau)\otimes I_n\right)\boldsymbol{x}(t+\tau)\mathsf{d}\tau\right]_{i=1}^{\nu} \end{aligned}$$

$$
= \left[ \left( \mathop{\mathsf{diag}}_{i=1}^{\nu} \mathfrak{L}_{\boldsymbol{f}_i}^{-1} M_i \mathfrak{L}_{\boldsymbol{h}_i} \right) \otimes I_n \right] \boldsymbol{\xi}(t), \tag{C.4a}
$$

$$
\begin{bmatrix} \boldsymbol{x}(t) \\ \boldsymbol{\chi}(t,-1) \end{bmatrix} = \begin{bmatrix} \boldsymbol{\chi}(t,0) \\ \boldsymbol{x}(t-r_\nu) \end{bmatrix} = \left[ \boldsymbol{x}(t-r_i) \right]_{i=0}^{\nu} \tag{C.4b}
$$

following from the relation in (6b) with $\boldsymbol{\xi}(t)$ as in (31), given that $\left(\mathfrak{L}_{\boldsymbol{f}_i}^{-1} \otimes I_n\right) F_i(\tau) = \mathfrak{L}_{\boldsymbol{f}_i}^{-1} \boldsymbol{f}_i(\tau) \otimes I_n$.

Taking into account the identities in (C.2)–(C.4b), we now differentiate (weak derivative) $\mathsf{v}(\mathbf{x}_t(\cdot))$ in (C.1a) along the trajectory of $\dot{\boldsymbol{x}}(t)$, $\boldsymbol{z}(t)$ in (32) and consider $\mathsf{s}(\boldsymbol{z}(t), \boldsymbol{w}(t))$ in (37). Then, it results in

$$
\begin{aligned}
&\widetilde{\forall} t \geq t_0, \dot{\mathsf{v}}(\mathbf{x}_t(\cdot)) - \mathsf{s}(\boldsymbol{z}(t), \boldsymbol{w}(t)) = \\
&\mathsf{Sy}\left( \boldsymbol{\eta}^\top(t) \begin{bmatrix} P_1 & P_2 \\ * & P_3 \end{bmatrix} \dot{\boldsymbol{\eta}}(t) \right) - \begin{bmatrix} \boldsymbol{z}(t) \\ \boldsymbol{w}(t) \end{bmatrix}^\top \begin{bmatrix} \widetilde{J}^\top J_1^{-1} \widetilde{J} & J_2 \\ * & J_3 \end{bmatrix} \begin{bmatrix} \boldsymbol{z}(t) \\ \boldsymbol{w}(t) \end{bmatrix} \\
&= \boldsymbol{\vartheta}^\top(t)\, \mathsf{Sy}\left( S^\top \begin{bmatrix} P_1 & P_2 \\ * & P_3 \end{bmatrix} \begin{bmatrix} \mathbf{A} + \mathbf{B}_1 \left[ (I_\beta \otimes K) \oplus \mathsf{O}_q \right] \\ \left[ \mathbf{M} \otimes I_n \quad \mathsf{O}_{dn,(\mu n+q)} \right] \end{bmatrix} \right. \\
&\quad \left. - \begin{bmatrix} \mathsf{O}_{(\beta n),m} \\ J_2^\top \end{bmatrix} \boldsymbol{\Sigma} \right) \boldsymbol{\vartheta}(t) + \boldsymbol{\chi}^\top(t,0) \left( \mathbf{Q} + \mathbf{R}\Lambda \right) \boldsymbol{\chi}(t,0) \\
&\quad - \boldsymbol{\chi}^\top(t,-1)\, \mathbf{Q} \boldsymbol{\chi}(t,-1) - \sum_{i=1}^{\nu} \int_{\mathcal{I}_i} \boldsymbol{x}^\top(t+\tau) R_i \boldsymbol{x}(t+\tau) \mathsf{d}\tau \\
&\quad - \boldsymbol{w}^\top(t) J_3 \boldsymbol{w}(t) - \boldsymbol{\vartheta}^\top(t) \boldsymbol{\Sigma}^\top \widetilde{J}^\top J_1^{-1} \widetilde{J} \boldsymbol{\Sigma} \boldsymbol{\vartheta}(t)
\end{aligned} \tag{C.5}
$$

with $\boldsymbol{\vartheta}(t)$ in (33f) and $\mathbf{A}, \mathbf{B}_1, K$ in (33a)–(33d), where $\boldsymbol{\Sigma}, S, \widehat{I}$ are given by the statements of Theorem 1. Note that $S, \widehat{I}$ in (39b)–(39d) are constructed as a result of the identity $\boldsymbol{\eta}(t) = S\boldsymbol{\vartheta}(t)$ which follows from $\mathfrak{z}(t) = \widehat{I}\boldsymbol{\xi}(t)$ in (C.3a) given the fact that vectors $\boldsymbol{x}(t)$, $\boldsymbol{\xi}(t)$ are part of $\boldsymbol{\vartheta}(t)$ and vectors $\boldsymbol{\xi}(t)$, $\mathfrak{z}(t)$ constitute $\boldsymbol{\eta}(t)$ in (C.1b).

Assuming that (38b) holds, we apply (B.2) with $\varpi(\tau) = 1$, $\mathcal{K}_i = \mathcal{I}_i$, $X_i = R_i \succ 0$, and the parameters in (8)–(10) such that $\mathfrak{q}_i(\tau) = \boldsymbol{\phi}_i(\tau)$, $\mathfrak{h}_i(\tau) = \boldsymbol{h}_i(\tau)$, $i \in \mathbb{N}_\nu$ to the summation $\sum_{i=1}^{\nu} \int_{\mathcal{I}_i} \boldsymbol{x}^\top(t+\tau) R_i \boldsymbol{x}(t+\tau) \mathsf{d}\tau$ in (C.5), which produces

$$
\begin{aligned}
&\sum_{i=1}^{\nu} \int_{\mathcal{I}_i} \boldsymbol{x}^\top(t+\tau) R_i \boldsymbol{x}(t+\tau) \mathsf{d}\tau \geq \\
&\boldsymbol{\xi}^\top(t) \left[ \mathop{\mathsf{diag}}_{i=1}^{\nu} I_{\varkappa_i} \otimes R_i \right] \boldsymbol{\xi}(t) + \boldsymbol{e}^\top(t) \left[ \mathop{\mathsf{diag}}_{i=1}^{\nu} I_{\mu_i} \otimes R_i \right] \boldsymbol{e}(t)
\end{aligned} \tag{C.6}
$$

with $\boldsymbol{\xi}(t), \boldsymbol{e}(t)$ in (31). We note that matrices $\mathfrak{H}_i^{-1} \succ 0$ and $\mathfrak{E}_i^{-1} \succ 0$ in (C.6) were factored using the Cholesky decompositions that are concealed inside $\boldsymbol{\xi}(t), \boldsymbol{e}(t)$, respectively, as a consequence of the identities $\mathfrak{H}_i^{-1} \otimes R_i = [*] \left(I_{\varkappa_i} \otimes R_i\right) \left(\mathfrak{L}_{\boldsymbol{h}_i}^{-1} \otimes I_n\right)$ and $\mathfrak{E}_i^{-1} \otimes R_i = [*] \left(I_{\mu_i} \otimes R_i\right) \left(\mathfrak{L}_{\boldsymbol{\varepsilon}_i}^{-1} \otimes I_n\right)$ by (A.1). Then utilizing (C.6) on (C.5) yields

$$
\begin{aligned}
&\widetilde{\forall} t \geq t_0, \quad \dot{\mathsf{v}}(\mathbf{x}_t(\cdot)) - \mathsf{s}(\boldsymbol{z}(t), \boldsymbol{w}(t)) \\
&\qquad \leq \boldsymbol{\vartheta}^\top(t) \left( \boldsymbol{\Psi} - \boldsymbol{\Sigma}^\top \widetilde{J}^\top J_1^{-1} \widetilde{J} \boldsymbol{\Sigma} \right) \boldsymbol{\vartheta}(t)
\end{aligned}
$$

with $\boldsymbol{\Psi}$ in (39a) and $\boldsymbol{\vartheta}(t)$ in (33f). Now it is evident that if (38b) and $\boldsymbol{\Psi} - \boldsymbol{\Sigma}^\top \widetilde{J}^\top J_1^{-1} \widetilde{J} \boldsymbol{\Sigma} \prec 0$ are true, then

$$
\exists \epsilon_3 > 0 : \widetilde{\forall} t \geq t_0, \dot{\mathsf{v}}(\mathbf{x}_t(\cdot)) - \mathsf{s}(\boldsymbol{z}(t), \boldsymbol{w}(t)) \leq -\epsilon_3 \left\| \boldsymbol{x}(t) \right\|^2. \tag{C.7}
$$

Inspecting the internal structure of $\boldsymbol{\Psi} - \boldsymbol{\Sigma}^\top \widetilde{J}^\top J_1^{-1} \widetilde{J} \boldsymbol{\Sigma}$ in view of $\boldsymbol{\vartheta}(t)$ as in (33f), it follows from (C.7) that

$$
\exists \epsilon_3 > 0, \ \widetilde{\forall} t \geq t_0, \ \dot{\mathsf{v}}(\mathbf{x}_t(\cdot)) \leq -\epsilon_3 \left\| \boldsymbol{x}(t) \right\|^2 \tag{C.8}
$$

for any $\boldsymbol{\psi}(\cdot) \in \mathcal{C}(\mathcal{J}; \mathbb{R}^n)$ in (32) with $\boldsymbol{w}(t) \equiv \mathbf{0}_q$. Note that $\mathbf{x}_t(\cdot)$ and $\boldsymbol{x}(t)$ in (C.8) follow the definition in (34b). Thus, if (38b) and $\boldsymbol{\Psi} - \boldsymbol{\Sigma}^\top \widetilde{J}^\top J_1^{-1} \widetilde{J} \boldsymbol{\Sigma} \prec 0$ are feasible, then there exists $\mathsf{v}(\mathbf{x}_t(\cdot))$ in (C.1a) satisfying (34b)–(35). Finally, applying the Schur complement to $\boldsymbol{\Psi} - \boldsymbol{\Sigma}^\top \widetilde{J}^\top J_1^{-1} \widetilde{J} \boldsymbol{\Sigma} \prec 0$ with $J_1^{-1} \prec 0$ produces (38c). Hence we have proven that if (38b)–(38c) are feasible, then there exists $\epsilon_3 > 0$ such that $\mathsf{v}(\mathbf{x}_t(\cdot))$ in (C.1a) satisfies (34b)–(35).

### Part II:

In what follows, we prove that there exist $\epsilon_1, \epsilon_2 > 0$ such that the KF in (C.1a) satisfies (34a) if (38a)–(38b) are feasible. Consider $\mathsf{v}(\mathbf{x}_t(\cdot))$ in (C.1a) at $t = t_0$, we get

$$
\begin{aligned}
&\exists \lambda > 0, \ \forall \boldsymbol{\psi}(\cdot) \in \mathcal{C}\left(\mathcal{J}; \mathbb{R}^n\right), \ \mathsf{v}(\mathbf{x}_{t_0}(\cdot)) = \mathsf{v}(\boldsymbol{\psi}(\cdot)) \\
&\leq \boldsymbol{\eta}^\top(t_0) \lambda \boldsymbol{\eta}(t_0) + \underbrace{\int_{-r}^{0} \boldsymbol{\psi}^\top(\tau) \lambda \boldsymbol{\psi}(\tau) \mathsf{d}\tau}_{\leq \lambda r \|\boldsymbol{\psi}(\cdot)\|_\infty^2} \\
&\leq \lambda \left\| \boldsymbol{\psi}(0) \right\|^2 + \lambda r \left\| \boldsymbol{\psi}(\cdot) \right\|_\infty^2 + \mathfrak{z}^\top(t_0) \left[ \mathop{\mathsf{diag}}_{i=1}^{\nu} I_{d_i} \otimes \lambda I_n \right] \mathfrak{z}(t_0) \\
&\leq \lambda \left\| \boldsymbol{\psi}(0) \right\|^2 + \lambda r \left\| \boldsymbol{\psi}(\cdot) \right\|_\infty^2 + \sum_{i=1}^{\nu} \int_{\mathcal{I}_i} \boldsymbol{\psi}^\top(\tau) \lambda \boldsymbol{\psi}(\tau) \mathsf{d}\tau \\
&\leq (\lambda + \lambda r) \left\| \boldsymbol{\psi}(\cdot) \right\|_\infty^2 + \lambda \int_{-r}^{0} \boldsymbol{\psi}^\top(\tau) \boldsymbol{\psi}(\tau) \mathsf{d}\tau \\
&\leq (\lambda + 2\lambda r) \left\| \boldsymbol{\psi}(\cdot) \right\|_\infty^2
\end{aligned} \tag{C.9}
$$

for any $\boldsymbol{\psi}(\cdot) \in \mathcal{C}\left(\mathcal{J}; \mathbb{R}^n\right)$ in (32), where the sequence of upper bound estimates in (C.9) follows from the property of quadratic forms $\forall X \in \mathbb{S}^n, \exists \lambda > 0 : \forall \mathbf{x} \in \mathbb{R}^n \setminus \{\mathbf{0}\}, \mathbf{x}^\top (\lambda I_n - X) \mathbf{x} > 0$ and an application of the second inequality in (B.2) with $\varpi(\tau) = 1$, $\mathcal{I}_i = \mathcal{K}_i$, $\mathfrak{h}_i(\tau) = \mathfrak{L}_{\boldsymbol{f}_i}^{-1} \boldsymbol{f}_i(\tau)$ and $X_i = \lambda I_n$. Next, applying the second inequality in (B.2) again to (C.1a) with $\varpi(\tau) = 1$, $\mathfrak{h}_i(\tau) = \mathfrak{L}_{\boldsymbol{f}_i}^{-1} \boldsymbol{f}_i(\tau)$ and $X_i = Q_i$ produces

$$
\begin{aligned}
&\sum_{i=1}^{\nu} \int_{\mathcal{I}_i} \boldsymbol{x}^\top(t+\tau) Q_i \boldsymbol{x}(t+\tau) \mathsf{d}\tau \\
&\geq \sum_{i=1}^{\nu} [*] \left( I_{d_i} \otimes Q_i \right) \int_{\mathcal{I}_i} \left( \mathfrak{L}_{\boldsymbol{f}_i}^{-1} \boldsymbol{f}_i(\tau) \otimes I_n \right) \boldsymbol{x}(t+\tau) \mathsf{d}\tau \\
&= \mathfrak{z}^\top(t) \left( \mathop{\mathsf{diag}}_{i=1}^{\nu} I_{d_i} \otimes Q_i \right) \mathfrak{z}(t)
\end{aligned} \tag{C.10}
$$

if (38b) is true. Moreover, by utilizing (C.10) on (C.1a) with $t = t_0$ and $\boldsymbol{x}(t_0 + \theta) = \boldsymbol{\psi}(\theta)$ in (32) and considering (38b) and the results in (C.9), we find that $\mathsf{v}(\bullet)$ in (C.1a) satisfies (34a) for some $\epsilon_1; \epsilon_2 > 0$ provided that the inequalities in (38a)–(38b) have feasible solutions.

In conclusion, we have shown that feasible solutions to (38a)–(38c) imply the existence of the KF in (C.1a) for some $\epsilon_1, \epsilon_2, \epsilon_3 > 0$ satisfying (35) and (34a)–(34b). ∎

**Remark 14.** The matrix inequality in (6c) guarantees the invertibility of matrices $\mathfrak{F}_i, \mathfrak{H}_i, \mathfrak{E}_i$. The fact that $\boldsymbol{f}_i(\cdot)$ in (21) and (C.1b) are identical illustrates our deliberate design of the KSD concept so that it can coordinate perfectly with the use of KF (C.1a). Given the proofs in Appendix C, we can clearly see that functions $\boldsymbol{f}_i(\cdot)$ in $\boldsymbol{g}_i(\cdot)$ in Proposition 1 are set apart as $\mathscr{W}^{1,2}$ functions so that $\boldsymbol{f}_i(\cdot)$ can be differentiated (weak) as the integral kernels of KF (C.1a). This enables us to utilize KF (C.1a) with $\mathscr{W}^{1,2}$ integral kernels which in turn can contain any number of WDLIFs even if none of the functions in $\boldsymbol{f}_i(\cdot)$ are included by any matrix in (20). Thus, our KF (C.1a) is significantly more general than those[2] in [144, 152, 162] even for the case of $\nu = 1$, as the generality of KFs is dominated by the generality of their integral kernels. Additionally, inequality (B.2) used in the proof of Theorem 1 is derived based on the least-squares principle to ensure minimal conservatism. Lastly, the relation in (6b) is specifically formulated to ensure that (C.5) is constructible in quadratic forms via $\boldsymbol{\vartheta}(t)$ in (33f) with $S$ in (39b) and $\mathbf{M}$ in (39f), considering the structure of $\mathfrak{z}(t)$ in (C.1b) and the expressions of $\frac{\mathsf{d}}{\mathsf{d}t}\mathfrak{z}(t)$ in (C.3b).

## Appendix D. Some Discussions on Conservatism

We present a qualitative analysis of the conservatism of our proposed method based on the properties of several relevant mathematical constructs. We show that the form of KF (C.1a) is a parameterization of the complete Krasovskiĭ functional [118] of LTDS $\dot{\boldsymbol{x}}(t) = \sum_{i=0}^{\nu} A_i \boldsymbol{x}(t - r_i)$ via basis functions $\boldsymbol{f}_i(\cdot)$. By analyzing various properties with appropriate comparisons, we demonstrate why the effectiveness of our approach is ensured and cannot be more conservative than those of the existing methods in terms of stability analysis.

Research has shown that the trivial solution to $\dot{\boldsymbol{x}}(t) = \sum_{i=0}^{\nu} A_i \boldsymbol{x}(t - r_i)$ is exponentially stable if and only if there exists a complete-type KF

$$\begin{aligned}\mathsf{v}(\mathbf{x}_t(\cdot)) &= \boldsymbol{x}^\top(t) P_1 \boldsymbol{x}(t) + 2\sum_{i=1}^{\nu} \boldsymbol{x}^\top(t) \int_{-r_i}^{0} \widetilde{P}_i^2(\tau) \boldsymbol{x}(t+\tau)\mathsf{d}\tau \\ &+ \sum_{i=1}^{\nu}\sum_{j=1}^{\nu} \int_{-r_j}^{0}\int_{-r_i}^{0} \boldsymbol{x}^\top(t+\theta) \widetilde{P}_{i,j}^3(\tau,\theta) \boldsymbol{x}(t+\tau)\mathsf{d}\tau\mathsf{d}\theta \\ &+ \sum_{i=1}^{\nu} \int_{-r_i}^{0} \boldsymbol{x}^\top(t+\tau)\left(Q_i + (\tau + r_i)R_i\right)\boldsymbol{x}(t+\tau)\mathsf{d}\tau \end{aligned} \tag{D.1}$$

satisfying conditions analogous to (34a)–(34b). Given the basic properties of integrals, (D.1) is equivalent to

$$\begin{aligned}\mathsf{v}(\mathbf{x}_t(\cdot)) &= \boldsymbol{x}^\top(t) P_1 \boldsymbol{x}(t) + 2\sum_{i=1}^{\nu} \boldsymbol{x}^\top(t) \int_{\mathcal{I}_i} \widetilde{P}_i^2(\tau) \boldsymbol{x}(t+\tau)\mathsf{d}\tau \\ &+ \sum_{i=1}^{\nu}\sum_{j=1}^{\nu} \int_{\mathcal{I}_i}\int_{\mathcal{I}_j} \boldsymbol{x}^\top(t+\theta)\, \widetilde{P}_{i,j}^3(\tau,\theta)\, \boldsymbol{x}(t+\tau)\mathsf{d}\tau\mathsf{d}\theta \\ &+ \sum_{i=1}^{\nu} \int_{\mathcal{I}_i} \boldsymbol{x}^\top(t+\tau)\left(Q_i + (\tau + r_i)R_i\right)\boldsymbol{x}(t+\tau)\mathsf{d}\tau \end{aligned} \tag{D.2}$$

where $P_1 \in \mathbb{S}^n$, $\widetilde{P}_i^2(\cdot) \in \mathscr{W}^{1,2}(\mathcal{I}_i; \mathbb{R}^{n\times n})$ and $\widetilde{P}_{i,j}^3(\cdot,\cdot) \in \mathscr{W}^{1,2}(\mathcal{J}\times\mathcal{J}; \mathbb{R}^{n\times n})$ satisfying $\widetilde{P}_{i,j}^3(\tau,\theta) = \widetilde{P}_{j,i}^{3,\top}(\theta,\tau)$. Now $\widetilde{P}_i^2(\cdot)$ and $\widetilde{P}_{i,j}^3(\cdot,\cdot)$ are variables with infinite dimension that are difficult to obtain numerically. To address this issue, we can employ a subspace parameterization

$$\begin{aligned}\widetilde{P}_i^2(\tau) &= \widehat{P}_i^2\left[\mathfrak{L}_{\boldsymbol{f}_i}^{-1}\boldsymbol{f}_i(\tau)\otimes I_n\right], \quad \widehat{P}_i^2 \in \mathbb{R}^{n\times d_i n} \\ \widetilde{P}_{i,j}^3(\tau,\theta) &= \left[\mathfrak{L}_{\boldsymbol{f}_j}^{-1}\boldsymbol{f}_j(\theta)\otimes I_n\right]^\top \widehat{P}_{i,j}^3\left[\mathfrak{L}_{\boldsymbol{f}_i}^{-1}\boldsymbol{f}_i(\tau)\otimes I_n\right]\end{aligned} \tag{D.3}$$

using normalized functions $\mathfrak{L}_{\boldsymbol{f}_i}^{-1}\boldsymbol{f}_i(\tau)$ with $\mathfrak{L}_{\boldsymbol{f}_i}$ and $\boldsymbol{f}_i(\cdot) \in \mathscr{W}^{1,2}(\mathcal{I}_i; \mathbb{R}^{d_i})$ in Theorem 1, where $\widehat{P}_{i,j}^3 \in \mathbb{R}^{d_j n\times d_i n}$ satisfies $\widehat{P}_{i,j}^3 = \widehat{P}_{j,i}^{3,\top}$. Consequently, CKF (D.2) becomes our KF in (C.1a) if we set $P_2 = \left\llbracket \widehat{P}_i^2 \right\rrbracket_{i=1}^{\nu}$ and $P_3 = \left\llbracket \left[\widehat{P}_{i,j}^3\right]_{i=1}^{\nu}\right\rrbracket_{j=1}^{\nu}$ in (D.2) because of relations

$$\begin{aligned}&\sum_{i=1}^{\nu} \boldsymbol{x}^\top(t) \int_{\mathcal{I}_i} \widetilde{P}_i^2(\tau)\boldsymbol{x}(t+\tau)\mathsf{d}\tau \\ &= \boldsymbol{x}^\top(t) \sum_{i=1}^{\nu} \widehat{P}_i^2 \int_{\mathcal{I}_i}\left[\mathfrak{L}_{\boldsymbol{f}_i}^{-1}\boldsymbol{f}_i(\tau)\otimes I_n\right]\boldsymbol{x}(t+\tau)\mathsf{d}\tau \\ &= \boldsymbol{x}^\top(t)\left\llbracket \widehat{P}_i^2 \right\rrbracket_{i=1}^{\nu} \mathfrak{z}(t) = \boldsymbol{x}^\top(t) P_2\, \mathfrak{z}(t) = \mathfrak{z}^\top(t) P_2^\top \boldsymbol{x}(t)\end{aligned}$$

and

$$\begin{aligned}&\sum_{i=1}^{\nu}\sum_{j=1}^{\nu}\int_{\mathcal{I}_i}\int_{\mathcal{I}_j} \boldsymbol{x}^\top(t+\theta)\widetilde{P}_{i,j}^3(\tau,\theta)\boldsymbol{x}(t+\tau)\mathsf{d}\tau\mathsf{d}\theta \\ &= \sum_{i=1}^{\nu}\sum_{j=1}^{\nu}\int_{\mathcal{I}_j} \boldsymbol{x}^\top(t+\theta)\left[\mathfrak{L}_{\boldsymbol{f}_j}^{-1}\boldsymbol{f}_j(\theta)\otimes I_n\right]^\top \mathsf{d}\theta\, \widehat{P}_{i,j}^3 \times \\ &\int_{\mathcal{I}_i}\left[\mathfrak{L}_{\boldsymbol{f}_i}^{-1}\boldsymbol{f}_i(\tau)\otimes I_n\right]\boldsymbol{x}(t+\tau)\mathsf{d}\tau = \mathfrak{z}^\top(t)\left\llbracket\left[\widehat{P}_{i,j}^3\right]_{i=1}^{\nu}\right\rrbracket_{j=1}^{\nu}\mathfrak{z}(t)\end{aligned}$$

with $\mathfrak{z}(t)$ in (C.1b), following from the definition of quadratic forms and identity $\left[\int_{\mathcal{I}_i} F_i(\tau)\boldsymbol{x}(t+\tau)\mathsf{d}\tau\right]_{i=1}^{\nu} = \left[\int_{\mathcal{I}_j} F_j(\theta)\boldsymbol{x}(t+\theta)\mathsf{d}\theta\right]_{j=1}^{\nu}$.

Consequently, our KF in (C.1a) can be regarded as a particular realization of the CKF in (D.2) using subspace parameterization via $\mathfrak{L}_{\boldsymbol{f}_i}^{-1}\boldsymbol{f}_i(\tau) = \left[\widetilde{f}_{ij}(\tau)\right]_{j=1}^{d_i}$. Let $\boldsymbol{f}_i(\tau) = \left[f_{ij}(\tau)\right]_{j=1}^{d_i}$ satisfying (6c), then by some fundamental properties of functional analysis [55], we have $\operatorname{span}_{j=1}^{d_i}\widetilde{f}_{ij}(\cdot) = \operatorname{span}_{j=1}^{d_i} f_{ij}(\cdot) \subset \operatorname{span}_{j=1}^{d_i+1} f_{ij}(\cdot) = \operatorname{span}_{j=1}^{d_i+1}\widetilde{f}_{ij}(\cdot)$. Therefore, we see that the generality of $\widetilde{P}_i^2(\tau)$, $\widetilde{P}_{i,j}^3(\tau,\theta)$ in (D.3) can be augmented indefinitely by adding an arbitrary finite number of new WDLIFs to $\boldsymbol{f}_i(\cdot)$, which in turn increases the generality

[2]The equivalent $\boldsymbol{f}(\cdot)$ functions in [144, 152, 162] are the special case of our $\boldsymbol{f}_i(\cdot)$ in this work.

of our KF in (C.1a).

The generality of our KF in (C.1a) can also be illustrated by comparing it to the KFs in the existing literature. The KF in [144] represents a special case of (C.1a) with $\nu = 1$ wherein $\boldsymbol{f}_1(\cdot)$ comprises Legendre polynomials. Likewise, the KFs in [152, 162] are specific instances of (C.1a) with $\nu = 1$ for the case of retarded systems, where all functions in $\boldsymbol{f}_1(\cdot)$ are the solutions to linear homogeneous differential equations with constant coefficients satisfying $\frac{\mathrm{d}\boldsymbol{f}_1(\tau)}{\mathrm{d}\tau} = N\boldsymbol{f}_1(\tau)$ for some $N \in \mathbb{R}^{d_1 \times d_1}$. The effectiveness of these methods has been demonstrated by the numerical experiments in [144, 152, 162] that these approaches can successfully detect the analytical delay margins of collections of challenging examples of LTDSs. Given that our KF and synthesis conditions generalize all the counterparts in [144, 152, 162], we can confidently claim that the effectiveness of our new method is guaranteed, as it cannot be more conservative than the above existing results.

On the other hand, we do not claim in this paper that $\boldsymbol{u}(t) = K\boldsymbol{x}(t)$ is capable of stabilizing (18) under any combination of state space parameters. Indeed, if (18) happens to be $\dot{\boldsymbol{x}}(t) = A\boldsymbol{x}(t) + B\boldsymbol{u}(t-r)$ in the absence of state delays, then it is preferable to employ a controller [231] that can fully eliminate the effects of input delay. Our main objective, however, is to address a system with multiple pointwise and general distributed delays simultaneously at the state, input and output, as described in (18), so that the methodology may accommodate a wide variety of applications. At the same time, the works by [161, 232] have shown that system $\dot{\boldsymbol{x}}(t) = A\boldsymbol{x}(t) + \int_{-r}^{0} B(\tau)\boldsymbol{x}(t+\tau)\mathrm{d}\tau$ can be stabilized by $\boldsymbol{u}(t) = K\boldsymbol{x}(t)$. Given that no predictor controllers have been proposed for systems with pointwise and general distributed delays at both the states and inputs, it makes $\boldsymbol{u}(t) = K\boldsymbol{x}(t)$ the most viable option at the moment for our problem. Furthermore, the use of $\boldsymbol{u}(t) = K\boldsymbol{x}(t)$ does offer distinct advantages over controllers with more complex structures, including lower implementation costs and greater simplicity, potentially leading to higher reliability in certain real-world applications. Therefore, it is meaningful to first formulate a solution to the DSFC problem of (18) using $\boldsymbol{u}(t) = K\boldsymbol{x}(t)$, as this solution can serve as a foundation for developing new methods in the future for the same control problem using controllers with more complex structures.

Finally, the investigated control problem differs considerably from the problem researched by [200] in several aspects; notably, systems with time-varying delays $r(t)$ necessitate fundamentally different analytical strategies. Since the expression of $r(t)$ is time-varying and unknown, the kernel functions in [200] could not be approximated over $[-r(t), 0]$, which makes the synthesis condition in [200] absent in $\mathfrak{C}_i$. Moreover, the strategies of constructing KFs for systems with an unknown $r(t)$ differ significantly from those employed when $r_i$ are constants. In fact, even treating $\boldsymbol{x}(t - r(t))$ and $\int_{-r(t)}^{0} D(\tau)\boldsymbol{x}(t+\tau)\mathrm{d}\tau$ together is fairly different from analyzing systems with $\int_{-r(t)}^{0} D(\tau)\boldsymbol{x}(t+\tau)\mathrm{d}\tau$ alone, as explained in [200] considering the use of integral inequalities in the contexts of constructing KFs. Thus, readers should not misinterpret the superficial similarity between (38c) and [200, Eq. (26)], as their internal structures are vastly different in light of the definitions of the matrices in Theorem 1.

## Appendix E. Proof of Theorem 2

*Proof.* First of all, (38c) can be rewritten as

$$\mathsf{Sy}\left(\mathbf{P}^\top \mathbf{\Pi}\right) + \mathbf{\Phi} = [*] \begin{bmatrix} \mathsf{O}_n & \mathbf{P} \\ * & \mathbf{\Phi} \end{bmatrix} \begin{bmatrix} \mathbf{\Pi} \\ I_{\beta n+q+m} \end{bmatrix} \prec 0. \tag{E.1}$$

whose structure is similar to those of the inequalities in Lemma 2. As there are two inequalities in (41b), a new matrix inequality must be constructed so that Lemma 2 can be utilized with (E.1). Observing that whenever (E.1) holds, its bottom-right principal sub-matrix yields

$$\Upsilon^\top \begin{bmatrix} \mathsf{O}_n & \mathbf{P} \\ * & \mathbf{\Phi} \end{bmatrix} \Upsilon = \begin{bmatrix} -J_3 - \mathsf{Sy}(D_2^\top J_2) & D_2^\top \widetilde{J} \\ * & J_1 \end{bmatrix} \prec 0 \tag{E.2}$$

where $\Upsilon^\top := \begin{bmatrix} \mathsf{O}_{(q+m),(n+\beta n)} & I_{q+m} \end{bmatrix}$. This is because the right-hand side in (E.2) is the $2 \times 2$ block matrix at the bottom-right corner of both $\mathsf{Sy}\left(\mathbf{P}^\top \mathbf{\Pi}\right) + \mathbf{\Phi}$ and $\mathbf{\Phi}$, hence (E.2) is implied by (E.1) or (38c). Also, we can verify the following relations

$$\begin{gathered} \begin{bmatrix} -I_n & \mathbf{\Pi} \end{bmatrix} \begin{bmatrix} \mathbf{\Pi} \\ I_{\beta n+q+m} \end{bmatrix} = \mathsf{O}_{n,(\beta n+q+m)}, \\ \widetilde{\Upsilon}\Upsilon = \begin{bmatrix} I_{n+\beta n} & \mathsf{O}_{(n+\beta n),(q+m)} \end{bmatrix} \Upsilon = \mathsf{O}_{(n+\beta n),(q+m)}, \\ \widetilde{\Upsilon}^\top = \begin{bmatrix} I_{n+\beta n} \\ \mathsf{O}_{(q+m),(n+\beta n)} \end{bmatrix}, \quad \widetilde{\Upsilon}_\perp := \Upsilon = \begin{bmatrix} \mathsf{O}_{(n+\beta n),(q+m)} \\ I_{q+m} \end{bmatrix}, \end{gathered} \tag{E.3}$$

where $\operatorname{rank} \Upsilon = q + m$ and $\operatorname{rank}\left(\begin{bmatrix} \mathbf{\Pi}^\top & I_{\beta n+q+m} \end{bmatrix}^\top\right) = \beta n + q + m$. Moreover, we see that $\mathsf{dim}(\operatorname{Ker} \widetilde{\Upsilon}) = q + m$ and $\mathsf{dim}(\operatorname{Ker} \begin{bmatrix} -I_n & \mathbf{\Pi} \end{bmatrix}) = \beta n + q + m$ by the rank nullity theorem [199]. With the relations in (E.3), it shows that matrices $\begin{bmatrix} \mathbf{\Pi}^\top & I_{\beta n+q+m} \end{bmatrix}^\top$, $\Upsilon$ contain the bases of the null space of $\begin{bmatrix} -I_n & \mathbf{\Pi} \end{bmatrix}$, $\widetilde{\Upsilon}$, respectively.

As a consequence, we can select $P := \begin{bmatrix} -I_n & \mathbf{\Pi} \end{bmatrix}$, $P_\perp := \begin{bmatrix} \mathbf{\Pi}^\top & I_{\beta n+q+m} \end{bmatrix}^\top$ and $Q := \widetilde{\Upsilon}$, $Q_\perp := \Upsilon$ for applying Lemma 2 to $\mathsf{Sy}\left(\mathbf{P}^\top \mathbf{\Pi}\right) + \mathbf{\Phi} \prec 0$ in (38c), which concludes that (E.1)–(E.2) are true if and only if there exists matrix $\mathbf{W} \in \mathbb{R}^{(n+\beta n) \times n}$ such that

$$\mathsf{Sy}\left(\widetilde{\Upsilon}^\top \mathbf{W} \begin{bmatrix} -I_n & \mathbf{\Pi} \end{bmatrix}\right) + \begin{bmatrix} \mathsf{O}_n & \mathbf{P} \\ * & \mathbf{\Phi} \end{bmatrix} \prec 0, \tag{E.4}$$

where the inequality is still non-convex due to the product between $\mathbf{W}$ and $\mathbf{\Pi}$. To convexify (E.4), we set

$$\mathbf{W} = \begin{bmatrix} W \\ [\alpha_i]_{i=1}^{\beta} W \end{bmatrix}, \quad W \in \mathbb{S}^n, \quad \beta = 1 + \nu + \kappa \tag{E.5}$$

in (E.4) with known scalars $\{\alpha_i\}_{i=1}^{\beta} \subset \mathbb{R}$, which gives

$$\boldsymbol{\Theta} = \mathsf{Sy}\left(\begin{bmatrix} W \\ [\alpha_i W]_{i=1}^{\beta} \\ \mathsf{O}_{(q+m),n} \end{bmatrix} \begin{bmatrix} -I_n & \boldsymbol{\Pi} \end{bmatrix}\right) + \begin{bmatrix} \mathsf{O}_n & \mathbf{P} \\ * & \boldsymbol{\Phi} \end{bmatrix} \prec 0. \tag{E.6}$$

We note that (E.6) is a sufficient condition implying (E.1) or (38c) due to the structural constraints in (E.5). Moreover, (E.6) implies that $W \in \mathbb{S}^n$ is invertible, as $-W - W^\top = -2W$ is the sole element in the first diagonal block of matrix $\boldsymbol{\Theta}$.

Let $X^\top = W^{-1}$, and apply congruent transformations to the matrix inequalities in (38a)–(38b) and (E.6). Then

$$\begin{aligned} &\left(I_{1+d} \otimes X^\top\right) \begin{bmatrix} P_1 & P_2 \\ * & P_3 \end{bmatrix} \left(I_{1+d} \otimes X\right) \succ 0, \quad X = X^\top \\ &\left(I_\nu \otimes X^\top\right) \mathbf{Q} \left(I_\nu \otimes X\right) \succ 0, \quad [*] \, \mathbf{R} \left(I_\nu \otimes X\right) \succ 0, \\ &\left[\left(I_{1+\beta} \otimes X^\top\right) \oplus I_{q+m}\right] \boldsymbol{\Theta} \left[\left(I_{1+\beta} \otimes X\right) \oplus I_{q+m}\right] \prec 0, \end{aligned} \tag{E.7}$$

hold if and only if (38a)–(38b) and (E.6) hold. Also, by the properties in (A.1) and variables

$$\begin{bmatrix} \acute{P}_1 & \acute{P}_2 \\ * & \acute{P}_3 \end{bmatrix} := [*] \begin{bmatrix} P_1 & P_2 \\ * & P_3 \end{bmatrix} \left(I_{1+d} \otimes X\right), \tag{E.8a}$$

$$\acute{\mathbf{Q}} := \underset{i=1}{\overset{\nu}{\mathsf{diag}}}\, X Q_i X, \quad \acute{\mathbf{R}} := \underset{i=1}{\overset{\nu}{\mathsf{diag}}}\, X R_i X, \tag{E.8b}$$

(E.7) can be rewritten as (42a)–(42b) and

$$\begin{aligned} \acute{\boldsymbol{\Theta}} &:= \left[\left(I_{1+\beta} \otimes X^\top\right) \oplus I_{q+m}\right] \boldsymbol{\Theta} \left[\left(I_{1+\beta} \otimes X\right) \oplus I_{q+m}\right] \\ &= \mathsf{Sy}\left(\begin{bmatrix} I_n \\ [\alpha_i I_n]_{i=1}^{\beta} \\ \mathsf{O}_{(q+m),n} \end{bmatrix} \begin{bmatrix} -X & \acute{\boldsymbol{\Pi}} \end{bmatrix}\right) + \begin{bmatrix} \mathsf{O}_n & \acute{\mathbf{P}} \\ * & \acute{\boldsymbol{\Phi}} \end{bmatrix} \prec 0, \end{aligned} \tag{E.9}$$

where $\acute{\mathbf{P}} = X\mathbf{P}\left[\left(I_\beta \otimes X\right) \oplus I_{q+m}\right] = \cdots$

$$\begin{bmatrix} \acute{P}_1 & \mathsf{O}_{n,\nu n} & \acute{P}_2 \widehat{I} & \mathsf{O}_{n,\mu n} & \mathsf{O}_{n,q} & \mathsf{O}_{n,m} \end{bmatrix}, \tag{E.10}$$

and $\acute{\boldsymbol{\Pi}} = \boldsymbol{\Pi}\left[\left(I_\beta \otimes X\right) \oplus I_{q+m}\right] = \cdots$

$$\begin{aligned} &= \begin{bmatrix} \mathbf{A}\left[\left(I_\beta \otimes X\right) \oplus I_q\right] + \mathbf{B}_1\left[\left(I_\beta \otimes KX\right) \oplus \mathsf{O}_q\right] & \mathsf{O}_{n,m} \end{bmatrix} \\ &= \begin{bmatrix} \mathbf{A}\left[\left(I_\beta \otimes X\right) \oplus I_q\right] + \mathbf{B}_1\left[\left(I_\beta \otimes V\right) \oplus \mathsf{O}_q\right] & \mathsf{O}_{n,m} \end{bmatrix} \end{aligned} \tag{E.11}$$

by defining $V = KX$ and

$$\begin{aligned} \acute{\boldsymbol{\Phi}} = \boldsymbol{\Phi}(\acute{P}_2, \acute{P}_3, \acute{\mathbf{Q}}, \acute{\mathbf{R}}, \acute{\boldsymbol{\Sigma}}) = \left[\left(I_\beta \otimes X^\top\right) \oplus I_{q+m}\right] \times \\ \boldsymbol{\Phi}(P_2, P_3, \mathbf{Q}, \mathbf{R}, \boldsymbol{\Sigma})\left[\left(I_\beta \otimes X\right) \oplus I_{q+m}\right]. \end{aligned} \tag{E.12}$$

We note that (E.9) is the same as (42c), and (E.9)–(E.12) are obtained using (E.8a)–(E.8b) and identities $\widehat{I}\left(I_\varkappa \otimes X\right) = \left(I_d \otimes X\right)\widehat{I}$ and

$$\begin{aligned} &\begin{bmatrix} \mathbf{M} \otimes I_n & \mathsf{O}_{dn,(q+m)} \end{bmatrix} \left[\left(I_\beta \otimes X\right) \oplus I_{q+m}\right] \\ &= \begin{bmatrix} I_d \mathbf{M} \otimes X I_n & \mathsf{O}_{dn,(q+m)} \end{bmatrix} \\ &= \left(I_d \otimes X\right)\begin{bmatrix} \mathbf{M} \otimes I_n & \mathsf{O}_{dn,(q+m)} \end{bmatrix}, \\ \acute{\boldsymbol{\Sigma}} &= \boldsymbol{\Sigma}\left[\left(I_\beta \otimes X\right) \oplus I_{q+m}\right] \\ &= \left(\mathbf{C} + \mathbf{B}_2\left[\left(I_\beta \otimes K\right) \oplus \mathsf{O}_q\right]\right)\left[\left(I_\beta \otimes X\right) \oplus I_{q+m}\right] \\ &= \mathbf{C}\left[\left(I_\beta \otimes X\right) \oplus I_q\right] + \mathbf{B}_2\left[\left(I_\beta \otimes V\right) \oplus \mathsf{O}_q\right] \end{aligned}$$

considering the properties in (A.1)–(A.2) with $\mathbf{M}$ in (39f). Since $-X - X^\top = -2X$ is the only element in the first diagonal block of $\acute{\boldsymbol{\Theta}}$ in (E.9) and (42c), it is clear to see that matrix $X$ is invertible if (42c) holds. This is consistent with the invertibility of $W \in \mathbb{S}^n$ implied by (E.6).

In conclusion, the inequalities in (38a)–(38b) and (42a)–(42b) are equivalent. Meanwhile, since (42c) amounts to (E.6) from which (38c) is inferred, it follows that (38a)–(38c) are satisfied if (42a)–(42c) hold for some matrix $X \in \mathbb{S}^n$ with known scalars $\{\alpha_i\}_{i=1}^{\beta} \subset \mathbb{R}$. ■

## Appendix F. A Brief Tutorial on the Application of KSD in light of the Proposed SDP Framework

By investigating numerical examples, we illustrate how to utilize the KSD concept in (21) to analyze the matrix kernels of general DDs in conjunction with the proposed SDP framework in Section 4.

### *Appendix F.1. The Case with $\boldsymbol{\phi}_1(\tau) = []_{0\times 1}$, $\mu_1 = 0$ and $\nu = 1$*

Consider the integral matrix kernel

$$\widetilde{A}(\tau) = \begin{bmatrix} 1 + \sin(\cos(10\tau)) & \tau^3 + \cos(10\tau) \\ 0 & \ln(\sin(20\tau) + 2) \end{bmatrix} \tag{F.1}$$

for $\tau \in [-r_1, 0]$. First of all, since (21) can be used without using approximations, we can simply choose

$$\boldsymbol{f}_1(\tau) = \begin{bmatrix} 1 & \tau^3 & \cos(10\tau) & \sin(\cos(10\tau)) & \ln(\sin(20\tau) + 2) \end{bmatrix}^\top \tag{F.2}$$

given that all functions in (F.1) are differentiable. Because

$$\begin{aligned} \boldsymbol{f}_1'(\tau) = \Big[ 0 \quad 3\tau^2 \quad -10\sin(10\tau) \quad \cdots \\ \cdots \quad -10\sin(10\tau)\cos(\cos(10\tau)) \quad \frac{20\cos 20\tau}{\sin(20\tau) + 2} \Big]^\top, \end{aligned} \tag{F.3a}$$

we can then select

$$\begin{aligned} \boldsymbol{\varphi}_1(\tau) = \boldsymbol{\alpha}(\tau) = \Big[ \tau^2 \quad \sin(10\tau) \quad \cdots \\ \cdots \quad \sin(10\tau)\cos(\cos(10\tau)) \quad \frac{\cos 20\tau}{\sin(20\tau) + 2} \Big]^\top, \end{aligned} \tag{F.3b}$$

which satisfies both (6c) and (6b) of the form

$$\frac{\mathsf{d}\boldsymbol{f}_1(\tau)}{\mathsf{d}\tau} = M_1 \boldsymbol{h}_1(\tau) = M_1 \begin{bmatrix} \boldsymbol{\varphi}_1(\tau) \\ \boldsymbol{f}_1(\tau) \end{bmatrix}$$

with

$$M_1 = \begin{bmatrix} \begin{smallmatrix} 0 & 0 & 0 & 0 \\ 3 & 0 & 0 & 0 \\ 0 & -10 & 0 & 0 \\ 0 & 0 & -10 & 0 \\ 0 & 0 & 0 & 20 \end{smallmatrix} & \mathsf{O}_5 \end{bmatrix}. \tag{F.4a}$$

Since all functions in (F.1) are factorized via $\boldsymbol{f}_1(\cdot)$ in (F.2), we let $\boldsymbol{\phi}_1(\tau) = []_{0\times 1}$, meaning that no functions are approximated. As a result, we have $\boldsymbol{g}_1(\tau) = \boldsymbol{h}_1(\tau)$. In view of the elements in

$\boldsymbol{g}_1(\tau)$, we then can construct

$$\widehat{A}_1 = \begin{bmatrix} \mathsf{O}_{2,8} & \begin{matrix} 1 & 0 & 0 & 1 & 0 & 1 & 1 & 0 & 0 & 0 \\ 0 & 0 & 0 & 0 & 0 & 0 & 0 & 0 & 0 & 1 \end{matrix} \end{bmatrix}$$

to satisfy $\forall \tau \in [-r_1, 0]$, $\widetilde{A}(\tau) = \widehat{A}_1(\boldsymbol{g}_1(\tau) \otimes I_n)$ for (F.1). This shows an example on the use of KSD without using approximations. Note, however, that the structure of $\boldsymbol{f}_1(\cdot)$ in (F.2) also determines the integral kernels of the KF in (C.1a). This can be reflected from the configurations of the matrices in Theorem 1 such as $\mathfrak{F}_i, \mathfrak{H}_i, M_i$ and $P_2, P_3$.

On the other hand, we can also select

$$\boldsymbol{\varphi}_1(\tau) = \begin{bmatrix} \widetilde{\boldsymbol{\alpha}}^\top(\tau) & \boldsymbol{\alpha}^\top(\tau) \end{bmatrix}^\top \tag{F.5}$$

for (6b) with $\boldsymbol{\alpha}(\tau)$ in (F.3b) and any function $\widetilde{\boldsymbol{\alpha}}(\cdot) \in \mathscr{L}^2\left([-r_1, 0]; \mathbb{R}^\lambda\right)$ such that $\boldsymbol{g}_1(\cdot)$ satisfies (6c) with $\boldsymbol{f}_1(\cdot)$ in (F.2). Accordingly, the corresponding (6b) is

$$\frac{\mathrm{d}\boldsymbol{f}_1(\tau)}{\mathrm{d}\tau} = \begin{bmatrix} \mathsf{O}_{5\times\lambda} & M_1 \end{bmatrix} \begin{bmatrix} \boldsymbol{\varphi}_1(\tau) \\ \boldsymbol{f}_1(\tau) \end{bmatrix}$$

with $M_1$ in (F.4a). By $\boldsymbol{g}_1(\tau) = \begin{bmatrix} \boldsymbol{\varphi}_1^\top(\tau) & \boldsymbol{f}_1^\top(\tau) \end{bmatrix}^\top$ with $\boldsymbol{\varphi}_1(\tau)$ in (F.5) and $\boldsymbol{f}_1(\tau)$ in (F.2), we can construct

$$\widehat{A}_1 = \begin{bmatrix} \mathsf{O}_{2,8+2\lambda} & \begin{matrix} 1 & 0 & 0 & 1 & 0 & 1 & 1 & 0 & 0 & 0 \\ 0 & 0 & 0 & 0 & 0 & 0 & 0 & 0 & 0 & 1 \end{matrix} \end{bmatrix}$$

to satisfy $\widetilde{A}(\tau) = \widehat{A}_1(\boldsymbol{g}_1(\tau) \otimes I_n)$ for the matrix-valued function in (F.1). This is an illustrative example to show that we can add an unlimited number of new $\mathscr{L}^2$ functions (linearly independent) to the $\boldsymbol{\varphi}_1(\tau)$ in $\boldsymbol{g}_1(\tau)$ by enlarging $\mathsf{dim}(\boldsymbol{\varphi}_1(\tau))$. In this case, $\boldsymbol{f}_1(\cdot)$ in (F.2) remains unchanged, so does the structure of the KF in (C.1a).

Finally, we show that an unlimited number of new WDLIFs can be added to $\boldsymbol{f}_1(\cdot)$ in (F.2) by augmenting its dimension. This usually takes place if we aim at increasing the generality of the integral kernels of KF (C.1a), thereby potentially improving the feasibility of the proposed SDP synthesis conditions. Let

$$\boldsymbol{f}_1(\tau) = \begin{bmatrix} \boldsymbol{q}^\top(\tau) & 1 & \tau^3 & \cos(10\tau) & \cdots \\ & \cdots & \sin(\cos(10\tau)) & \ln(\sin(20\tau)+2) \end{bmatrix}^\top \tag{F.6}$$

with any $\boldsymbol{q}(\cdot) \in \mathscr{W}^{1,2}\left([-r_1, 0]; \mathbb{R}^k\right)$ so that inequality $\int_{-r_1}^{0} \boldsymbol{f}_1(\tau)\boldsymbol{f}_1^\top(\tau)\mathrm{d}\tau \succ 0$ holds with $\mathsf{dim}(\boldsymbol{f}_1(\tau)) = k + 5$. Given the result in (F.3a) and the fact that $\boldsymbol{q}'(\cdot) \in \mathscr{L}^2\left([-r_1, 0]; \mathbb{R}^k\right)$, we can always find an appropriate function $\boldsymbol{p}(\cdot) \in \mathscr{L}^2\left([-r_1, 0]; \mathbb{R}^c\right)$ with constant matrix $N \in \mathbb{R}^{(k+5)\times(c+k+5)}$ such that

$$\frac{\mathrm{d}\boldsymbol{f}_1(\tau)}{\mathrm{d}\tau} = N \begin{bmatrix} \boldsymbol{p}(\tau) \\ \boldsymbol{f}_1(\tau) \end{bmatrix}, \quad \int_{-r_1}^{0} [*] \begin{bmatrix} \boldsymbol{p}^\top(\tau) & \boldsymbol{f}_1^\top(\tau) \end{bmatrix} \mathrm{d}\tau \succ 0$$

which is in line with the relations in (6b)–(6c) with the $\boldsymbol{f}_1(\cdot)$ in (F.6). By setting $\boldsymbol{\phi}_1(\tau) = [\,]_{0\times 1}$ and

$$\boldsymbol{\varphi}_1(\tau) = \boldsymbol{p}(\tau),\ M_1 = N,\ \boldsymbol{g}_1(\tau) = \begin{bmatrix} \boldsymbol{\varphi}_1^\top(\tau) & \boldsymbol{f}_1^\top(\tau) \end{bmatrix}^\top, \tag{F.7}$$

we can construct

$$\widehat{A}_1 = \begin{bmatrix} \mathsf{O}_{2,2k+2c} & \begin{matrix} 1 & 0 & 0 & 1 & 0 & 1 & 1 & 0 & 0 & 0 \\ 0 & 0 & 0 & 0 & 0 & 0 & 0 & 0 & 0 & 1 \end{matrix} \end{bmatrix}$$

to satisfy $\widetilde{A}(\tau) = \widehat{A}_1(\boldsymbol{g}_1(\tau) \otimes I_n)$ for (F.1). This example demonstrates how an unlimited number of new WDLIFs can be added to $\boldsymbol{f}_1(\cdot)$ using the KSD concept. This is a useful property that shows the feasibility of the proposed SDP framework may be increased indefinitely, even if none of the new WDLIFs in $\boldsymbol{q}(\cdot)$ are included by $\widetilde{A}(\cdot)$ in (F.1).

*Appendix F.2. The Case with $\boldsymbol{\phi}_i(\tau) \neq [\,]_{0\times 1}$, $\mu_i \neq 0$ and $\nu \geq 2$*

In the case of $\boldsymbol{\phi}_i(\tau) \neq [\,]_{0\times 1}$, $\mu_i \neq 0$ with $\nu \geq 2$, readers may refer to the numerical examples in Section 6 to see how the KSD concept can be further leveraged to address DD matrices with general structures, wherein $\boldsymbol{\phi}_i(\cdot)$ is approximated by $\boldsymbol{h}_i(\cdot)$. Similar to the analysis in Section Appendix F.1, an unlimited number of WDLIFs and $\mathscr{L}^2$ functions can be added to $\boldsymbol{f}_i(\cdot)$ and $\boldsymbol{\varphi}_i(\cdot)$, respectively, as long as the matrix inequality in (6c) is satisfied. Additionally, it is possible to move one of the functions in $\boldsymbol{\varphi}_i(\cdot)$ to $\boldsymbol{\phi}_i(\cdot)$, or vice versa. In summary, users have substantial freedom to decide which DD kernels are approximated, and which of those are directly factorized. Meanwhile, the feasibility of our SDP synthesis conditions may be improved indefinitely by adding an unlimited number of WDLIFs to $\boldsymbol{f}_i(\cdot)$ to increase the generality of KF (C.1a).